\documentclass[11pt]{article}

\usepackage{graphicx}
\usepackage{latexsym,amsmath,amsfonts,amscd, amsthm, dsfont}
\usepackage{bm,color}
\usepackage{epsfig,verbatim,epstopdf,graphics}
\usepackage{subfigure}
\usepackage{changebar}
\usepackage{scalerel}
\usepackage{multirow}
\allowdisplaybreaks

\newtheoremstyle{plainNoItalics}{}{}{\normalfont}{}{\bfseries}{.}{ }{}

\theoremstyle{plain}
\newtheorem{thm}{Theorem}[section]

\theoremstyle{plainNoItalics}

\newtheorem{rem}[thm]{Remark}
\newtheorem{prop}[thm]{Proposition}

\newcommand{\bx}{{\bf x}}
\newcommand{\bv}{{\bf v}}

\newcommand{\bE}{{\bf E}}

\newcommand{\be}{\begin{eqnarray}}
\newcommand{\ee}{\end{eqnarray}}
\newcommand{\beno}{\begin{eqnarray*}}
\newcommand{\eeno}{\end{eqnarray*}}

\makeatletter

\newcommand{\Rmnum}[1]{\expandafter\@slowromancap\romannumeral #1@}
\makeatother

\begin{document}

\baselineskip=1.8pc


\begin{center}
{\bf
A WENO-based Method of Line Transpose Approach for Vlasov Simulations
}
\end{center}

\vspace{.2in}
\centerline{
Andrew Christlieb\footnote{
 Department of Mathematics and Electrical Engineering, Michigan State University, East Lansing, MI, 48824. E-mail: andrewch@math.msu.edu
},
 Wei Guo \footnote{
Department of Mathematics, Michigan State University, East Lansing, MI, 48824. E-mail:
wguo@math.msu.edu},
Yan Jiang\footnote{Department of Mathematics, Michigan State University, East Lansing, MI, 48824. E-mail: jiangyan@math.msu.edu}
}

\bigskip
\noindent
{\bf Abstract.}

 In this paper, a high order implicit Method of Line Transpose (MOL$^T$) method based on a weighted essentially non-oscillatory (WENO) methodology is developed for one-dimensional linear transport equations and further applied to the Vlasov-Poisson (VP) simulations via dimensional splitting. In the MOL$^T$ framework, the time variable is first discretized by a diagonally implicit strong-stability-preserving Runge-Kutta method, resulting in a boundary value problem (BVP) at the discrete time levels. Then an integral formulation coupled with a high order WENO methodology is employed to solve the BVP. As a result, the proposed scheme is high order accurate in both space and time and free of oscillations even though the solution is discontinuous or has sharp gradients. Moreover, the scheme is able to take larger time step evolution compared with an explicit MOL WENO scheme with the same order of accuracy. The desired positivity-preserving (PP) property of the scheme is further attained by incorporating a newly proposed high order PP limiter. We perform numerical experiments on several benchmarks including linear advection, solid body rotation problem; and on the Landau damping, two-stream instabilities, bump-on-tail, and plasma sheath by solving the VP system. The efficacy and efficiency of the proposed scheme is numerically verified.

\vfill

\newpage

\section{Introduction}
Transport modeling plays an important role in a wide range of applications, including fluid dynamics, plasma physics, atmospheric sciences, among many others. The main focus of this paper is to develop of efficient high order transport solvers for the Vlasov simulations. As a high-dimensional kinetic transport equation, the Vlasov-Poisson (VP) system is considered a fundamental model in plasma physics which describes the dynamics of charged particles due to the self-consistent electric force:
\begin{subequations}
	\label{eq:vp}
	\begin{align}
	&f_t + \mathbf{v}\cdot\nabla_\bx f + \bE(\bx,t)\cdot\nabla_\bv f =0, \quad \bx\times \bv\in\Omega_\bx\times\Omega_\bv\label{eq:vlasov}\\
	& \bE(\bx,t)=-\nabla_\bx\phi(\bx,t), \ \ -\Delta_\bx\phi(\bx,t)=-1+\rho(\bx,t) \label{eq:possion}
	\end{align}
\end{subequations}
where $f(\bx,\bv,t)$ describes the probability of finding a particle with velocity $\bv$ at position $\bx$ at time $t$, $\bE$ is the electrostatic field, $\phi$ is the self-consistent electrostatic potential, and $\rho(\bx,t)=\int_{\Omega_\bv} f(\bx,\bv,t)d\bv$ is the electron charge density and the 1 represents the uniformly distributed infinitely massive ions in the background. 

The VP simulation is challenging. Besides the famous curse of dimensionality, the underlying system may develop filamentation solution structures, since the long-range forces drive
 the solution out of thermo-equilibrium. Hence, it is of importance that the used numerical scheme is able to effectively capture filamentation structures without producing spurious oscillations. Further, the VP system conserves several physical quantities including the total mass, momentum, and energy, which requires that the numerical scheme converses those invariants on the discrete level. In addition, as a probability density function in phase space, the unknown function $f$ always remains non-negative. In practice, the positivity preserving property is another highly desired property for the design of Vlasov solver, since negative solutions may trigger unphysical phenomenon in plasmas.
  
 There exist a large amount of successful VP solvers in the literature, such as the prominent particle-in-cell
 (PIC) method \cite{birdsall2005plasma,hockney2010computer,verboncoeur2005particle} , which evolves the solution by following the trajectories of some sampled macro-particles, Eulerian type approaches \cite{FilbetS,zhou2001numerical,cheng2012study,cheng2014ts}, for which the state variable is advanced according to the underlying equations on a fixed numerical grid, and semi-Lagrangian approach \cite{cheng1975integration,parker1997convected,begue1999two,sonnendrucker1999semi,filbet2001conservative,besse2005semi,crouseilles2009forward,crouseilles2010conservative,rossmanith2011positivity,qiu2011positivity,christlieb2014high,guo2013hybrid}, which updates the solution by tracking the characteristics on a fixed mesh grid. However, these popular schemes still have limitations. For example, the PIC method is known to suffer from low order sampling noise; Eulerian approaches are subject to stringent time step restriction associated with explicit time stepping; general boundary conditions are difficult to handle in a semi-Lagrangian setting. Also note that many of these existing methods are in the Method of Line (MOL) framework, meaning that the spatial variable is first discretized, then the numerical solution is updated in time by coupling a suitable time integrator. 
%
In this paper, we consider an alternative approach to advance the solution, which is called the Method of Line Transpose (MOL$^T$) or Rothes's method in the literature \cite{salazar2000theoretical,schemann1998adaptive}. As its name implies, such an approach is defined in an orthogonal fashion of the MOL schemes, i.e., the discretization is first carried out for temporal variable, resulting in a boundary value problem (BVP) at the discrete time levels. Then a preferred BVP solver can be further applied to advance the solution. For example, in \cite{causley}, an implicit integral formulation based on the MOL$^T$ framework for one-dimensional wave equations is developed and analyzed. It is proven that the scheme is A-stable, hence allowing for extra large time step evolution. In \cite{causley2013method}, the computational cost of the scheme is reduced from $\mathcal{O}(N^2)$ to $\mathcal{O}(N)$, where $N$ is the number of discrete mesh points, by taking advantage of an important property of the analytical solution of the wave equations. More recently, a successive convolution  technique \cite{causley2014higher} is developed to enhance the temporal accuracy to arbitrary high order.  
A notable advantage of such a scheme is that, even though it is implicit in time, we do not need to explicitly solve a linear system, while the BVP is inverted analytically in an integral formulation. 
The method is further extended to multi-dimensional wave equations  by utilizing the alternating direction implicit strategy at the expense of introducing splitting errors \cite{causley}. However, the simplicity of the scheme is preserved. The recent applications of this MOL$^T$ approach include the heat equation, the Allen-Cahn equation \cite{Causley2015method}, and Maxwell's equations \cite{Cheng2015asymptotic}.

In this paper, we consider formulating a numerical scheme suited for the Vlasov equation in the MOL$^T$ framework in order to take advantage of its good properties such as the simplicity and the ability to take large time steps for evolution. In order to accommodate the MOT$^T$ approach, the nonlinear multi-dimensional Vlasov equation is first split into several lower-dimensional linear transport equations. Such an idea has been widely used in semi-Lagrangian Vlasov simulations \cite{cheng1975integration}.
As we mentioned earlier, the extra difficulty of VP simulations lies in the distinctive filamentation solution structures, which requires that the employed numerical scheme is able to effectively capture sharp gradient of the solution without producing spurious oscillations even on an under-resolved mesh. Note that we cannot directly apply the existing MOL$^T$ schemes to the VP system, since they are all linear schemes and hence exhibit oscillations when the solution is discontinuous or has a sharp gradient. In order to control undesired oscillations, we propose to further couple two numerical ingredients. Firstly, in the MOL$^T$ framework, the time variable is discretized via an implicit high order strong-stability-preserving (SSP) Runge-Kutta (RK) method \cite{gottlieb2001strong}. Secondly, a robust weighted essential non-oscillatory (WENO) methodology is incorporated when inverting the BVP. We expect that the resulting scheme is efficient, high order accurate, robust, and  able to take large time step evolution (larger than the explicit counterpart). Furthermore, a high order positivity-preserving (PP) limiter suited for the proposed MOL$^T$ scheme is developed  to ensure such a high desired property for the numerical solution.

We also would like to point out that the proposed MOL$^T$ WENO scheme shares some similarity with the semi-Lagrangian WENO scheme developed in \cite{Qiu_Christlieb,Qiu_Shu2}, in the sense that both schemes are high order accurate and able to 
effectively capture the filamentation structures with large time step. On the other hand, compared with semi-Lagrangian WENO schemes, the proposed method is able to conveniently handle general boundary conditions, which is not trivial for semi-Lagrangian approaches.

The rest of the paper is organized as follows. In Section 2, we formulate the numerical scheme for one-dimensional linear transport equations by coupling several main numerical ingredients such as the MOL$^T$ approach, the WENO methodology and the high order PP limiter. In Section 3, the method is extended to two-dimensional transport equations via dimensional splitting. We present several benchmark tests in Section 4 to verify the performance of the proposed scheme, including linear transport, rigid body rotation, as well as Landau damping, two-stream instabilities, bump-on-tail and plasma sheath for VP simulations.
We conclude the paper in Section 5.

\section{Formulation of the Scheme}

We start with the following one-dimensional advection equation
\begin{equation}\label{eq:1dadv}
  u_t + cu_x = 0,\quad x\in[a,b],
\end{equation}
where $c$ is constant and denotes the wave propagation speed. The boundary condition is assumed to be periodic,  Dirichlet boundary condition
$$u(a,t)=g_{1}(t),\ \text{for} \ c>0, \ \ \ \text{or} \ \ \ \ u(b,t)=g_{2}(t),\ \text{for} \ c<0, $$
or  Neumann boundary condition
$$u_{x}(a,t)=h_{1}(t),\ \text{for} \ c>0, \ \ \ \text{or} \ \ \ \ u_{x}(b,t)=h_{2}(t),\ \text{for} \ c<0. $$

\subsection{The method of line transpose approach}

As discussed in the introduction, in a MOL$^T$ framework, we first discretize the time variable while leaving the spatial variable continuous. For example, we can apply the first order backward Euler method to
\eqref{eq:1dadv}, and obtain
$$
\frac{u^{n+1}-u^n}{\Delta t} + cu^{n+1}_x =0,
$$
i.e.,
\begin{equation}\label{eq:mol1}
  u^{n+1}_x +\text{sgn}(c)\cdot\alpha u^{n+1} =\text{sgn}(c)\cdot\alpha u^n,
\end{equation}
where $\Delta t$ denotes the time step and $\alpha=\frac{1}{|c\Delta t|}$. Note that \eqref{eq:mol1} is an initial value problem (IVP), for which we can obtain an explicit representation for solution $u^{n+1}$ at time level $t^{n+1}$. If we assume $c>0$, $u^{n+1}$ can be represented by
\begin{equation}\label{eq:representationL}
  u^{n+1}(x) = I^L[u^n,\alpha](x) + A^{n+1}e^{-\alpha(x-a)},
\end{equation}
where
\begin{equation}\label{eq:convolutionL}
  I^L[v,\alpha](x) \doteq \alpha\int_a^x e^{-\alpha(x-y)}v(y)dy,
\end{equation}
and $A^{n+1}$ is a constant and should be determined by the boundary condition. The superscript $L$ indicates that the characteristics traverse from the left to the right.
If $c<0$, then the characteristics traverse from the right to the left. Due to the symmetry, we can advance the solution as
\begin{equation}\label{eq:representationR}
  u^{n+1}(x) = I^R[u^n,\alpha](x) + B^{n+1}e^{-\alpha(b-x)},
\end{equation}
when the first order backward Euler scheme is used. Here,
\begin{equation}\label{eq:convolutionR}
  I^R[v,\alpha](x) \doteq \alpha\int_x^b e^{-\alpha(y-x)}v(y)dy,
\end{equation}
 and again the constant $B^{n+1}$ is determined by the boundary condition.

We need to further discretize the spatial variable based on the semi-discrete schemes \eqref{eq:representationL} and \eqref{eq:representationR}. Following the idea in \cite{causley,causley2013method}, the domain $[a,b]$ is discretized by $M+1$ uniformly distributed grid points
$$a=x_0<x_1<\cdots<x_{M-1}<x_M=b,$$
and $\Delta x=\frac{b-a}{M}$ denotes the mesh size. The solution is accordingly discretized at points $\{x_i\}_{i=0}^M$, and we use  $u^n_i$ to denote the numerical solution at location $x_i$ at time level $t^n$. The main part of the algorithm is to evaluate the convolution integral $I^{L}[v,\alpha](x)$ in \eqref{eq:representationL}, \eqref{eq:representationR} or $I^{R}[v,\alpha]$ in \eqref{eq:convolutionL}, \eqref{eq:convolutionR} based on the point values $\{v_i\}_{i=0}^M$.

We first consider the case of $c>0$. Note that,
\begin{equation}\label{eq:recursive}
I^L[v,\alpha](x) = I^L[v,\alpha](x-\delta_x)e^{-\alpha\delta_x} + J^L[v,\alpha,\delta_x](x),
\end{equation}
where
\begin{equation}\label{eq:JL}
J^L[v,\alpha,\delta_x](x) \doteq  \alpha\int_{x-\delta_x}^{x} v(y)e^{-\alpha (x-y)}dy,
\end{equation}
and $\delta_x$ is a positive constant.
Therefore, in order to obtain $I^L[v,\alpha](x_i)$ denoted by $I^L_i$,
we only need to approximate the increment $J^L[v,\alpha,\Delta x](x_i)$ denoted by $J^L_i$, and obtain $I^L_i$ via a recursive relation, i.e.,  
\begin{equation}\label{eq:recursive1}
I^L_i = I^L_{i-1}e^{-\alpha\Delta x} + J^L_i,\quad i=1,\cdots,M, \quad I^L_0 = 0.
\end{equation}
For the case of $c<0$, again thanks to the symmetry, we have the following identity
\begin{equation}\label{eq:recursive}
I^R[v,\alpha](x) = I^R[v,\alpha](x+\delta_x)e^{-\alpha\delta_x} + J^R[v,\alpha,\delta_x](x),
\end{equation}
where
\begin{equation}\label{eq:JR}
J^R[v,\alpha,\delta_x](x) \doteq  \alpha\int_{x}^{x+\delta_x} v(y)e^{-\alpha (y-x)}dy,
\end{equation}
and hence a similar  recursive relation can be obtained
\begin{equation}\label{eq:recursive2}
I^R_i = I^R_{i+1}e^{-\alpha\Delta x} + J^R_i,\quad i=0,\cdots,M-1, \quad I^R_M = 0,
\end{equation}
where $I^R_i$ and $J^R_i$ denote the approximate values of $I^R[v,\alpha](x_i)$ and $J^R[v,\alpha,\Delta x](x_i)$, respectively.

After computing the contributions from the left or right characteristic, we further need to find the global coefficient $A$ in \eqref{eq:representationL} or $B$ in \eqref{eq:representationR}, which is used to enforce the boundary condition.
Since $I^{L}[u^{n},\alpha](a)=I^{R}[u^{n},\alpha](b)=0$, for a Dirichlet boundary condition, the global coefficient can be obtained by
\begin{equation}
\label{eq:bcdir}
A^{n+1}=g_{1}(t^{n+1}), \ \text{or} \ \  B^{n+1}=g_{2}(t^{n+1}).
\end{equation}
Similarly, for a Neumann boundary condition, because $(I^{L}[v,\alpha])'(x)=-\alpha I^{L}[v,\alpha](x)+\alpha v(x)$ and $(I^{R}[v,\alpha])'(x)=\alpha I^{R}[v,\alpha](x)-\alpha v(x)$, we can get the coefficients by
\begin{align}
\label{eq:bcneu}
A^{n+1}=v(a)-\frac{1}{\alpha}h_{1}(t^{n+1}), \ \  \text{or} \ \ 
B^{n+1}=v(b)+\frac{1}{\alpha}h_{2}(t^{n+1}).
\end{align}
For the periodic boundary condition, the coefficient can be determined by using the  mass conservation property of the solution, i.e,
\begin{equation*}
\int_{a}^{b}u(x,t)dx=\int_{a}^{b}u(x,0)dx.
\end{equation*} 
On the discrete level, we enforce
\begin{equation*}
\sum_{i=0}^{M-1}u^{n}_{i}=\sum_{i=0}^{M-1}u^{n+1}_{i}
=\sum_{i=0}^{M-1}I^{L}_{i}+A^{n+1}\sum_{i=0}^{M-1}e^{-i\alpha\Delta x}, \ \ \
u^{n+1}_{M}=u^{n+1}_{0}, \ \ \ for \ c>0,
\end{equation*}
 or 
\begin{equation*}
\sum_{i=1}^{M}u^{n}_{i}=\sum_{i=1}^{M}u^{n+1}_{i}
=\sum_{i=1}^{M}I^{R}_{i}+B^{n+1}\sum_{i=1}^{M}e^{-(M-i)\alpha\Delta x} , \ \ \
u^{n+1}_{0}=u^{n+1}_{M}, \ \ \ for \ c<0.
\end{equation*}
Then the coefficient is obtained
\begin{equation}
\label{eq:bcL}
A^{n+1} =\left(\sum_{i=0}^{M-1}u^{n}_{i}-\sum_{i=0}^{M-1}I^{L}_{i}\right)\bigg/\sum_{i=0}^{M-1}e^{-i\alpha\Delta x},
\end{equation}
or
\begin{equation}
\label{eq:bcR}
B^{n+1}=\left(\sum_{i=1}^{M}u^{n}_{i}- \sum_{i=1}^{M}I^{R}_{i}\right)\bigg/ \sum_{i=1}^{M}e^{-(M-i)\alpha\Delta x}.
\end{equation}

\subsection{WENO methodology}

The existing MOL$^T$ schemes mainly use linear interpolation methods to compute increment $J^L_i$ and $J^R_i$ such as \cite{causley}, which work well for smooth problems. However, the schemes may generate spurious oscillations when the solution has discontinuities because of the famous Gibbs phenomenon. In order to control such undesired oscillations, we will incorporate a WENO methodology in the MOL$^T$ framework to evaluate  $J^L_i$ and $J^R_i$ below.

The first WENO scheme was introduced by Liu {\em et al.} \cite{liu1994weighted}, as a third order accurate finite volume scheme to solve hyperbolic conservation laws in one dimension. Later, in \cite{jiang1996efficient}, a general framework for arbitrary order accurate finite difference WENO schemes was designed, which is more efficient for multi-dimensional calculations. The schems have the advantage of attaining high order accuracy in smooth regions of the solution while maintaining sharp and essentially non-oscillatory transitions of discontinuities. More information about WENO can be found in \cite{shu2009high}.

In our work, a WENO integration procedure is used, based on the idea proposed in \cite{chou2006high,chou2007high,liu2009positivity}. 
Below, a $(2k-1)$-th order approximation of $J^{L}_{i}$ is provided as an example, where $k$ is an integer. The corresponding coefficients of the third order and  fifth order schemes are given in the Appendix. The whole procedure consists of four steps.

\begin{enumerate}
\item
On each small stencil $S_{r}(i)=\{x_{i-r-1},\ldots,x_{i-r-1+k}\}$, $r=0,\ldots,k-1$, which contains $x_{i-1}$ and $x_{i}$, there is a unique polynomial $p_{r}(x)$ of degree at most $k$ which interpolates $v(x)$ at the nodes in $S_{r}(i)$. Then we can get the integral 
\begin{equation}
\label{eq:weno1}
J^{L}_{i,r} =\alpha \int_{x_{i-1}}^{x_{i}}e^{-\alpha(x_{i}-y)}p_{r}(y)dx
=\sum_{j=0}^{k}c^{(r)}_{i-r-1+j}v_{i-r-1+j},
\end{equation}
where the coefficients $c^{(r)}_{i-r-1+j}$ depend on $\alpha$ and the cell size $\Delta x$, but not on $v$.

\item
Similarly, on the big stencil $S(i)=\{x_{i-k},\ldots,x_{i+k-1}\}$, there is a unique polynomial $p(x)$ of degree at most $2k-1$ interpolating $v(x)$ at the nodes in $S(i)$. Then we have
\begin{equation}
\label{eq:weno2}
J^{L}_{i} =\alpha \int_{x_{i-1}}^{x_{i}}e^{-\alpha(x_{i}-y)}p(y)dx
=\sum_{j=0}^{2k-1}c_{i-k+j}v_{i-k+j},
\end{equation}
and linear weights $d_{r}$, such that
\begin{equation}
J^{L}_{i}=\sum_{r=0}^{k-1}d_{r}J^{L}_{i,r},
\end{equation}
where $\sum_{r=0}^{k-1}d_{r}=1$.

\item
Change the linear weights $d_{r}$ into nonlinear weights $\omega_{r}$, using
\begin{equation}
\omega_{r}=\frac{\tilde{\omega}_{r}}{\sum_{s=0}^{k-1}\tilde{\omega}_{s}}, \ \ r=0,\ldots,k-1,
\end{equation}
with
\begin{equation*}
\tilde{\omega}_{r}=\frac{d_{r}}{(\epsilon+\beta_{r})^2}.
\end{equation*}
Here, $\epsilon>0$ is a small number to avoid the denominator becoming zero, and we take $\epsilon=10^{-6}$ in our numerical tests. The smoothness indicator $\beta_{r}$, measuring the relative smoothness of the function $v(x)$ in the stencil $S_{r}(i)$, is defined by
\begin{equation}
\beta_{r}=\sum_{l=1}^{k} \int_{x_{i-1}}^{x_{i}}\Delta x^{2l-1} (\frac{\partial^{l}p_{r}(x)}{\partial x^{l}})^2 dx,
\end{equation}

\item
Lastly, we have 
\begin{equation}
J^{L}_{i}=\sum_{r=0}^{k-1}\omega_{r}J^{L}_{i,r}.
\end{equation}

\end{enumerate}

The process to obtain $J^{R}_{i}$ is mirror-symmetric to that for $J^{L}_{i}$, with respect to the point $x_{i}$.

When dealing with a Dirichlet or a Neumann boundary condition, we further need to extrapolate point values at the ghost points $\{x_{-(k-1)},\ldots,x_{-1}\}$ and $\{x_{M+1},\ldots,x_{M+k-1}\}$ using the point values in the domain, where
$$x_{-i}=x_{0}-i\Delta x, \ \ x_{M+i}=x_{M}+i\Delta x.$$  
Following the idea in \cite{tan2012efficient}, a WENO extrapolation strategy is designed to get more robust results. For instance, when approximating the left ghost point $x_{-i}$, we consider the big stencil $S_{2k-2}=\{x_{0},\ldots,x_{2k-2}\}$:

\begin{enumerate}
\item On each small stencil $S_{r}=\{x_{0},\ldots,x_{r}\}$, $r=0,\ldots,2k-2$, we can construct a polynomial $p_{r}(x)$ of degree at most $r$ interpolating $v(x)$ at the nodes in $S_{r}$. And $v_{-i}$ can be extrapolated by
\begin{align}
v_{-i}=\sum_{r=0}^{2k-2}d_{r}p_{r}(x_{-i}),
\end{align}
where $d_{0}=\Delta x^{2k-2}$, $d_{1}=\Delta x^{2k-3}$, $\ldots$, $d_{2k-3}=\Delta x$, and $d_{2k-2}=1-\sum_{r=0}^{2k-3}d_{r}$.

\item Change the weights $d_{r}$ into nonlinear weights $\omega_{r}$
\begin{equation}
\omega_{r}=\frac{\tilde{\omega}_{r}}{\sum_{s=0}^{2k-2}\tilde{\omega}_{s}}, \ \ r=0,\ldots,2k-2,
\end{equation}
with
\begin{equation*}
\tilde{\omega}_{r}=\frac{d_{r}}{(\epsilon+\beta_{r})^2}.
\end{equation*}
Again, we take $\epsilon=10^{-6}$ in our numerical tests to avoid the denominator becoming zero. And the smoothness indicator $\beta_{r}$ is defined by
\begin{align*}
& \beta_{0}=\Delta x^2,\\
& \beta_{r}=\sum_{l=1}^{r} \int_{x_{-1}}^{x_{0}}\Delta x^{2l-1} (\frac{\partial^{l}p_{r}(x)}{\partial x^{l}})^2 dx, \ \ r=1,\ldots,2k-2.
\end{align*}

\item Finally, we have
\begin{align*}
v_{-i}=\sum_{r=0}^{2k-2}\omega_{r}p_{r}(x_{-i}).
\end{align*}

\end{enumerate} 

Similar, we can get the WENO extrapolation on right ghost points $x_{M+i}$ based on $S_{2k-2}=\{x_{M-2k+2},\ldots,x_{M}\}$.

\subsection{High order time discretization}
Note that scheme \eqref{eq:mol1} is first order accurate in time. Below, we will extend our schemes to high order accuracy via the implicit SSP RK method \cite{gottlieb2001strong,ketcheson2009optimal}. 

To approximate the solutions of ordinary differential equation
\begin{equation}
u_{t}=F(u),
\end{equation}
an $s$-stage and $k$-th order RK method, denoted as $RK(s,k)$, is defined as
\begin{subequations}
\label{eq:rk}
\begin{align}
& u^{(i)}=u^{n}+\Delta t\sum_{j=1}^{s}a_{ij}F(u^{(j)},t_{n}+c_{j}\Delta t),   \ \ \ i\leq j\leq s,\\
& u^{n+1}=u^{n}+\Delta t\sum_{j=1}^{s}b_{j}F(u^{(j)},t_{n}+c_{j}\Delta t), 
\end{align}
\end{subequations}
with the standard assumption
\begin{equation*}
c_{i}=\sum_{j=1}^{s}a_{ij}.
\end{equation*}
Denote the matrix $\textbf{A}=(a_{ij})_{s\times s}$ and the vector $\textbf{b}=(b_{i})_{s\times 1}$.
In our scheme, we will use the diagonally implicit (DIRK) methods, for which $\textbf{A}$ is lower triangular.  Denote $\textbf{1}=[1,\ldots,1]^{T}\in R^{s}$, $\textbf{U}^{rk}=[u^{(1)},\ldots, u^{(s)}] ^{T}$, $\textbf{F}=[\Delta t F(u^{(1)},t_{n}+c_{1}\Delta t),\ldots,\Delta tF(u^{(s)},t_{n}+c_{s}\Delta t)]^T$, and $\Lambda=diag(a_{11},\ldots,a_{ss})$. Then, we can get
$\textbf{F}=\textbf{A}^{-1}(\textbf{U}^{rk}-u^{n}\cdot\textbf{1})$. Hence
\begin{subequations}
\begin{align}
& \textbf{U}^{rk}-\Lambda\cdot\textbf{F}
=u^{n}\cdot\Lambda\cdot\textbf{A}^{-1}\cdot\textbf{1}+ (\textbf{A}-\Lambda)\cdot\textbf{A}^{-1} \textbf{U}^{rk}
\label{eq:rk1}\\
& u^{n+1}=u^{n}+\textbf{b}^{T}\cdot\textbf{F}
=u^{n}+\textbf{b}^{T}\cdot\textbf{A}^{-1}(\textbf{U}^{rk}-u^{n}\cdot\textbf{1})
\label{eq:rk2}
\end{align}
\end{subequations}
Note that $\textbf{A}$ is lower triangular. So $(\textbf{A}-\Lambda)\cdot\textbf{A}^{-1}$ is also lower triangular, with all diagnal elements equal to zero.  

For problem \eqref{eq:1dadv}, $F(u)=-cu_{x}$. Corresponding, each stage of \eqref{eq:rk1} can be rewritten as
\begin{align}
\label{eq:rkmolt}
u^{(i)}_{x} +\text{sgn}(c\cdot a_{ii}) \alpha_{i} u^{(i)} =\tilde{\alpha}_{i}u^{n}+ \sum_{j=1}^{i-1}\tilde{\alpha}_{ij}u^{(j)},
\end{align}
where $\alpha_{i} = 1/|c\cdot a_{ii}\Delta t|$. We can solve \eqref{eq:rkmolt} with MOL$^T$ method mentioned above. And finally, $u^{n+1}$ is obtained by \eqref{eq:rk2}.

When problem \eqref{eq:1dadv} combining with the periodic boundary condition, the global coefficients $A^{rk,i}$ or $B^{rk,i}$ on each stage can be obtained by \eqref{eq:bcL} or \eqref{eq:bcR}. For a Dirichlet boundary condition, the scheme may lose accuracy if we define $A^{rk,i}=g_{1}(t_{rk}^{i})$ or $B^{rk,i}=g_{2}(t_{rk}^{i})$. To avoid order reduction, we follow the idea in \cite{Isaia_rkbc}. In particular, denote $\textbf{c}^{l}=[c_{1}^{l},\ldots,c_{s}^{l}]^{T}$ for $l>0$. To achieve $k$-th order of accuracy, boundary values at intermediate stages are defined as
\begin{align}
\textbf{A}^{rk}
=& \sum_{l=0}^{k-2}(\text{A}^{l} \textbf{1} \otimes \Delta t^{l} \textbf{c}^{l})\frac{d^{l}}{dt^{l}} g_{1}(t_{n})
+(\textbf{A}^{k-1}\otimes \Delta t^{k-1}\textbf{c}^{k-1})\frac{d^{k-1}}{dt^{k-1}} \textbf{G}_{1}^{n+c}
\end{align}
for $c>0$, where
\begin{align*}
& \textbf{A}^{rk}=[A^{rk,1},\ldots,A^{rk,s}]^T, \\
& \textbf{G}_{1}^{n+c}=[g_{1}(t^{n}+c_{1}\Delta t),\ldots,g_{1}(t^{n}+c_{s}\Delta t)]^T,
\end{align*}
or
\begin{align}
\textbf{B}^{rk}
=& \sum_{l=0}^{k-2}(\text{A}^{l} \textbf{1} \otimes \Delta t^{l} \textbf{c}^{l})\frac{d^{l}}{dt^{l}} g_{2}(t_{n})
+(\textbf{A}^{k-1}\otimes \Delta t^{k-1}\textbf{c}^{k-1})\frac{d^{k-1}}{dt^{k-1}} \textbf{G}_{2}^{n+c}
\end{align}
for $c<0$, where
\begin{align*}
& \textbf{B}^{rk}=[B^{rk,1},\ldots,B^{rk,s}]^T, \\
& \textbf{G}_{2}^{n+c}=[g_{2}(t^{n}+c_{1}\Delta t),\ldots,g_{2}(t^{n}+c_{s}\Delta t)]^T.
\end{align*}

The similar idea is used when dealing with a Neumann boundary condition. 

In our numerical tests, we will use $RK(2,3)$ and $RK(4,4)$ schemes to obtain third order accuracy and fourth order accuracy in time, respectively. Corresponding coefficients are listed as appendix.

\subsection{Positivity-preserving limiter}

The exact solution of \eqref{eq:1dadv} satisfies the positivity preserving principle, meaning that if the initial condition $u(x,0)\geq0$, then the exact solution $u(x,t)\geq0$ for any time $t>0$. In practice, it is highly desirable to have numerical schemes also satisfy such a property on the discrete level. It is even more crucial  for the Vlasov simulation, since negative values may trigger some nonphysical precess in plasmas \cite{Qiu20118386}. Recently, Zhang {\em el at.} developed a high order positivity-preserving limiters for finite volume WENO methods and finite element DG methods, which is able to maintain high order accuracy of the original scheme and easy to implement \cite{zhang2010maximum,zhang2010positivity,zhang2011maximum,Zhang_JSC_2010}. Later, Xu {\em et al.} \cite{xu2014parametrized,liang2014parametrized,xiong2013parametrized,xiong2014parametrized,xiong2014high} developed a parameter maximum principle preserving flux limiter for finite difference WENO schemes by limiting the high order flux towards a first order monotone flux. However, both limiters are exclusively designed for MOL type methods and hence cannot be directly applied in the MOL$^T$ setting. Below, we will formulate a new positivity-preserving (PP) limiter suited for the proposed WENO MOL$^T$ scheme, and also show that such a limiter is able to preserve positive while maintaining the original high order accuracy.

Assume that we have $u^{n}_{i}\geq0$ for all $i$ at time $t^n$. After obtaining the numerical solutions $u^{n+1}_{i}$, we can rewrite our scheme as a conservative form  
$$u^{n+1}_{i}=u^{n}_{i}-\Delta t/\Delta x (\hat{f}_{i+1/2}-\hat{f}_{i-1/2}),$$
in which the artificial numerical fluxes $\hat{f}_{i+1/2}$ are recursively defined by
\begin{align*}
\hat{f}_{i+1/2}=\left\{\begin{array}{ll}
0, & i = -1,\\
\hat{f}_{i-1/2}-\Delta x/\Delta t(u^{n+1}_{i}-u^{n}_{i}), & i=0,\ldots,M.\\
\end{array}
\right.
\end{align*}

If we are able to modify the numerical flux $\hat{f}_{i+1/2}$ into $\tilde{f}_{i+1/2}$, such that
\begin{align}
\label{eq:pplimiter}
u^{new}_{i}=u^{n}_{i}-\Delta t/\Delta x (\tilde{f}_{i+1/2}-\tilde{f}_{i-1/2})\geq0, \ \ i=0,\ldots,M,
\end{align}
then the modified solution $u^{new}$ at time level $t=t^{n+1}$ will be non-negative automatically. Roughly speaking, the main idea of the limiter is to enforce the negative solution value back to zero by decreasing the outflow flux. For instance, when $c>0$, we will reset the numerical flux $\hat{f}_{i+1/2}$ if $u^{n+1}_{i}<0$. Details of the limiter are shown below. 

Firstly, we consider the case $c>0$ and a Dirichlet boundary condition. Numerical flux $\hat{f}_{i+1/2}$  is modified into $\tilde{f}_{i+1/2}$, starting from $i=1$ to $i=M$,
\begin{align}
\label{eq:limit}
\tilde{f}_{i+1/2}=\left\{\begin{array}{ll}
\tilde{f}_{i-1/2}+ \Delta x/\Delta t u^{n}_{i}, & \text{if}  \ u^{n+1,1}_{i}<0, \\
\hat{f}_{i+1/2}, & \text{otherwise,}\\ 
\end{array}
\right.
\end{align}
where
$$u_{i}^{n+1,1}=u^{n}_{i}-\Delta t/\Delta x (\hat{f}_{i+1/2}-\tilde{f}_{i-1/2}).$$
$\tilde{f}_{-1/2}=\hat{f}_{-1/2}$ and $\tilde{f}_{1/2}=\hat{f}_{1/2}$, due to the Dirichlet boundary condition. For a Nuemann boundary condition, we can also use \eqref{eq:limit} to modify fluxes. Once the modified fluxes are obtained, the positivity-preserving solution is computed by \eqref{eq:pplimiter}.

For periodic boundary conditions, however, we need an additional step in order to preserve the positivity and periodicity of the solution. In particular, we first modify  the numerical flux $\hat{f}_{i+1/2}$ into $\hat{\hat{f}}_{i+1/2}$ by \eqref{eq:limit} i.e., starting from $i=0$ into $i=M-1$, with
\begin{align*}
\hat{\hat{f}}_{i+1/2}=\left\{\begin{array}{ll}
\hat{f}_{i-1/2}+\Delta x/\Delta t u^{n}_{i}, & if \ u^{n+1,1}_{i}<0,\\
\hat{f}_{i+1/2}, & \text{otherwise},\\ 
\end{array}
\right.
\end{align*}
where
$$u_{i}^{n+1,1}=u^{n}_{i}-\Delta t/\Delta x (\hat{f}_{i+1/2}-\hat{\hat{f}}_{i-1/2}).$$
For periodicity, we further set $\tilde{f}_{-1/2}=\hat{\hat{f}}_{M-1/2}$. Note that since the flux at $x_{-1/2}$ is modified, the solution at $x_{0}$ may become negative. Hence, we need screen the solution one more time by the following procedure. Starting $i=0$, we compute
$$u_{i}^{n+1,2}=u^{n}_{i}-\Delta t/\Delta x (\hat{\hat{f}}_{i+1/2}-\tilde{f}_{i-1/2}).$$
If $u_{i}^{n+1,2}<0,$   we let 
$$\tilde{f}_{i+1/2}=\tilde{f}_{i-1/2}+\Delta x/\Delta t u^{n}_{i}.$$
Otherwise, let $i_0=i$ and for any $j$ with $i_0\leq j\leq M-1$, let 
$\tilde{f}_{j+1/2}=\hat{\hat{f}}_{j+1/2}.$ Finally, we are able to obtain the positivity-preserving solution $u^{new}$ via \eqref{eq:pplimiter}. Note that $i_0$ must be less than $M-1$ and hence $\tilde{f}_{M-1/2}=\tilde{f}_{-1/2}$.

When $c<0$, wind blows from right to left, and we will reset the numerical flux $\hat{f}_{i-1/2}$ if $u^{n+1}_{i}<0$. Similar idea can be used for different boundary conditions. 

Below, we show that the PP-limiter can preserve positivity while maintaining the previous high order. In particular, we have the following proposition.

\begin{prop} Consider the WENO MOL$^T$scheme is used for solving \eqref{eq:1dadv} and assume that the numerical solution $u^{n+1}_{i}$ is $p$-th order accurate in space and $q$-th order accurate in time, i.e., the global error
	$$e_i^{n+1} = |u_i^{n+1}-u(x_i,t^{n+1})|=\mathcal{O}(\Delta x^p+\Delta t^q),\quad \forall i,$$
	where $u(x_{i},t^{n+1})$ is the exact solution. If  the  proposed PP-limiter \eqref{eq:pplimiter} with \eqref{eq:limit} is applied, then 
		$$\tilde{e}_i^{n+1} = |u_i^{new}-u(x_i,t^{n+1})|=\mathcal{O}(\Delta x^{p-1}+\Delta x^{-1}\Delta t^q),\quad \forall i$$
		for the worst case, when $\Delta t=\mathcal{O}(\Delta x)$. 
\end{prop}

{\bf Proof.} Without loss of generality, we only consider \eqref{eq:1dadv} with $c>0$ and a Dirichlet boundary condition is imposed for brevity. Let $u^{n+1}_i$ be the first solution value that is less than zero starting from $u^{n+1}_{1}$, i.e., $u^{n+1}_i<0$ and  $u^{n+1}_j\ge0,\,1\le j<i.$ Then we have $\tilde{f}_{j-1/2}=\hat{f}_{j-1/2},\, 0\le j\le i,$ and $\tilde{f}_{i+1/2} = \tilde{f}_{i-1/2}+ \Delta x/\Delta t u^{n}_{i},$
such that $u^{new}_i=0.$ Recalling that $u(x_i,t^{n+1})\ge 0$ and $e_i^{n+1}=\mathcal{O}(\Delta x^p+\Delta t^q),$ and hence $\tilde{e}_i^{n+1}=\mathcal{O}(\Delta x^p+\Delta t^q).$ Also note that the modification of flux 
$$\tilde{f}_{i+1/2}-\hat{f}_{i+1/2} =\Delta x/\Delta t u^{n+1}_{i}=\mathcal{O}(\Delta x^{p}+\Delta t^{q})$$ 
is a negative high order quantity.

Now, we consider the effect of the modification of flux $\tilde{f}_{i+1/2}$ to the next point value in terms of accuracy. First, we need to compute $$u_{i+1}^{n+1,1}=u^{n}_{i+1}-\Delta t/\Delta x (\hat{f}_{i+3/2}-\tilde{f}_{i+1/2}).$$ Several cases need to be taken into account.

\noindent{\bf Case (1)} If $u_{i+1}^{n+1,1}\ge0.$ Then, $\tilde{f}_{i+3/2}=\hat{f}_{i+3/2}$ by \eqref{eq:limit}, and $u^{new}_{i+1}=u_{i+1}^{n+1,1}.$ We have
$$|u^{new}_{i+1}-u^{n+1}_{i+1}|=|\Delta t/\Delta x (\tilde{f}_{i+1/2}-\hat{f}_{i+1/2})|=\mathcal{O}(\Delta x^{p}+\Delta t^{q}).$$
Hence, $$\tilde{e}_{j+1}^{n+1}=|u^{new}_{i+1}-u(x_{i+1},t^{n+1})|\leq |u^{new}_{i+1}-u^{n+1}_{i+1}|+e_{i+1}^{n+1}=\mathcal{O}(\Delta x^{p}+\Delta t^{q}).$$ 

\noindent{\bf Case (2)}
If $u_{i+1}^{n+1,1}<0.$ Then, $\tilde{f}_{i+3/2}=\tilde{f}_{i+1/2}+ \Delta x/\Delta t u^{n}_{i+1}$ by \eqref{eq:limit}, and $u^{new}_{i+1}=0$. There are two cases we need to treat differently. If $u_{i+1}^{n+1}<0$, then by the same argument we used to treat the case at $x_{i}$, we have  $u^{n+1}_{i+1}=\mathcal{O}(\Delta x^{p}+\Delta t^{q})$  and hence $e^{n+1}_{i+1}=\mathcal{O}(\Delta x^{p}+\Delta t^{q})$. If $u_{i+1}^{n+1}\ge0$, then we have 
$$|u_{i+1}^{n+1,1}-u_{i+1}^{n+1}|=|\Delta t/\Delta x(\tilde{f}_{i+1/2}-\hat{f}_{i+1/2})|=\mathcal{O}(\Delta x^{p}+\Delta t^{q}).$$ 
Further, since  $u_{i+1}^{n+1,1}<0$ and  $u_{i+1}^{n+1}>0$, combining the above estimate, we have
$$|u^{n+1}_{i+1}|=\mathcal{O}(\Delta x^{p}+\Delta t^{q}).$$
Consequently,   $$\tilde{e}_{j+1}^{n+1}=|u^{new}_{i+1}-u(x_{i+1},t^{n+1})|=|u(x_{i+1},t^{n+1})|\leq |u^{n+1}_{i+1}|+e_{j+1}^{n+1}=\mathcal{O}(\Delta x^{p}+\Delta t^{q}).$$
We also have the following modification of the flux
\begin{equation}\label{eq:modi2}
	\tilde{f}_{i+3/2}-\hat{f}_{i+3/2} =\Delta x/\Delta t u^{n+1}_{i} + \tilde{f}_{i+1/2}-\hat{f}_{i+1/2} = \Delta x/\Delta t (u^{n+1}_{i+1} + u^{n+1}_{i})=\mathcal{O}(\Delta x^{p}+\Delta t^{q}).
\end{equation}

For both cases, we have shown that $\tilde{e}_{j+1}^{n+1}=\mathcal{O}(\Delta x^{p}+\Delta t^{q}),$ i.e., the PP-limiter preserves the order of accuracy.  Moreover, for case (1)
since the flux $\hat{f}_{i+3/2}$ is not changed, we can move on to consider the solution value at $x_{i+2}$ with the same procedure as the case at $x_{i}$. The situation of case (2) is more complicated when considering the subsequent effect of the modification of flux \eqref{eq:modi2}. Again, we need to consider two cases.

\noindent{\bf Case (2.1)} If $u_{i+2}^{n+1,1}\ge0.$ Similar to case (1) above, the flux at $\tilde{f}_{i+5/2}$ should remain the same, i.e., $\tilde{f}_{i+5/2}=\hat{f}_{i+5/2}$. By \eqref{eq:modi2}
$$|u^{new}_{i+2}-u^{n+1}_{i+2}|=|\Delta t/\Delta x (\tilde{f}_{i+3/2}-\hat{f}_{i+3/2})|=\mathcal{O}(\Delta x^{p}+\Delta t^{q}),$$
and consequently, $$\tilde{e}_{i+2}^{n+1}=\mathcal{O}(\Delta x^{p}+\Delta t^{q}).$$

\noindent{\bf Case (2.2)} If $u_{i+2}^{n+1,1}<0$. We can repeat the analysis of case (2) and obtain that 
$$|u^{n+1}_{i+2}|=\mathcal{O}(\Delta x^{p}+\Delta t^{q}),$$
which gives the estimate
$$\tilde{e}_{i+2}^{n+1}=\mathcal{O}(\Delta x^{p}+\Delta t^{q}),$$ and the modification of the flux
\begin{equation}\label{eq:modi3}
\tilde{f}_{i+5/2}-\hat{f}_{i+5/2} = \Delta x/\Delta t (u^{n+1}_{i+2}+u^{n+1}_{i+1} + u^{n+1}_{i})=\mathcal{O}(\Delta x^{p}+\Delta t^{q}).
\end{equation}

More generally, assume that the fluxes are successively modified with at most constant many times, i.e., there exists a constant $K_0$ independent of $\Delta x$ and $\Delta t$, such that if the fluxes $\tilde{f}_{i+k+1/2}$ for $0\le k\le K$ are successively changed based on
$$\tilde{f}_{i+k+1/2} = \hat{f}_{i+k+1/2} +\Delta x/\Delta t \sum_{l=0}^ku^{n+1}_{i+l},$$
then $K$ is always bounded by $K_0$. We have $$\tilde{e}^{n+1}_i=\mathcal{O}(\Delta x^{p}+\Delta t^{q}).$$ For the worst cases, the fluxes may be successively modified with $\mathcal{O}(1/\Delta x)$ times, then the accumulation of error leads to one order reduction in space, i.e., 
$$\tilde{e}^{n+1}_i=\mathcal{O}(\Delta x^{p-1}+\Delta x^{-1}\Delta t^{q}).$$

\begin{rem}
	In our simulations, the worst cases never happened: $K_0$ is always a constant independent of mesh sizes, and optimal accuracy is observed for all our accuracy tests.
\end{rem}
\section{Two-dimensional implementation}

In this section, we extend the algorithm to two-dimensional problems in the following form:
\begin{equation}
\label{eq:2d}
u_{t}+f(y,t)u_{x}+g(x,t)u_{y}=0
\end{equation}
In particular, \eqref{eq:2d} is solved in a splitting framework. First, we decouple the 2D advection equation into a sequence of  independent  1D advection equations
\begin{subequations}
\begin{align}
u_t + f(y,t)u_x &= 0,\label{eq:splitx}\\
u_t + g(x,t)u_y &=0.\label{eq:splity}
\end{align}
\end{subequations}
Let $Q_{1}$ and $Q_{2}$ be the approximate solvers for each transport problem, that is $u^{n+1}=Q_{1}(\Delta t)u^{n}$ is an approximation of \eqref{eq:splitx}, and $u^{n+1}=Q_{2}(\Delta t)u^{n}$ is an approximation of \eqref{eq:splity}. The splitting form is first order accurate in time if one evolves \eqref{eq:splitx} for a full time step, then evolves \eqref{eq:splity} for a full time step, denoted as $u^{n+1}=Q_{2}(\Delta t)\cdot Q_{1}(\Delta t) u^{n}$. High order semi-discrete scheme for one step evolution from $u^{n}$ to $u^{n+1}$ is formulated in \cite{forest1990fourth}. For example, the second order splitting, i.e., Strange splitting, is
\begin{align}
\label{eq:sp2}
u^{n+1}=Q_{1}(\frac{1}{2}\Delta t)\cdot Q_{2}(\Delta t)\cdot Q_{1}(\frac{1}{2}\Delta t) u^{n}.
\end{align}
The third order splitting is
\begin{align}
\label{eq:sp3}
u^{n+1}=Q_{2}(\Delta t)\cdot Q_{1}(-\frac{1}{24}\Delta t)\cdot Q_{2}(-\frac{2}{3}\Delta t)\cdot Q_{1}(\frac{3}{4}\Delta t)\cdot  Q_{2}(\frac{2}{3}\Delta t)\cdot Q_{1}(\frac{7}{24}\Delta t) u^{n}.
\end{align}
And the forth order splitting is
\begin{align}
\label{eq:sp4}
u^{n+1}=& Q_{2}((\alpha+1/2)\Delta t)\cdot Q_{1}((2\alpha+1)\Delta t)\cdot Q_{2}(-\alpha\Delta t)\cdot Q_{1}(-(4\alpha+1)\Delta t)\cdot\nonumber\\
&  Q_{2}(-\alpha\Delta t)\cdot Q_{1}((2\alpha+1)\Delta t)\cdot Q_{2}((\alpha+1/2)\Delta t) u^{n},
\end{align}
where $\alpha=(2^{1/3}+2^{-1/3}-1)/6$. 
 
 Note that the 1D1V VP system \eqref{eq:vp} can be split into two 1D advection equations:
 \begin{subequations}
 \begin{align}
& f_{t}+vf(x)=0, \label{eq:vpx}\\
& f_{t}+E(x,t)f_{v}=0, \label{eq:vpv}
 \end{align}
 \end{subequations}
 known as spatial advection and velocity acceleration, respectively.
Similarly, we can get the high order approximation by using \eqref{eq:sp2}, \eqref{eq:sp3} or \eqref{eq:sp4}. Further, the electrostatic field $E(x,t)$ is obtained by solving Poisson's equation \eqref{eq:possion} when it is needed.
\section{Numerical Examples}

In this section, we present the results of our numerical experiments for the schemes described in the previous sections. We will use the third and the fifth order WENO integrations for spatial discretization, denoted as WENO3 and WENO5. And the time discretization $RK(2,3)$ and $RK(4,4)$ are used for WENO3 and WENO5, respectively, to guarantee third order and fourth order accuracy in time. For 1D problems, the time steps are defined as 
$$\Delta t=CFL\Delta x/|c|,$$ 
where $CFL=1.5$ for $RK(2,3)$ and $CFL=2.9$ for $RK(4,4)$. While for 2D problems, the time steps are chosen as
$$\Delta t=CFL/\max(|c_{x}|/\Delta x,|c_{y}|/ \Delta y)$$
where $|c_{x}|$ and $|c_{y}|$ are the maximum wave propagation speeds in the $x$- and $y$-directions, and $\Delta x$ and $\Delta y$ are grid sizes in the $x$- and $y$-directions, respectively. For $RK(2,3)$, we use the third order time splitting scheme \eqref{eq:sp3} with $CFL=1.5$. While, for $RK(4,4)$, we use the forth order time splitting scheme \eqref{eq:sp4} with $CFL=1.6$.

To avoid the effect of rounding error, we change the condition \eqref{eq:pplimiter} into
$$u^n_{i}-\frac{\Delta t}{\Delta x}(\tilde{f}_{i+1/2}-\tilde{f}_{i-1/2})\geq 10^{-16}.$$

\subsection{Basic 1-D and 2-D tests}

\textbf{ Example 4.1} (one-dimensional linear advection equation)
\begin{align}
\label{eq:prob1}
u_{t}+u_{x}=0, \ \ \  x\in [-\pi,\pi]
\end{align}
with two different initial conditions:
\begin{itemize}
\item continuous initial condition $u(x,0)=\cos(x)^4$, 
\item discontinuous initial condition:
$
u(x,0)=\left\{\begin{array}{ll}
1, & x\in[-\frac{1}{4}\pi,\frac{1}{4}\pi]\\
0, & \text{otherwise}.\\
\end{array}
\right.
$
\end{itemize}

\begin{figure}
\centering
\subfigure[Periodic boundary condition.]{
\includegraphics[width=0.4\textwidth]{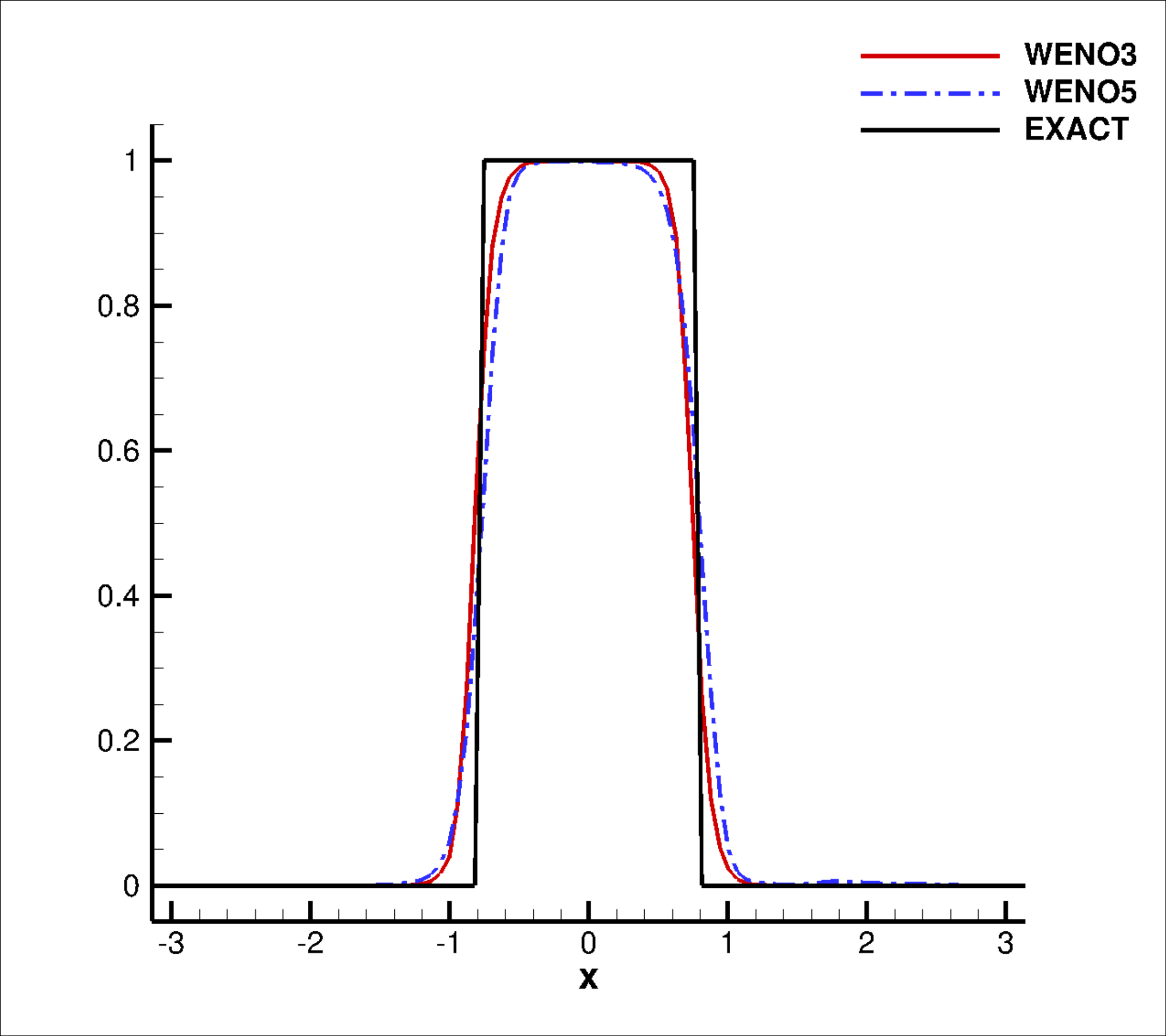}}
\subfigure[Dirichlet boundarycondition.]{
\includegraphics[width=0.4\textwidth]{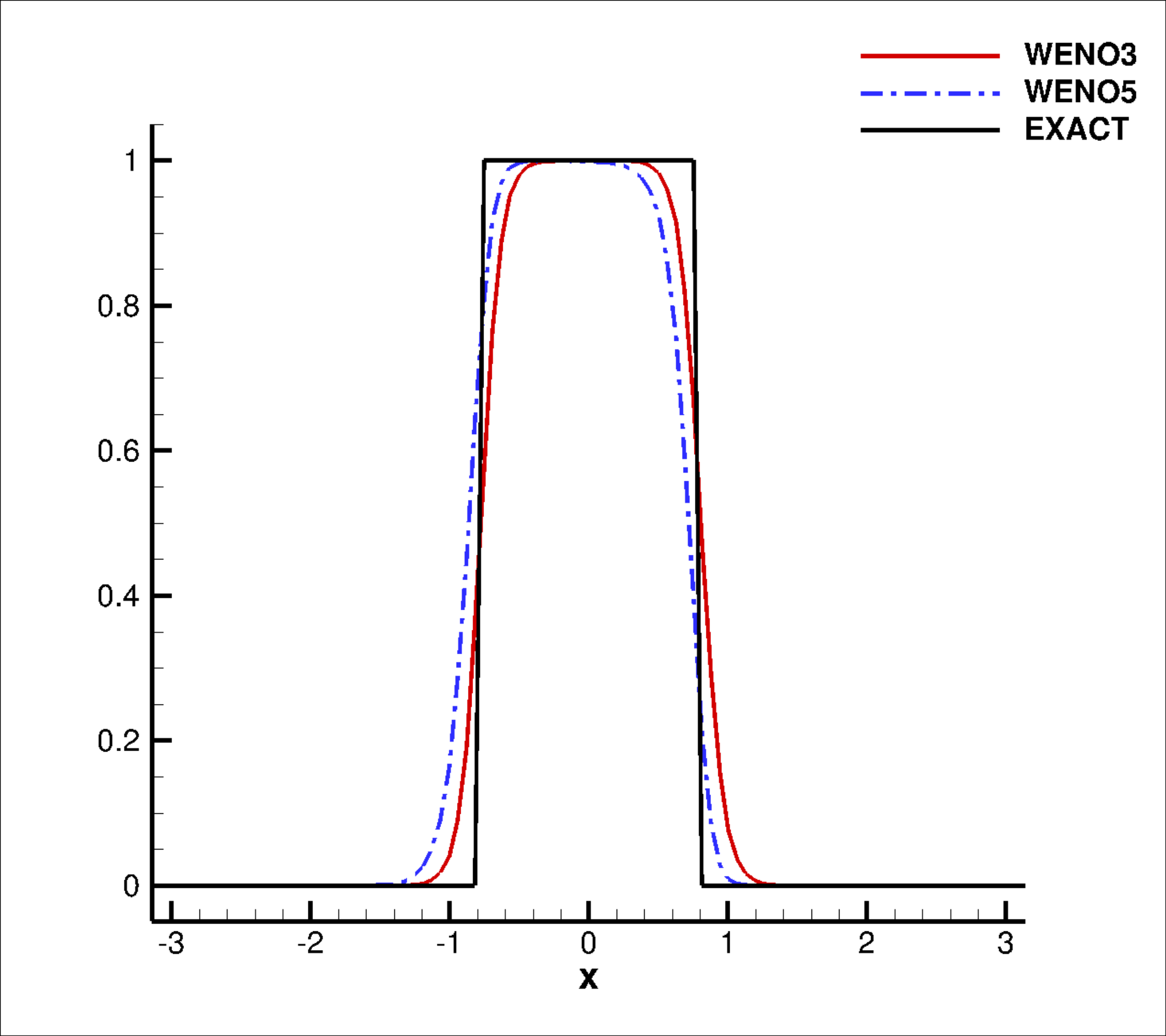}}
\caption{\em The one-dimensional linear advection equation with discontinuous initial condition. $T=2\pi$. $100$ grid points.}
\label{Fig1}
\end{figure}

The initial condition can be extended periodically on $R$ with period $2\pi$. The exact solution is $u(x,t)=u(x-t,0)$. For the continuous problem, we test with a $2\pi$-periodic boundary condition, a Dirichlet boundary condition $u(-\pi,t)=\cos(-\pi-t)^4$ and a Neumann boundary condition $u_{x}(-\pi,t)=-4\cos(-\pi-t)^3 \sin(-\pi-t)$. The $L_{1}$ and $L_{\infty}$ errors and orders of accuracy are listed in Table \ref{tab1}, \ref{tab2} and \ref{tabadd1}, and the PP-limiter is applied. We observe that WENO3 with $RK(2,3)$ exhibits third order accuracy, and WENO5 with $RK(4,4)$ exhibits the fourth order accuracy. The PP-limiter can keep the solutions non-negative without loss of accuracy. While for the discontinuous case (Figure \ref{Fig1}), we consider the $2\pi$-periodic boundary condition and the following Dirichlet boundary condition
\begin{align*}
u(-\pi,t)=\left\{\begin{array}{ll}
1, & t\in[\frac{3}{4}\pi,\frac{5}{4}\pi],\\
0, & t\in(0,\frac{3}{4}\pi)\cup (\frac{5}{4}\pi,2\pi).\\
\end{array}
\right.
\end{align*} 
It is observed that the WENO methodology used effectively suppress the unphysical oscillations for both schemes.

\begin{table}[htb]
\caption{\label{tab1}\em Accuracy on the one-dimensional linear advection equation with continuous initial condition. $T=2\pi$. Periodic boundary condition. }
\centering
\bigskip
\begin{small}
\begin{tabular}{|c|c|c|c|c|c|c|c|}
  \hline
 &  &  $N_x$ &  $L_1$ errors &   order  &  $L_\infty$ error  & order & min value \\\hline
 \multirow{12}{1.5cm}{WENO3} & \multirow{6}{2cm}{Without\\ PP-limiters}  &
           20 &  6.29E-01  &    --    &  2.36E-01  &     --    &  7.11E-03  \\
 & &    40 &  1.21E-01  &  2.38  &  4.29E-02  &  2.46  &  -7.30E-03  \\
 & &    80 &  1.32E-02  &  3.20  &  4.77E-03  &  3.17  &  -1.48E-03  \\
 & &  160 &  1.25E-03  &  3.39  &  5.21E-04  &  3.19  &  -2.32E-04  \\
 & &  320 &  1.17E-04  &  3.42  &  4.40E-05  &  3.57  &  -2.74E-05  \\
 & &  640 &  1.50E-05  &  2.96  &  4.83E-06  &  3.19  &  -2.89E-06  \\\cline{2-8}
 & \multirow{6}{2cm}{With\\ PP-limiters}  &
 		  20 &  6.28E-01  &    --    &  2.33E-01  &    --     &  1.35E-02  \\
 & &    40 &  1.22E-01  &  2.37  &  4.38E-02  &  2.41  &  1.11E-16  \\
 & &    80 &  1.29E-02  &  3.24  &  4.95E-03  &  3.15  &  1.11E-16  \\
 & &  160 &  1.19E-03  &  3.44  &  4.87E-04  &  3.35  &  1.11E-16  \\
 & &  320 &  1.17E-04  &  3.35  &  5.92E-05  &  3.04  &  1.11E-16  \\
 & &  640 &  1.50E-05  &  2.96  &  8.77E-06  &  2.76  &  1.11E-16  \\\hline
  \multirow{12}{1.5cm}{WENO5} & \multirow{6}{2cm}{Without\\ PP-limiters}  &
           20 &  4.17E-01  &    --    &  1.57E-01  &     --    &  -4.28E-02  \\
 & &    40 &  3.63E-02  &  3.52  &  1.95E-02  &  3.01  &  -5.50E-03  \\
 & &    80 &  2.26E-03  &  4.01  &  1.32E-03  &  3.89  &  -4.79E-05  \\
 & &  160 &  1.24E-04  &  4.19  &  8.48E-05  &  3.96  &  -6.48E-06  \\
 & &  320 &  5.49E-06  &  4.50  &  2.05E-06  &  5.37  &  -3.83E-07  \\
 & &  640 &  3.26E-07  &  4.08  &  8.86E-08  &  4.53  &  -1.01E-09  \\\cline{2-8}
 & \multirow{6}{2cm}{With\\ PP-limiters}  &
 		  20 &  5.06E-01  &    --    &  2.00E-01  &    --     &  1.59E-02 \\
 & &    40 &  5.10E-02  &  3.31  &  2.93E-02  &  2.77  &  7.06E-04  \\
 & &    80 &  6.01E-03  &  3.08  &  4.86E-03  &  2.59  &  3.36E-04  \\
 & &  160 &  1.55E-04  &  5.28  &  1.28E-04  &  5.25  &  7.28E-06  \\
 & &  320 &  5.43E-06  &  4.84  &  1.75E-06  &  6.19  &  4.75E-08  \\
 & &  640 &  3.26E-07  &  4.06  &  8.86E-08  &  4.30  &  1.73E-10  \\\hline
\end{tabular}
\end{small}
\end{table}

\begin{table}[htb]
\caption{\label{tab2}\em Accuracy on the one-dimensional linear advection equation with continuous initial condition. $T=2\pi$. Dirichlet boundary condition. }
\centering
\bigskip
\begin{small}
\begin{tabular}{|c|c|c|c|c|c|c|c|}
  \hline
 &  &  $N_x$ &  $L_1$ errors &   order  &  $L_\infty$ error  & order & min value \\\hline
 \multirow{12}{1.5cm}{WENO3} & \multirow{6}{2cm}{Without\\ PP-limiters}  &
     	  20 &  4.55E-01  &    --     &  3.11E-01  &    --    &  -4.42E-03  \\
 & &    40 &  9.09E-02  &  2.32  &  6.66E-02  &  2.22  &  -1.28E-03  \\
 & &    80 &  1.10E-02  &  3.05  &  1.15E-02  &  2.54  &  -3.71E-04  \\
 & &  160 &  1.08E-03  &  3.35  &  2.46E-03  &  2.22  &  -1.09E-04  \\
 & &  320 &  8.87E-05  &  3.60  &  1.95E-04  &  3.66  &  -1.82E-05  \\
 & &  640 &  9.63E-06  &  3.20  &  6.09E-06  &  5.00  &  -2.15E-06  \\\cline{2-8}
  & \multirow{6}{2cm}{With\\ PP-limiters}  &
      	  20 &  4.53E-01  &    --     &  3.11E-01  &    --    &  1.11E-16  \\
 & &    40 &  9.07E-02  &  2.32  &  6.66E-02  &  2.22  &  5.55E-17  \\
 & &    80 &  1.10E-02  &  3.04  &  1.15E-02  &  2.54  &  5.55E-17  \\
 & &  160 &  1.06E-03  &  3.38  &  2.46E-03  &  2.22  &  5.55E-17  \\
 & &  320 &  8.87E-05  &  3.58  &  1.95E-04  &  3.66  &  1.11E-16  \\
 & &  640 &  9.63E-06  &  3.20  &  6.61E-06  &  4.88  &  5.55E-17  \\
\hline
 \multirow{12}{1.5cm}{WENO5} & \multirow{6}{2cm}{Without\\ PP-limiters}  &
     	  20 &  3.59E-01  &    --    &  2.14E-01  &     --    &  1.18E-02  \\
 & &    40 &  4.37E-02  &  3.04  &  3.24E-02  &  2.73  &  9.44E-04  \\
 & &    80 &  2.65E-03  &  4.05  &  5.79E-03  &  2.48  &  1.98E-05  \\
 & &  160 &  1.08E-04  &  4.61  &  1.69E-04  &  5.10  &  -1.52E-06  \\
 & &  320 &  2.91E-06  &  5.22  &  1.70E-06  &  6.64  &  -1.99E-07  \\
 & &  640 &  1.69E-07  &  4.10  &  8.22E-08  &  4.37  &  6.37E-09  \\\cline{2-8}
  & \multirow{6}{2cm}{With\\ PP-limiters}  &
  	      20 &  3.58E-01  &    --     &  2.15E-01  &    --    &  1.18E-02  \\
 & &    40 &  4.37E-02  &  3.04  &  3.24E-02  &  2.73  &  1.29E-03  \\
 & &    80 &  2.83E-03  &  3.95  &  5.79E-03  &  2.48  &  3.53E-05  \\
 & &  160 &  1.18E-04  &  4.57  &  1.75E-04  &  5.05  &  6.26E-06  \\
 & &  320 &  2.89E-06  &  5.36  &  1.67E-06  &  6.71  &  1.11E-16  \\
 & &  640 &  1.69E-07  &  4.10  &  8.22E-08  &  4.35  &  6.37E-09  \\
\hline
\end{tabular}
\end{small}
\end{table}

\begin{table}[htb]
\caption{\label{tabadd1}\em Accuracy on the one-dimensional linear advection equation with continuous initial condition. $T=2\pi$. Neumann boundary condition. }
\centering
\bigskip
\begin{small}
\begin{tabular}{|c|c|c|c|c|c|c|c|}
  \hline
 &  &  $N_x$ &  $L_1$ errors &   order  &  $L_\infty$ error  & order & min value \\\hline
 \multirow{12}{1.5cm}{WENO3} & \multirow{6}{2cm}{Without\\ PP-limiters}  &
           20 &  4.52E-01  &    --    &  3.09E-01  &     --    &  -7.98E-03  \\
 & &    40 &  9.14E-02  &  2.31  &  6.65E-02  &  2.22  &  -1.22E-03  \\
 & &    80 &  1.11E-02  &  3.05  &  1.15E-02  &  2.54  &  -4.02E-04  \\
 & &  160 &  1.09E-03  &  3.34  &  2.46E-03  &  2.22  &  -1.10E-04  \\
 & &  320 &  9.09E-05  &  3.58  &  1.95E-04  &  3.66  &  -1.83E-05  \\
 & &  640 &  9.93E-06  &  3.20  &  6.09E-06  &  5.00  &  -2.16E-06  \\\cline{2-8}
 & \multirow{6}{2cm}{With\\ PP-limiters}  &
 		  20 &  4.50E-01  &    --    &  3.09E-01  &    --     &  1.11E-16  \\
 & &    40 &  9.12E-02  &  2.30  &  6.64E-02  &  2.22  &  5.55E-17  \\
 & &    80 &  1.12E-02  &  3.04  &  1.15E-02  &  2.54  &  5.55E-17  \\
 & &  160 &  1.07E-03  &  3.37  &  2.46E-03  &  2.22  &  5.55E-17  \\
 & &  320 &  9.10E-05  &  3.56  &  1.95E-04  &  3.66  &  1.11E-16  \\
 & &  640 &  9.93E-06  &  3.20  &  6.65E-06  &  4.87  &  5.55E-17  \\\hline
  \multirow{12}{1.5cm}{WENO5} & \multirow{6}{2cm}{Without\\ PP-limiters}  &
           20 &  3.92E-01  &    --    &  2.22E-01  &     --    &  -2.93E-02  \\
 & &    40 &  4.75E-02  &  3.05  &  3.18E-02  &  2.80  &  -1.87E-03  \\
 & &    80 &  5.02E-03  &  3.24  &  5.74E-03  &  2.47  &  -4.70E-05  \\
 & &  160 &  3.13E-04  &  4.01  &  6.20E-04  &  3.21  &   9.64E-06  \\
 & &  320 &  4.11E-06  &  6.25  &  1.91E-06  &  8.34  &  -4.50E-07  \\
 & &  640 &  2.51E-07  &  4.03  &  1.20E-07  &  4.00  &  -1.00E-08  \\\cline{2-8}
 & \multirow{6}{2cm}{With\\ PP-limiters}  &
 		  20 &  3.95E-01  &    --    &  2.45E-01  &    --     &  4.85E-03 \\
 & &    40 &  5.74E-02  &  2.78  &  4.93E-02  &  2.31  &  0.00E+00  \\
 & &    80 &  6.32E-03  &  3.18  &  6.80E-03  &  2.86  &  1.56E-04  \\
 & &  160 &  2.35E-04  &  4.75  &  4.00E-04  &  4.09  &  1.28E-05  \\
 & &  320 &  3.93E-06  &  5.90  &  1.91E-06  &  7.71  &  1.11E-16  \\
 & &  640 &  2.29E-07  &  4.10  &  1.20E-07  &  4.00  &  0.00E+00  \\\hline
\end{tabular}
\end{small}
\end{table}

\textbf{Example 4.2} (Rigid body rotation)
\begin{align}
\label{eq:prob2}
u_{t}+yu_{x}-xu_{y}=0, \ \ \ x,y\in\Omega
\end{align}
with two different initial conditions:
\begin{itemize}
\item
continuous initial condition, with $\Omega=[-\frac{1}{2}\pi,\frac{1}{2}\pi]^2$
\begin{align*}
u(x,y,0)=0.5B(\sqrt{x^2+8y^2})+0.5B(\sqrt{8x^2+y^2}),
\end{align*}
where
$
B(r)=\left\{\begin{array}{ll}
\cos(r)^6,& \text{if} \ r\leq \frac{1}{2}\pi,\\
0, & \text{otherwise}.\\
\end{array}
\right.
$

\item dicontinuous initial condition, with $\Omega=[-1,1]^2$
\begin{align*}
u(x,y,0)=\left\{\begin{array}{ll}
1, & (x,y)\in[-0.75,0.75]\times[-0.25,0.25]\bigcup[-0.25,0.25]\times[-0.75,0.75],\\
0, & \text{otherwise}.\\
\end{array}
\right.
\end{align*}
\end{itemize}

\begin{figure}
\centering
\subfigure[WENO3.]{
\includegraphics[width=0.4\textwidth]{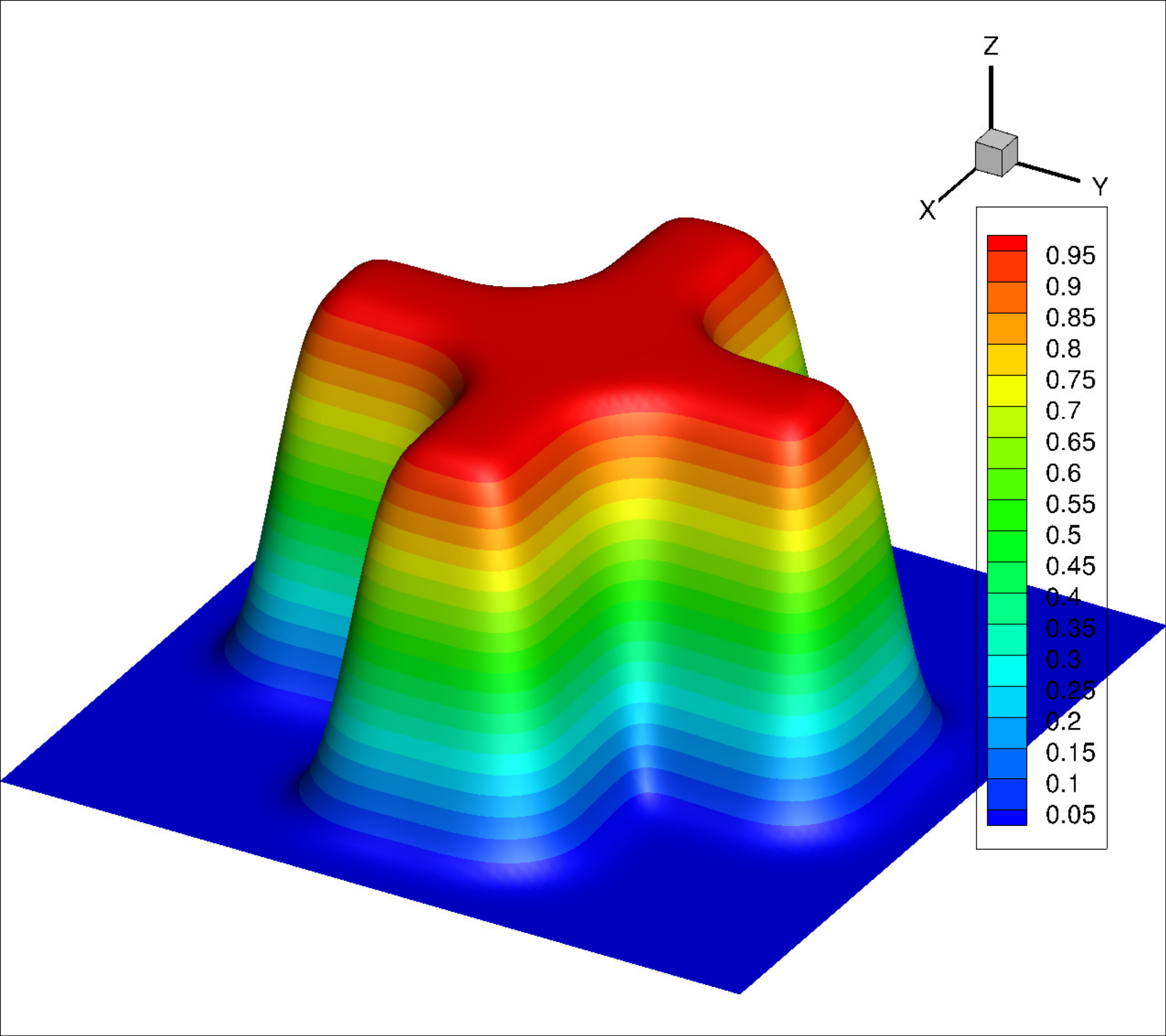}}
\subfigure[WENO5.]{
\includegraphics[width=0.4\textwidth]{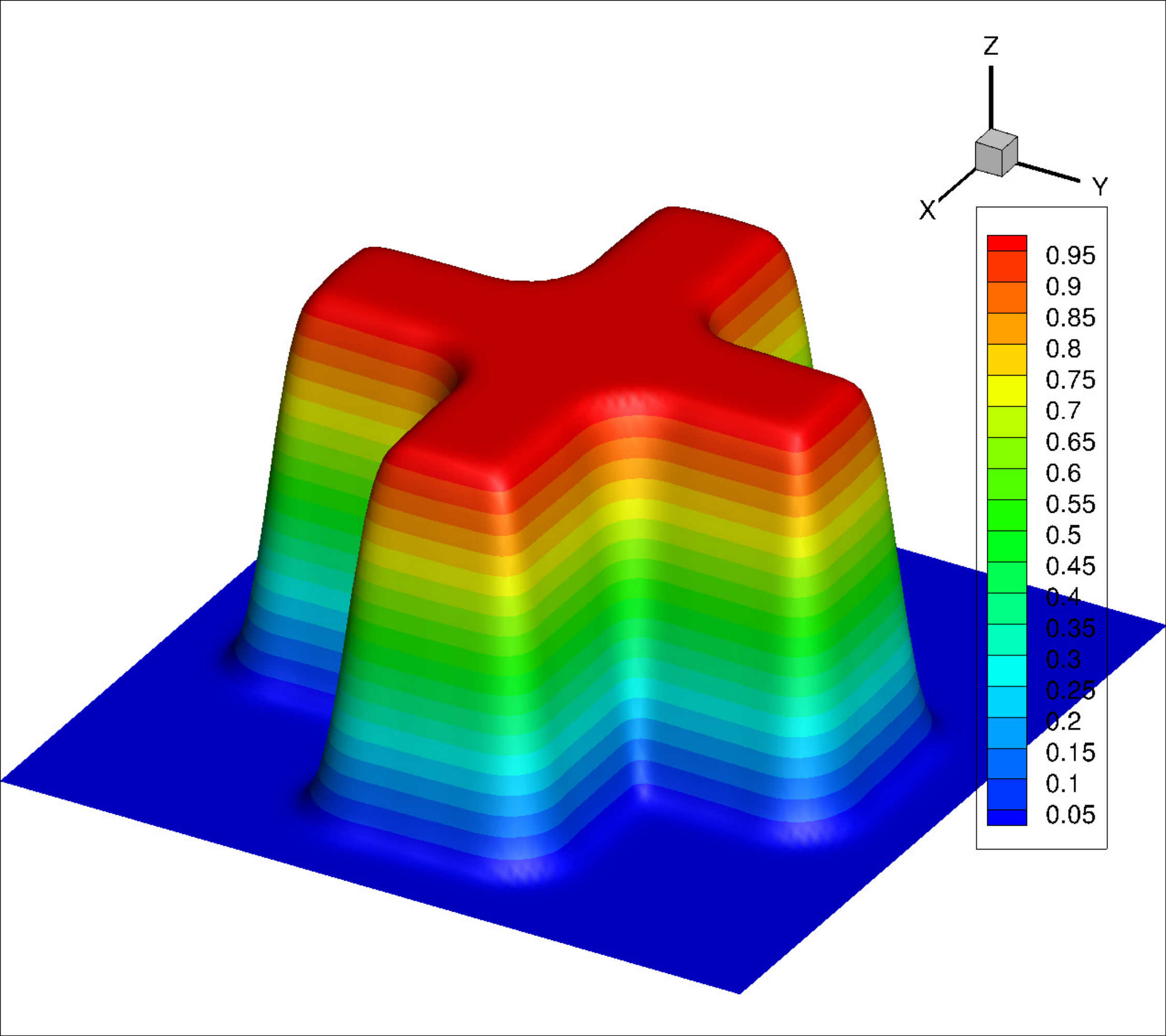}}
\caption{\em The rigid body rotation with discontinuous initial condition. $T=2\pi$. $100\times100$ grid points.}
\label{Fig2}
\end{figure}

We simulate the problems with the zero Dirichlet boundary conditions for both problems. In Table \ref{tab4}, we summarize the convergence study and minimum solution values at the final time $T=2\pi$. It is observed  that our schemes can achieve the designed order of accuracy and meanwhile the PP-limiter is able to keep positivity without destroying the accuracy of the original scheme. Furthermore, similar to 1D case, WENO methodology removes unphysical oscillations as expected (Figure \ref{Fig2}).

\begin{table}[htb]
\caption{\label{tab4}\em Accuracy on the rigid body rotation with continuous initial condition. $T=2\pi$.}
\centering
\bigskip
\begin{small}
\begin{tabular}{|c|c|c|c|c|c|c|c|}
  \hline
 &  &  $N_x\times N_y$ &  $L_1$ errors &   order  &  $L_\infty$ error  & order & min value \\\hline
 \multirow{10}{1.5cm}{WENO3} & \multirow{5}{2cm}{Without\\ PP-limiters}  &
 	       $20\times20$ &  1.84E-01  &     --    &  4.25E-01  &   --     &  -2.22E-05  \\
 & &     $40\times40$ &  9.84E-02  &  0.90  &  1.72E-01  &  1.31  &  -1.41E-04  \\
 & &     $80\times80$ &  2.66E-02  &  1.89  &  3.35E-02  &  2.36  &  -2.46E-04  \\
 & & $160\times160$ &  3.41E-03  &  2.96  &  4.40E-03  &  2.93  &  -2.84E-04  \\
 & & $320\times320$ &  3.96E-04  &  3.10  &  5.20E-04  &  3.08  &  -7.25E-05  \\\cline{2-8}
 & \multirow{5}{2cm}{With\\ PP-limiters}  &
 	       $20\times20$ &  1.84E-01  &     --    &  4.25E-01  &    --     &  0.00E+00  \\
 & &     $40\times40$ &  9.84E-02  &  0.90  &  1.72E-01  &  1.31  &  0.00E+00  \\
 & &     $80\times80$ &  2.65E-02  &  1.89  &  3.35E-02  &  2.36  &  0.00E+00  \\
 & & $160\times160$ &  3.32E-03  &  3.00  &  4.40E-03  &  2.93  &  0.00E+00  \\
 & & $320\times320$ &  3.92E-04  &  3.09  &  5.20E-04  &  3.08  &  0.00E+00  \\\hline
 \multirow{10}{1.5cm}{WENO5} & \multirow{5}{2cm}{Without\\ PP-limiters}  &
     	   $20\times20$ &  1.43E-01  &    --    &  3.32E-01  &    --     &  -9.51E-07  \\  	
 & &     $40\times40$ &  6.08E-02  &  1.24  &  7.43E-02  &  2.16  &  -1.67E-05  \\
 & &     $80\times80$ &  4.98E-03  &  3.61  &  7.49E-03  &  3.31  &  -2.50E-05  \\
 & & $160\times160$ &  1.34E-04  &  5.22  &  2.75E-04  &  4.77  &  -4.76E-06  \\
 & & $320\times320$ &  3.73E-06  &  5.16  &  7.56E-06  &  5.18  &  -1.03E-07  \\\cline{2-8}
 & \multirow{5}{2cm}{With\\ PP-limiters}  &
     	   $20\times20$ &  1.43E-01  &    --     &  3.32E-01  &    --    &  0.00E+00  \\  	
 & &     $40\times40$ &  6.09E-02  &  1.24  &  7.43E-02  &  2.16  &  0.00E+00  \\
 & &     $80\times80$ &  4.98E-03  &  3.61  &  7.49E-03  &  3.31  &  0.00E+00  \\
 & & $160\times160$ &  1.34E-04  &  5.22  &  2.75E-04  &  4.77  &  0.00E+00  \\
 & & $320\times320$ &  3.72E-06  &  5.17  &  7.56E-06  &  5.18  &  0.00E+00  \\
  \hline
\end{tabular}
\end{small}
\end{table}

\subsection{The VP systems} 

In this subsection, we consider the Vlasov-Poisson equation \eqref{eq:vp} with the following four initial conditions:
\begin{itemize}

\item Strong Landau damping:
\begin{align*}
f(x,v,0)=\frac{1}{\sqrt{2\pi}}(1+\alpha\cos(kx))\exp(-\frac{v^2}{2}), 
x\in[0,L], \ v\in[-V_{c},V_{c}]
\end{align*}
where $\alpha=0.5$, $k=0.5$, $L=4\pi$ and $V_{c}=2\pi$.

\item Two-stream instability I:
\begin{align*}
f(x,v,0)=\frac{2}{7\sqrt{2\pi}} (1+5v^2)(1+\alpha((\cos(2kx)+\cos(3kx))/1.2 +cos(kx))) \exp(-\frac{v^2}{2}), \\
x\in[0,L], \ v\in[-V_{c},V_{c}]
\end{align*}
where $\alpha=0.01$, $k=0.5$, $L=4\pi$ and $V_{c}=2\pi$.

\item Two-stream instability II:
\begin{align*}
f(x,v,0)=\frac{1}{\sqrt{2\pi}} (1+\alpha\cos(kx))v^2 \exp(-\frac{v^2}{2}), \ \ \ 
x\in[0,L], \ v\in[-V_{c},V_{c}]
\end{align*}
where $\alpha=0.05$, $k=0.5$,  $L=4\pi$ and $V_{c}=2\pi$.

\item Bump-on-tail instability:
\begin{align*}
f(x,v,0)=\frac{1}{\sqrt{2\pi}}(1+\alpha\cos(kx)) (0.9\exp(-0.5v^2)+0.2\exp(-4(v-4.5)^2)) \nonumber\\
x\in[-L,L], \ v\in[-V_{c},V_{c}]
\end{align*}
where $\alpha=0.04$, $k=0.3$,  $L=\frac{10}{3}\pi$ and $V_{c}=10$.
\end{itemize}

In our numerical tests, periodic boundary conditions are imposed in both directions, and a fast Fourier transform (FFT) is used to solve the 1-D Poisson equation. The VP system is time-reversible, meaning that if the solution $f$ is first evolved to time $T/2$, at when we reverse the solution in $v$-direction, i.e., $\tilde{f}(x,v,\frac{1}{2}T)=f(x,-v,\frac{1}{2}T)$, and continue to evolve the reserved solution $\tilde{f}$ to time $T$, then we can recover the initial condition with reverse velocity field, i.e., $\tilde{f}(x,v,T)=f(x,-v,0)$. Below, we make use of such a property of the VP system to perform the accuracy test.
In Table \ref{tab5}-\ref{tab8}, we list the errors of $\tilde{f}(x,-v,T)$ and associated orders of accuracy at time $T=10$, without or with the PP-limiter. As expected, all schemes achieve the designed order accuracy and the limiter does not adversely affect the accuracy.

Below, we briefly recall some classical preservation results in the VP system. We hope that our numerical solutions can preserve these classical conserved quantities as much as possible.
\begin{align*}
& Mass=\int_{v}\int_{x} f(x,v,t) dxdv,\\
& \|f\|_{L_{p}}=(\int_{v}\int_{x} |f(x,v,t)|^{p} dxdv)^{1/p}, \ \ \ 1\leq p< \infty, \\
& Energy= \frac{1}{2}(\int_{v}\int_{x}f(x,v,t)v^2dxdv+\int_{x}E^{2}(x,t)dx), \\
& Momentum= \int_{v}\int_{x} f(x,v,t)vdxdv,
\end{align*}

\begin{figure}
\centering
\subfigure[$T=10$]{
\includegraphics[width=0.4\textwidth]{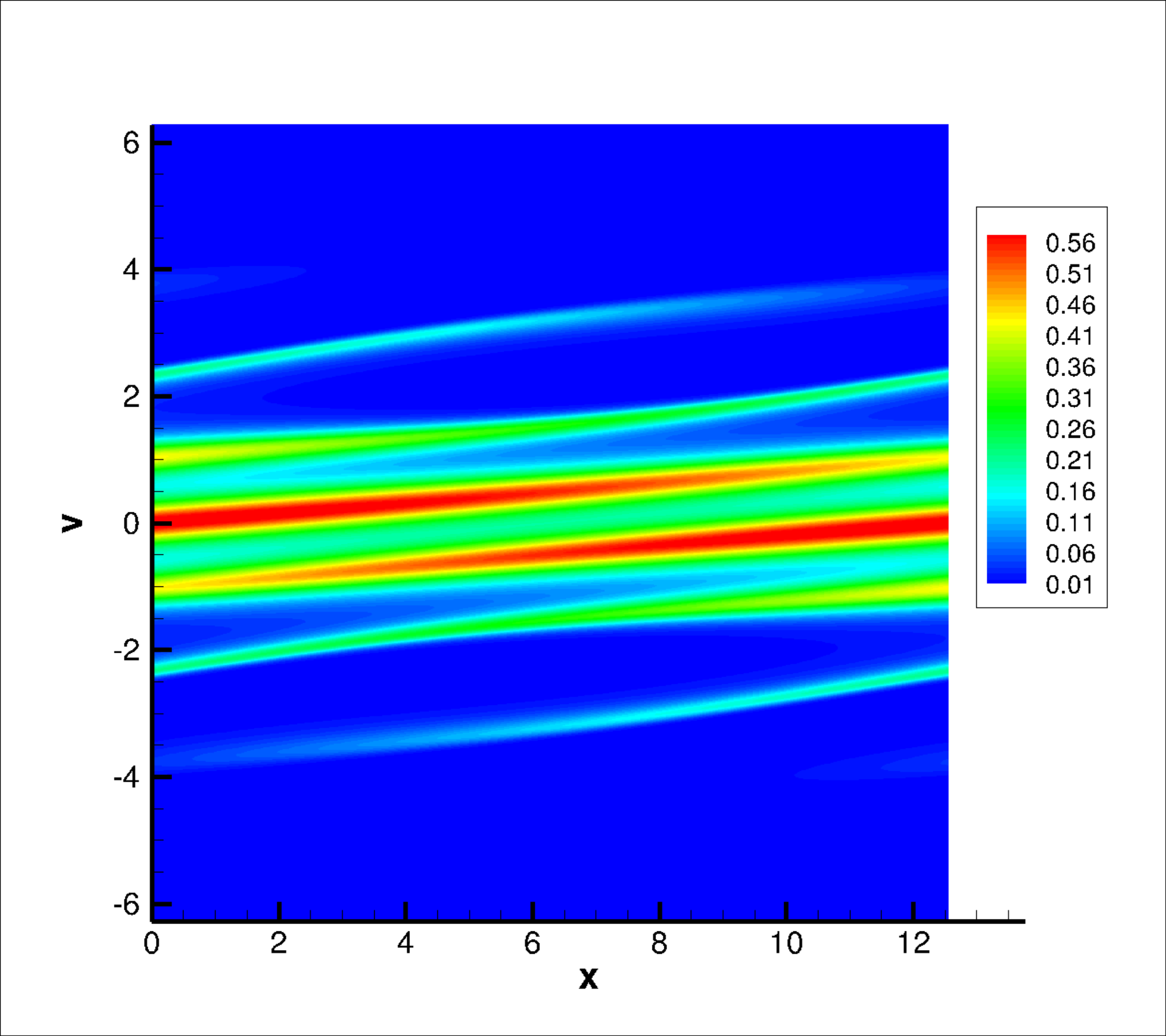}}
\subfigure[$T=20$]{
\includegraphics[width=0.4\textwidth]{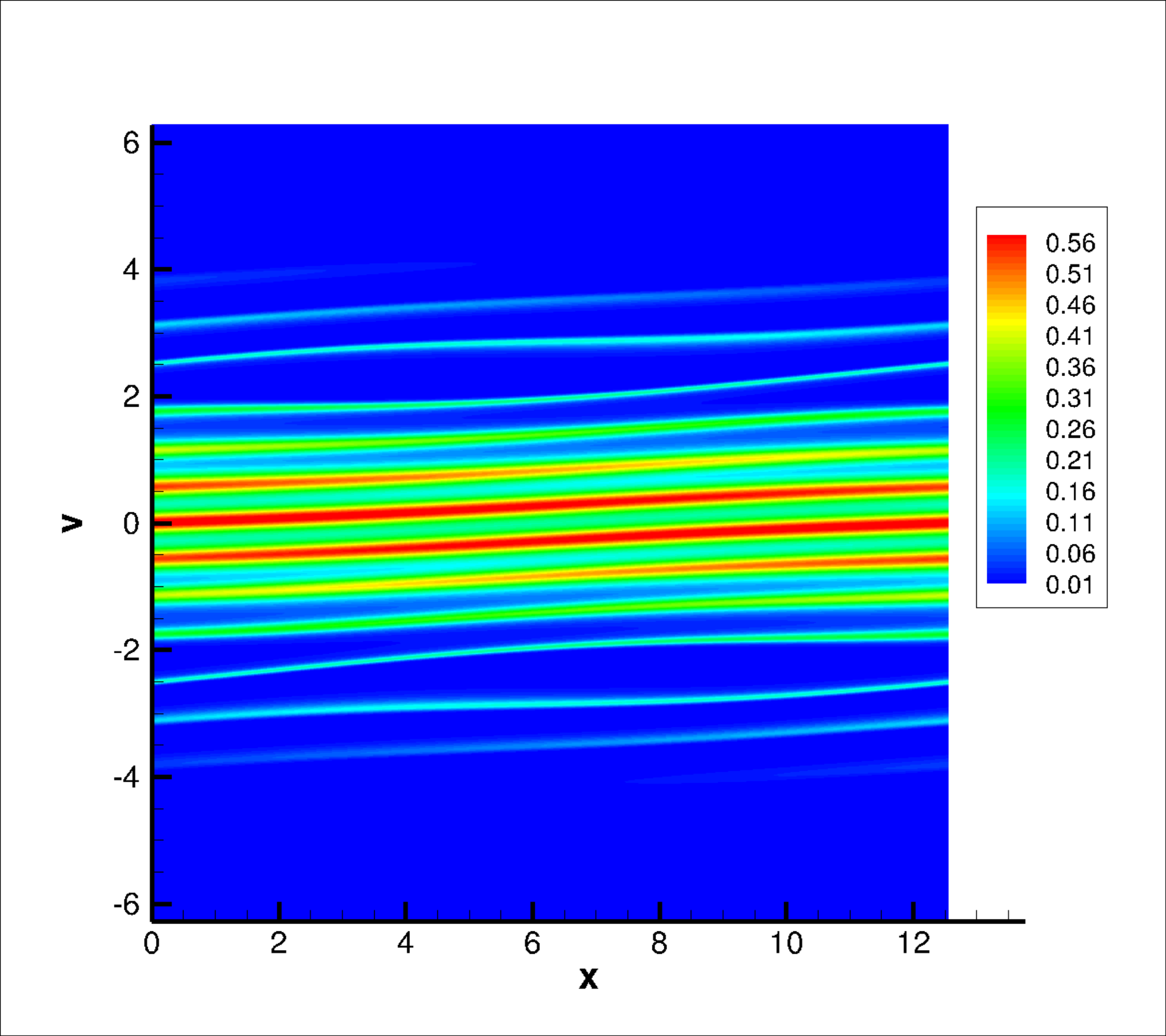}}
\subfigure[$T=30$]{
\includegraphics[width=0.4\textwidth]{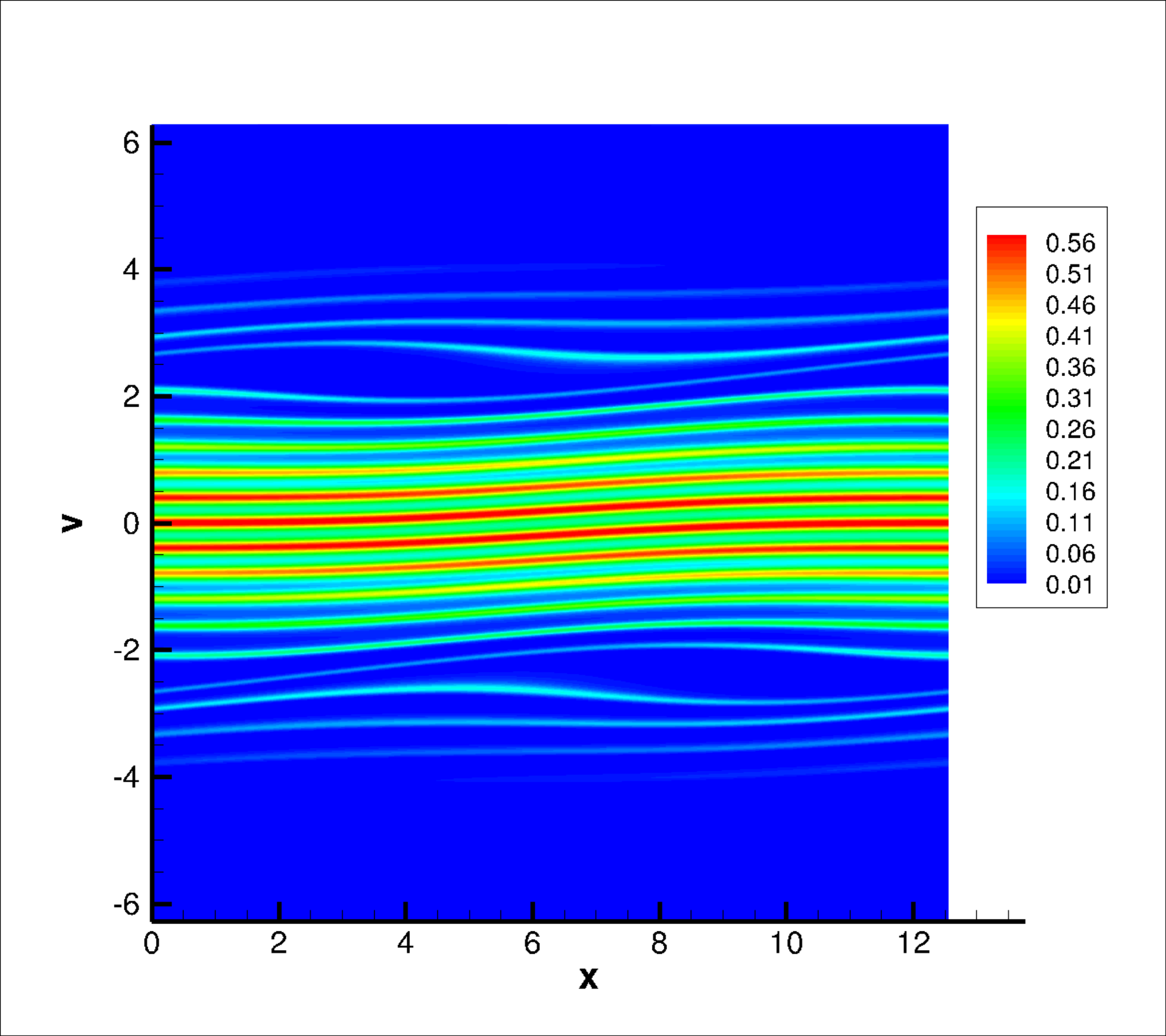}}
\subfigure[$T=40$]{
\includegraphics[width=0.4\textwidth]{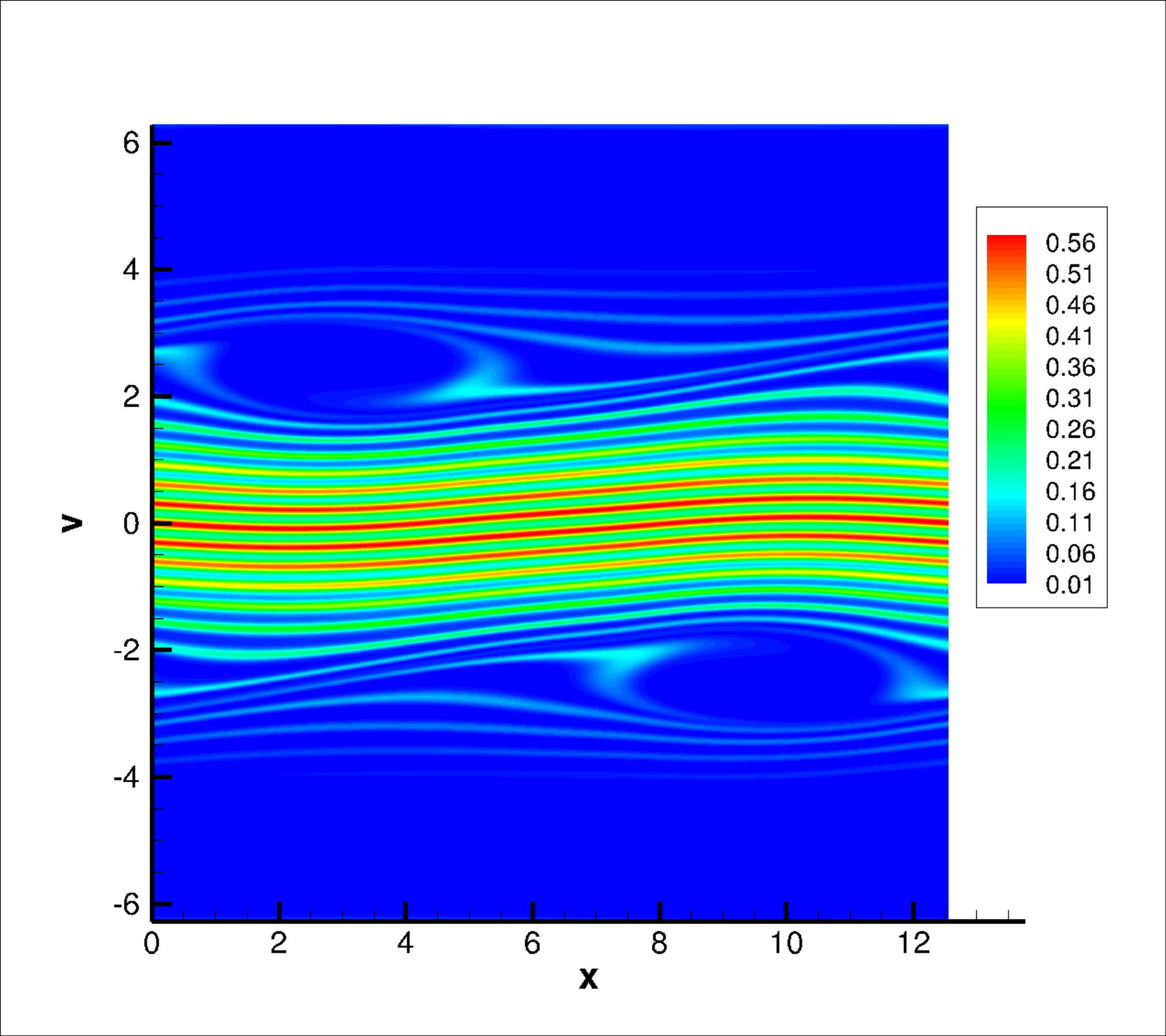}}
\caption{\em Strong Landau Damping. WENO3, with $256\times1024$ grid points.}
\label{Fig3}
\end{figure}

\begin{figure}
\centering
\subfigure[$T=10$]{
\includegraphics[width=0.4\textwidth]{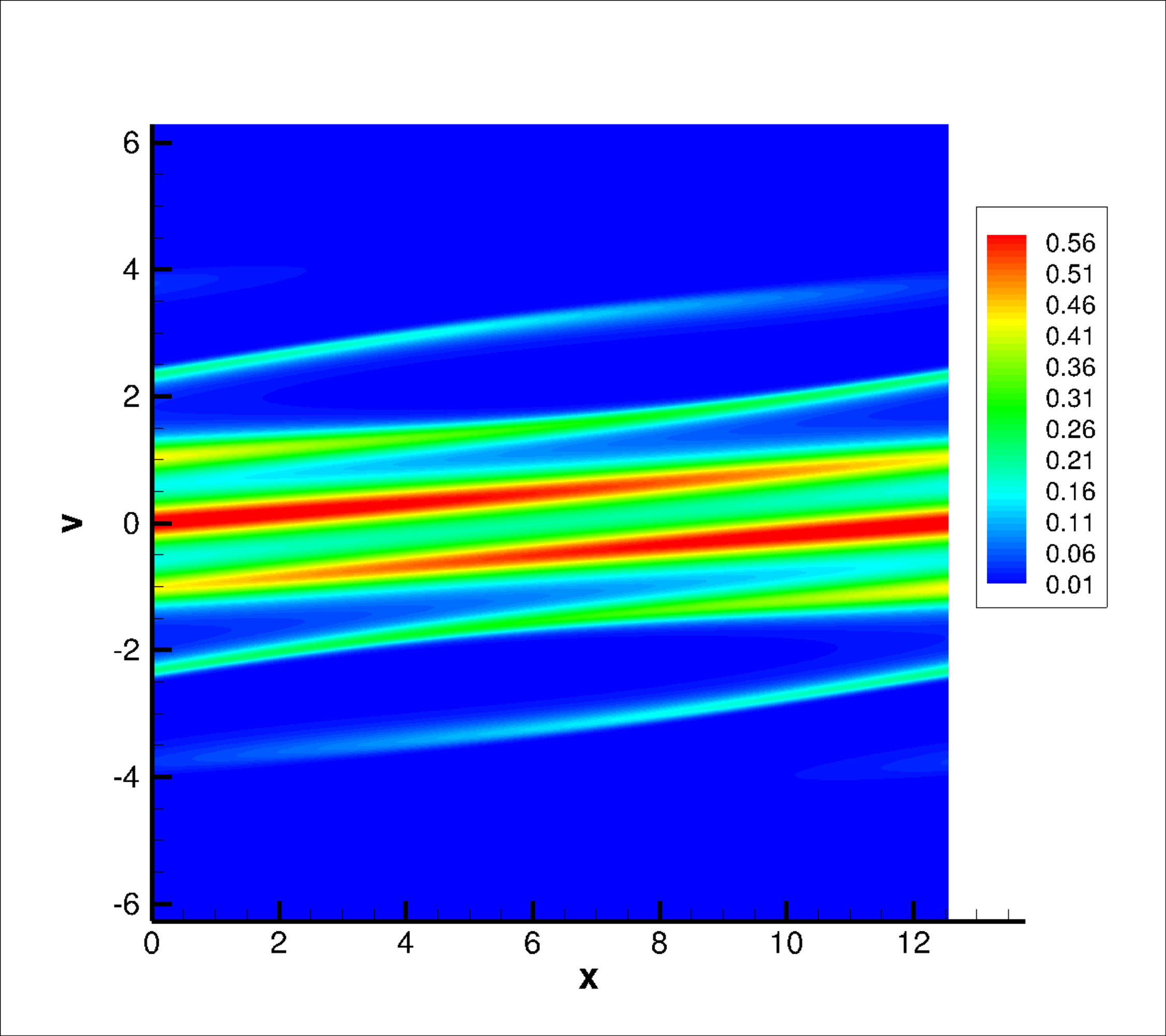}}
\subfigure[$T=20$]{
\includegraphics[width=0.4\textwidth]{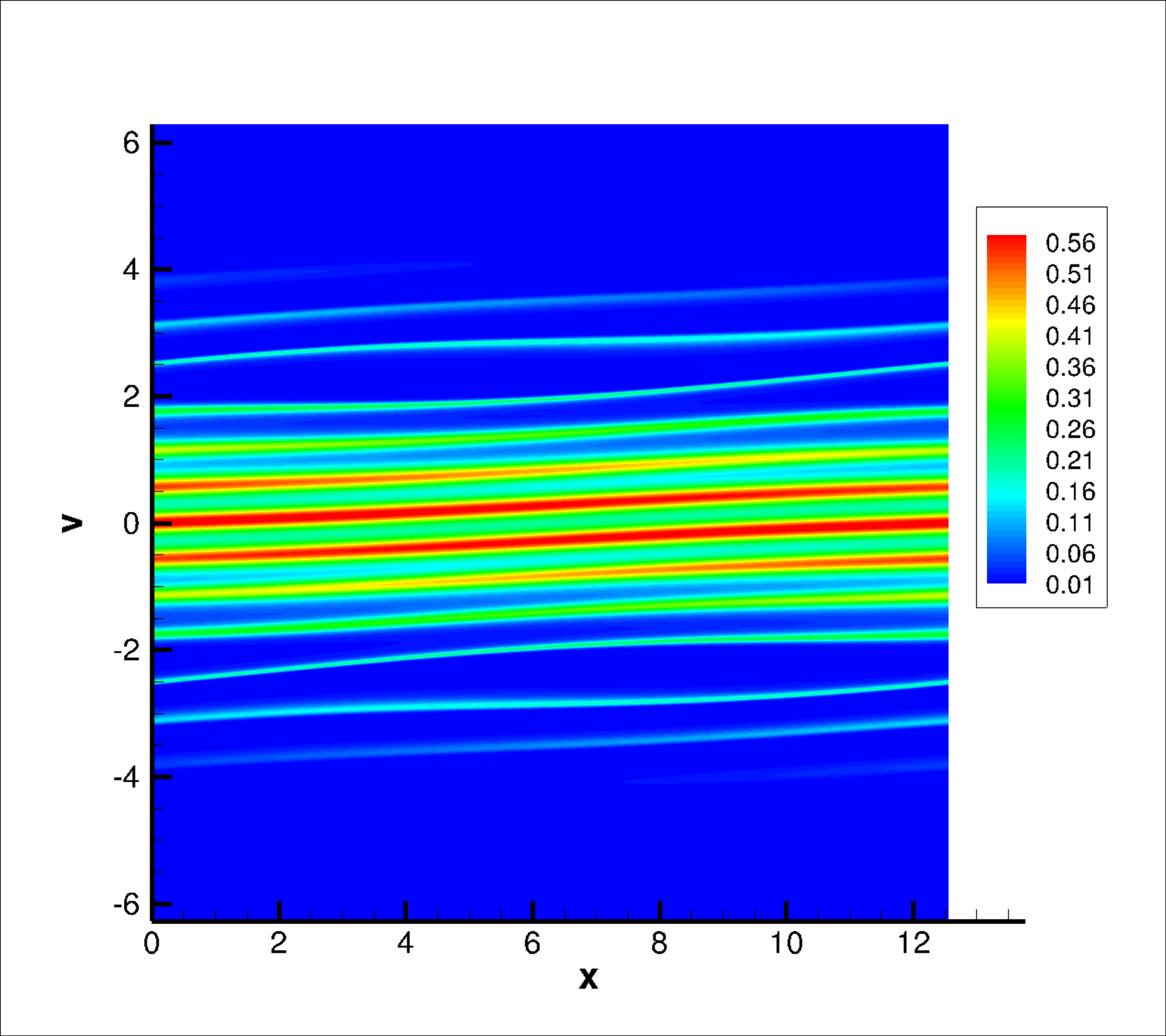}}
\subfigure[$T=30$]{
\includegraphics[width=0.4\textwidth]{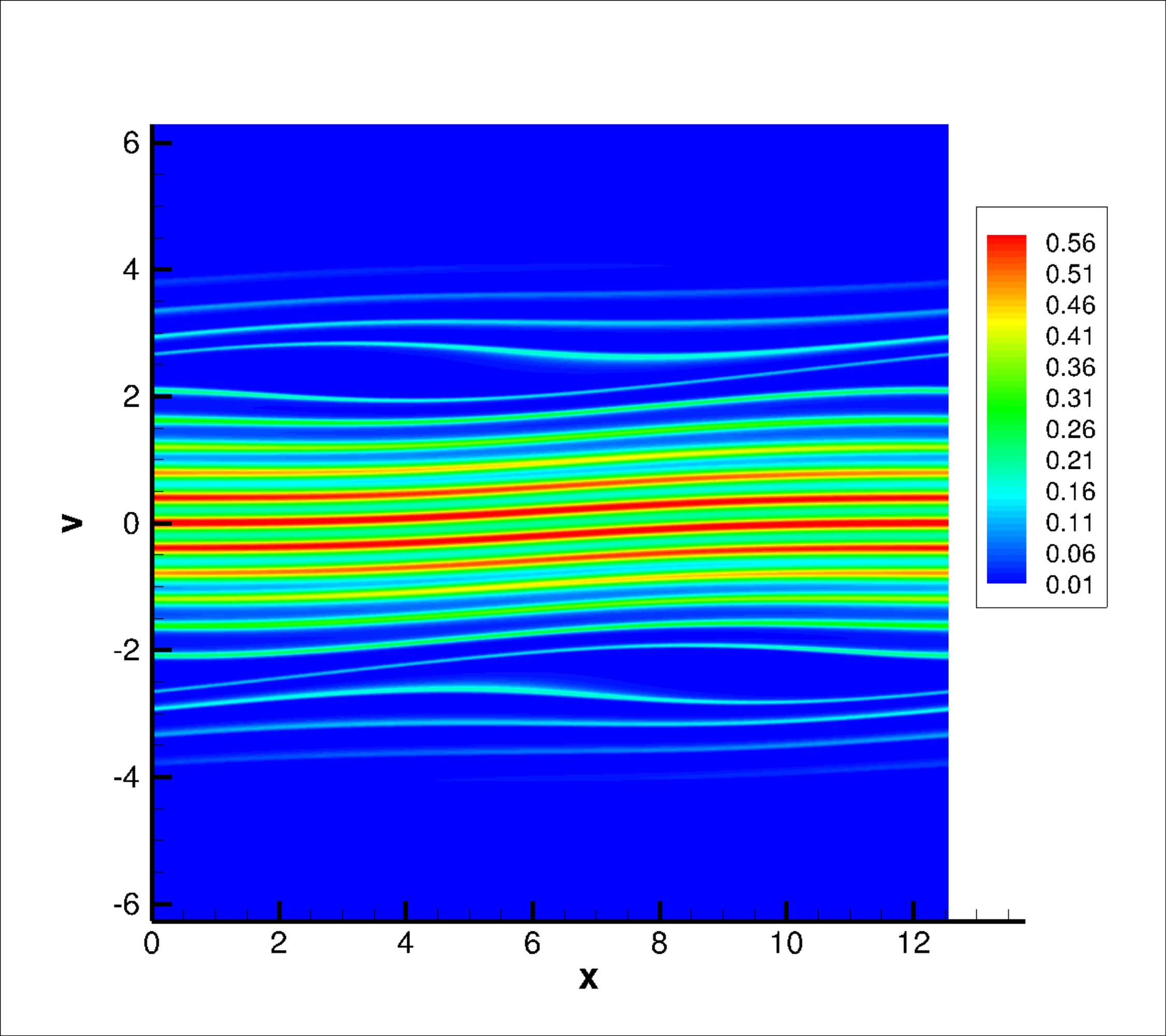}}
\subfigure[$T=40$]{
\includegraphics[width=0.4\textwidth]{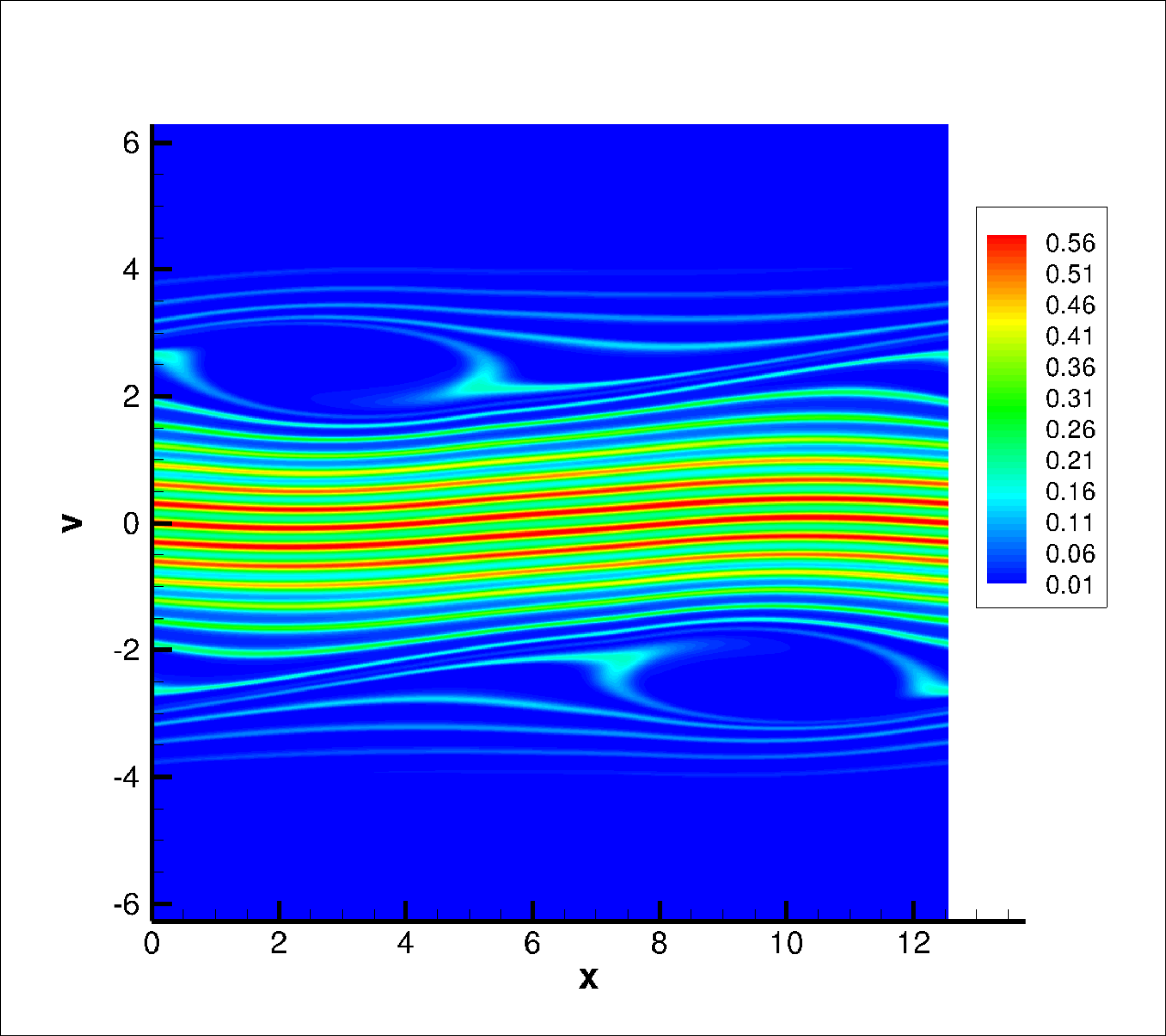}}
\caption{\em Strong Landau Damping. WENO5, with $256\times1024$ grid points.}
\label{Fig4}
\end{figure}

\begin{figure}
\centering
\subfigure[$T=10$]{
\includegraphics[width=0.4\textwidth]{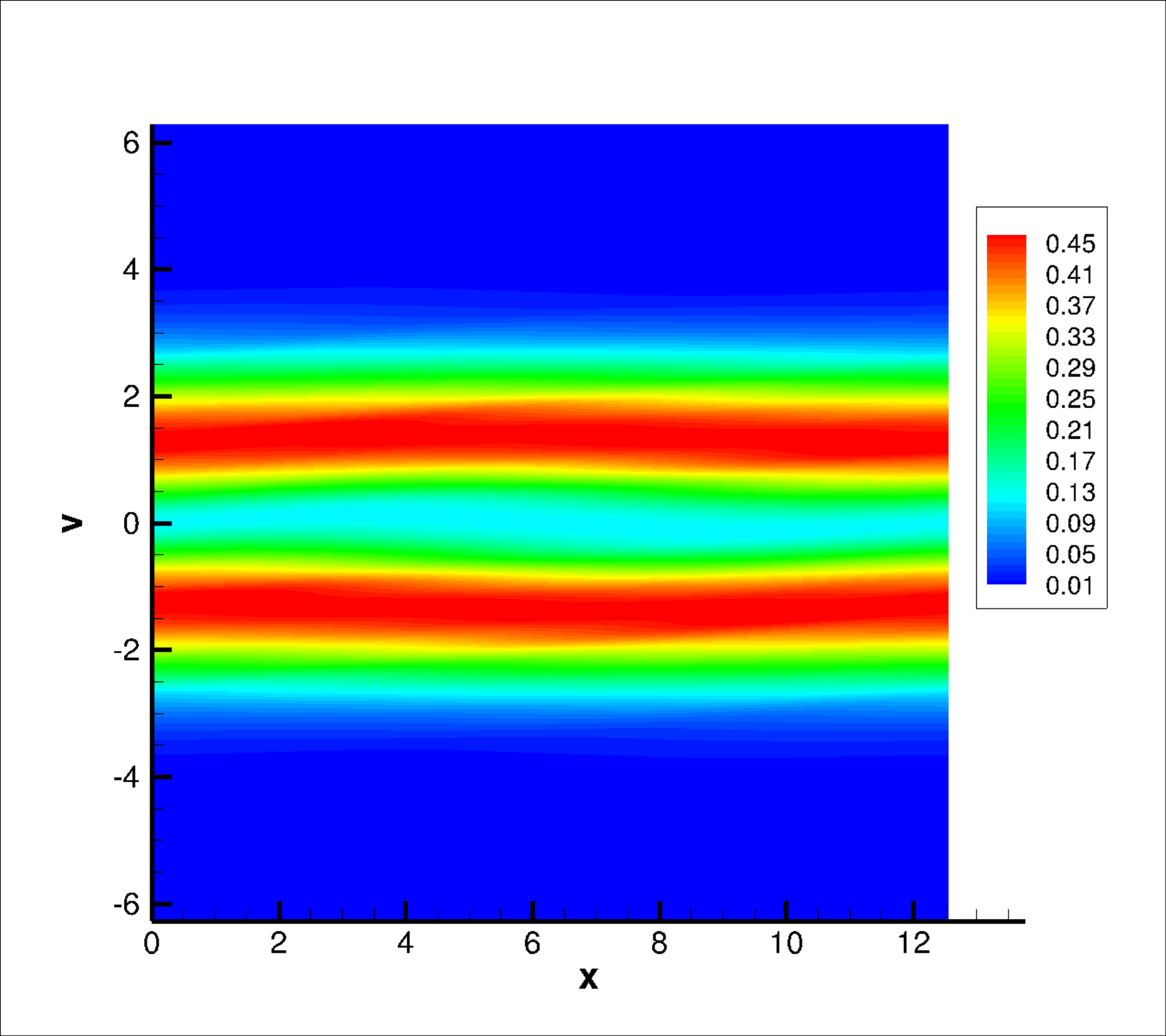}}
\subfigure[$T=20$]{
\includegraphics[width=0.4\textwidth]{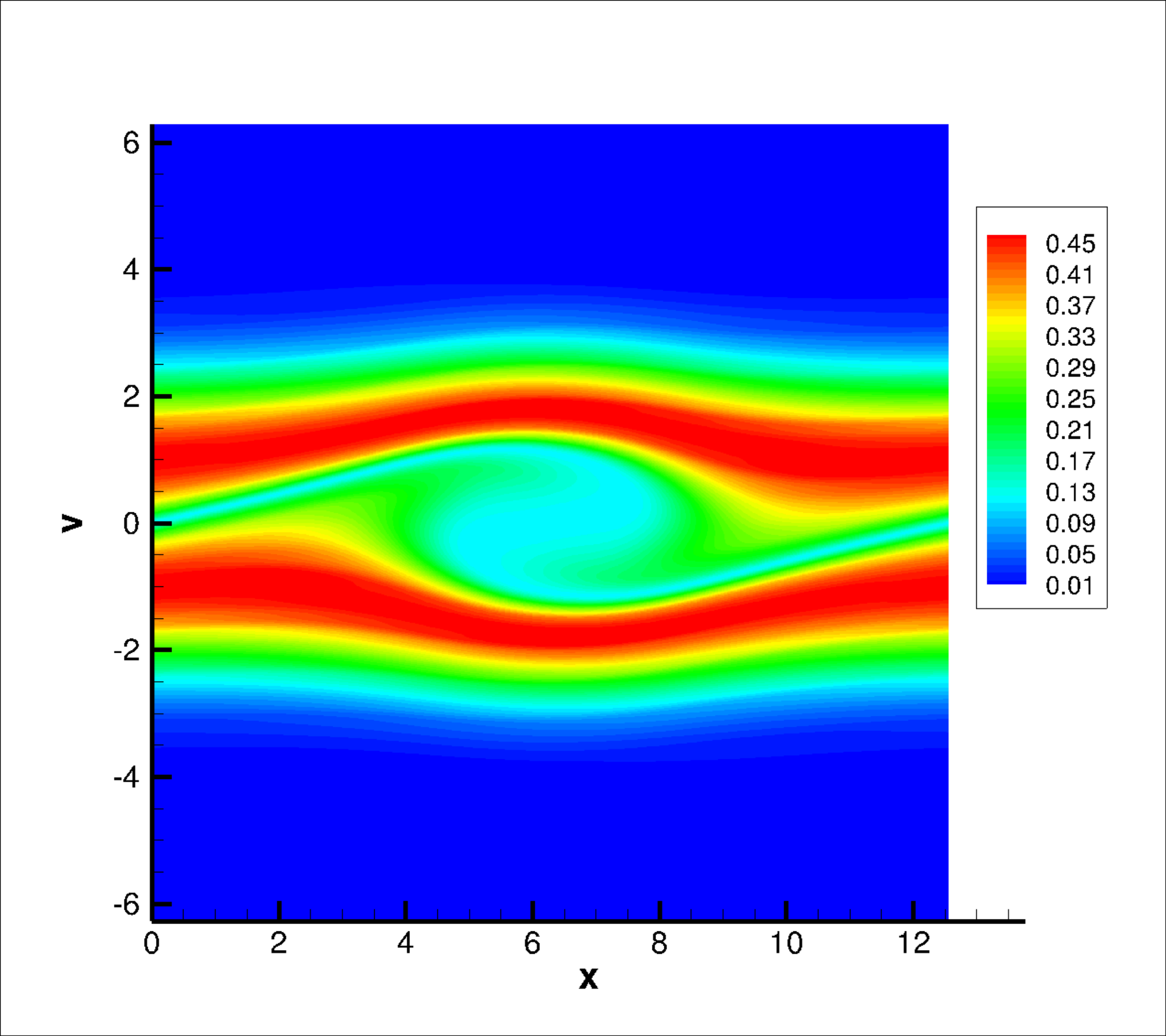}}
\subfigure[$T=30$]{
\includegraphics[width=0.4\textwidth]{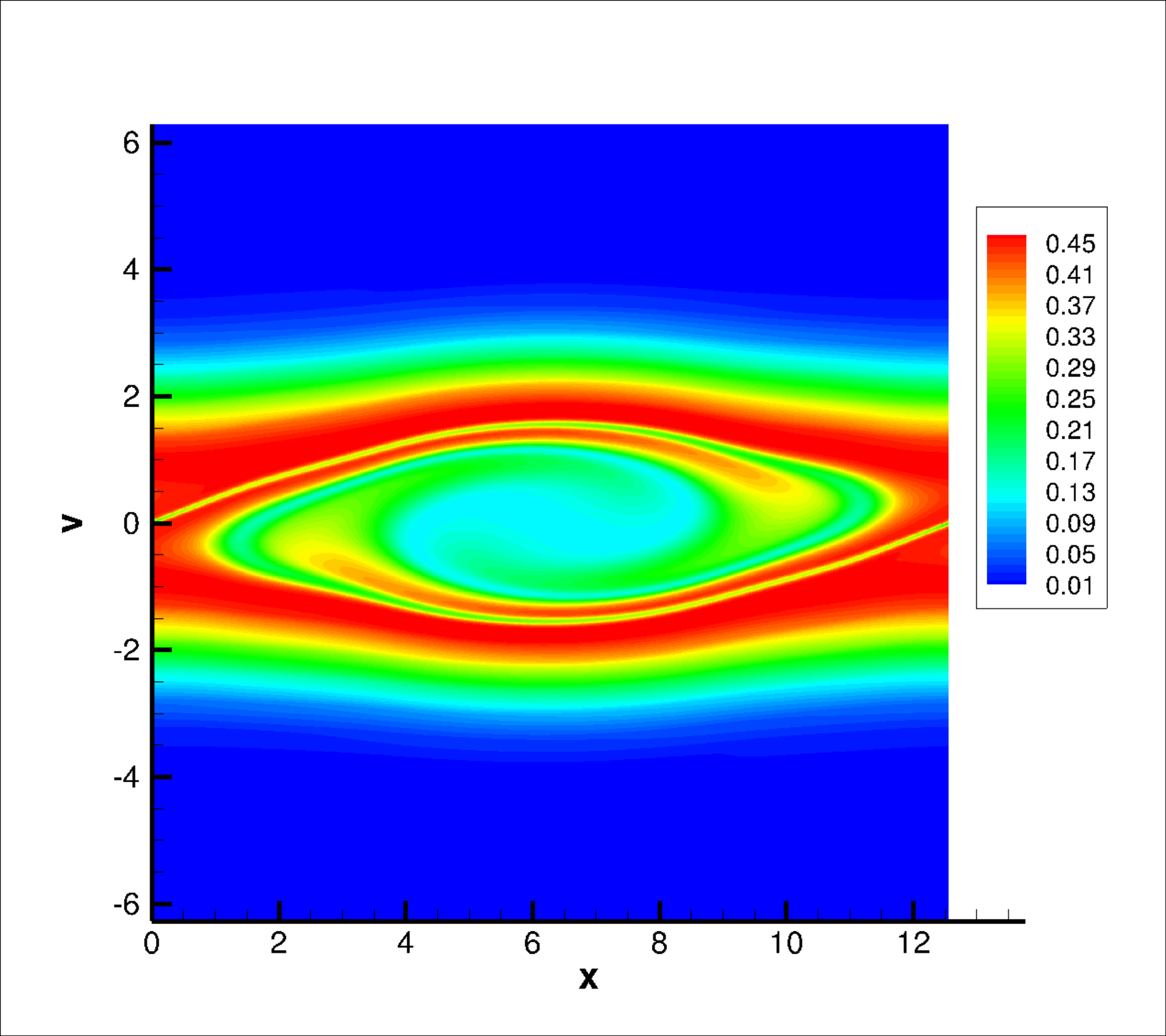}}
\subfigure[$T=40$]{
\includegraphics[width=0.4\textwidth]{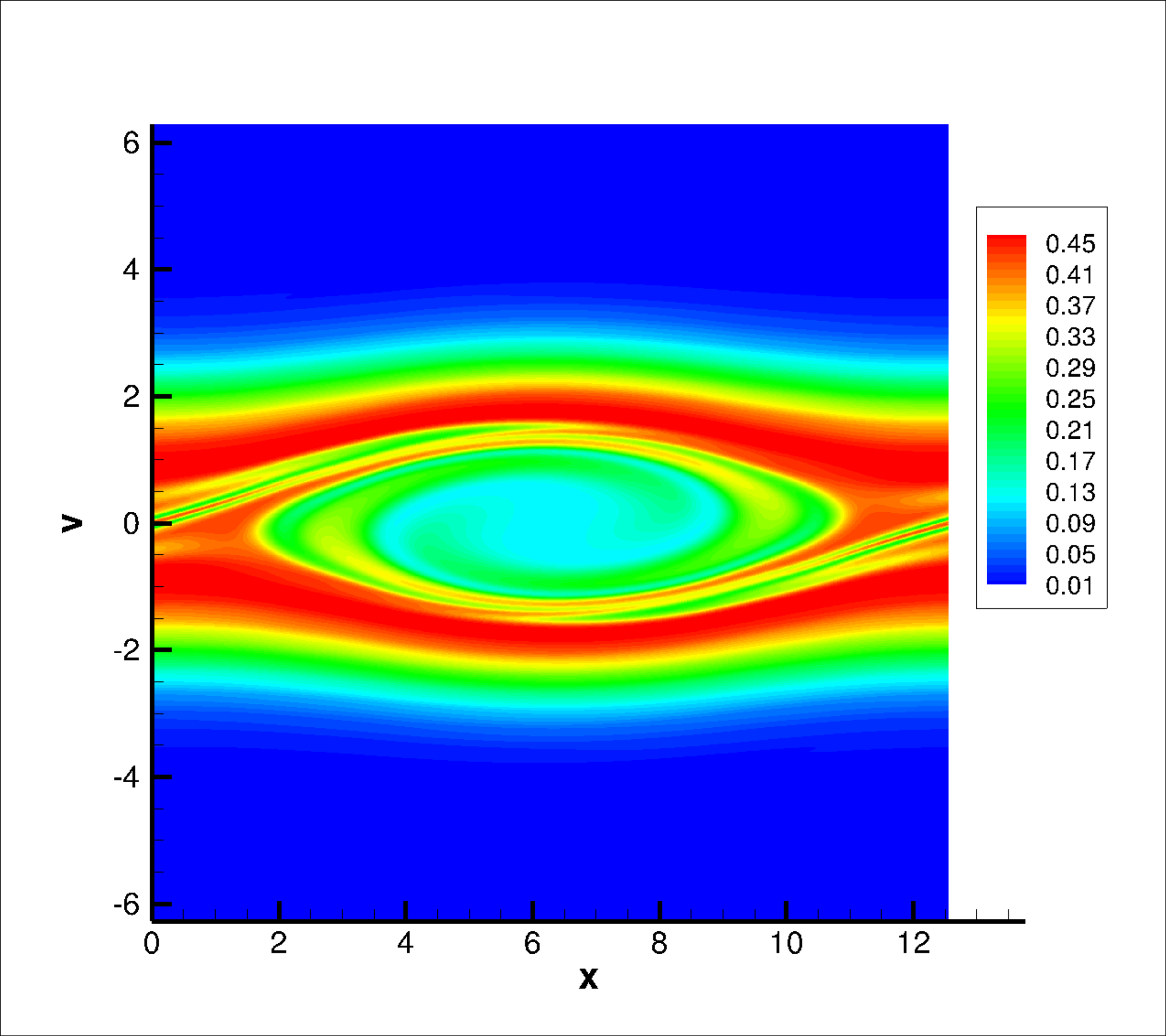}}
\caption{\em Two-stream instability I. WENO3, with $256\times1024$ grid points.}
\label{Fig5}
\end{figure}

\begin{figure}
\centering
\subfigure[$T=10$]{
\includegraphics[width=0.4\textwidth]{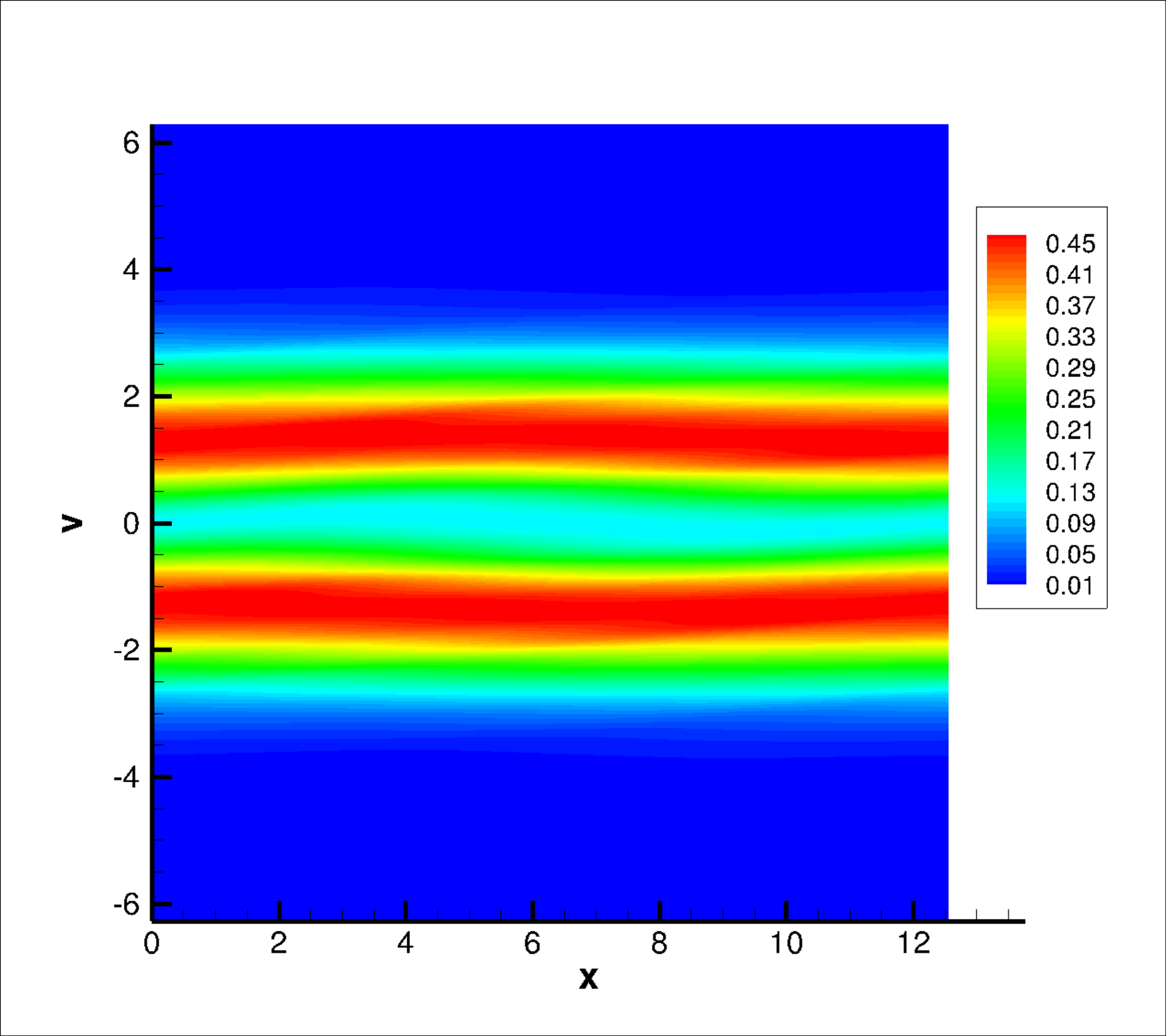}}
\subfigure[$T=20$]{
\includegraphics[width=0.4\textwidth]{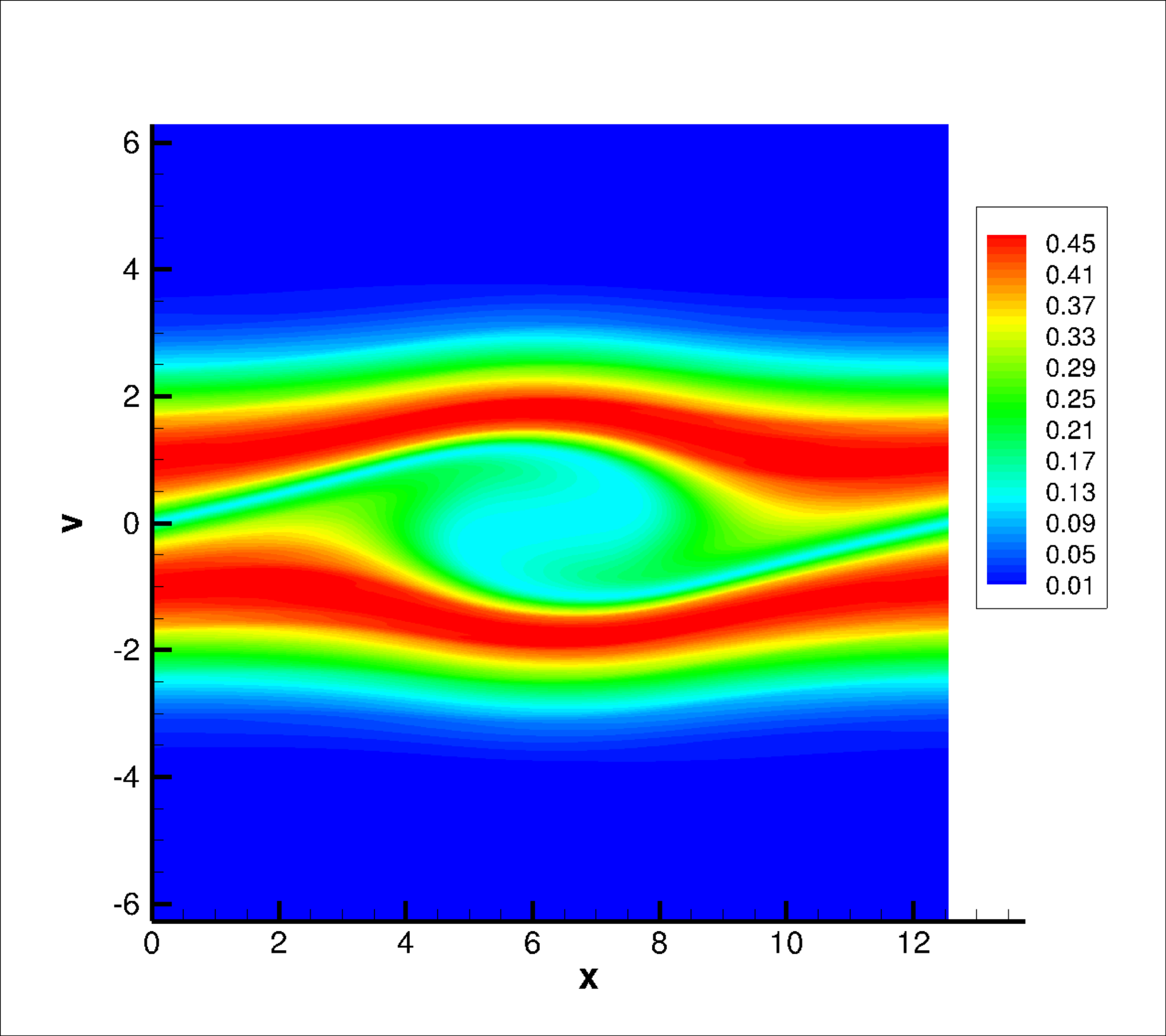}}
\subfigure[$T=30$]{
\includegraphics[width=0.4\textwidth]{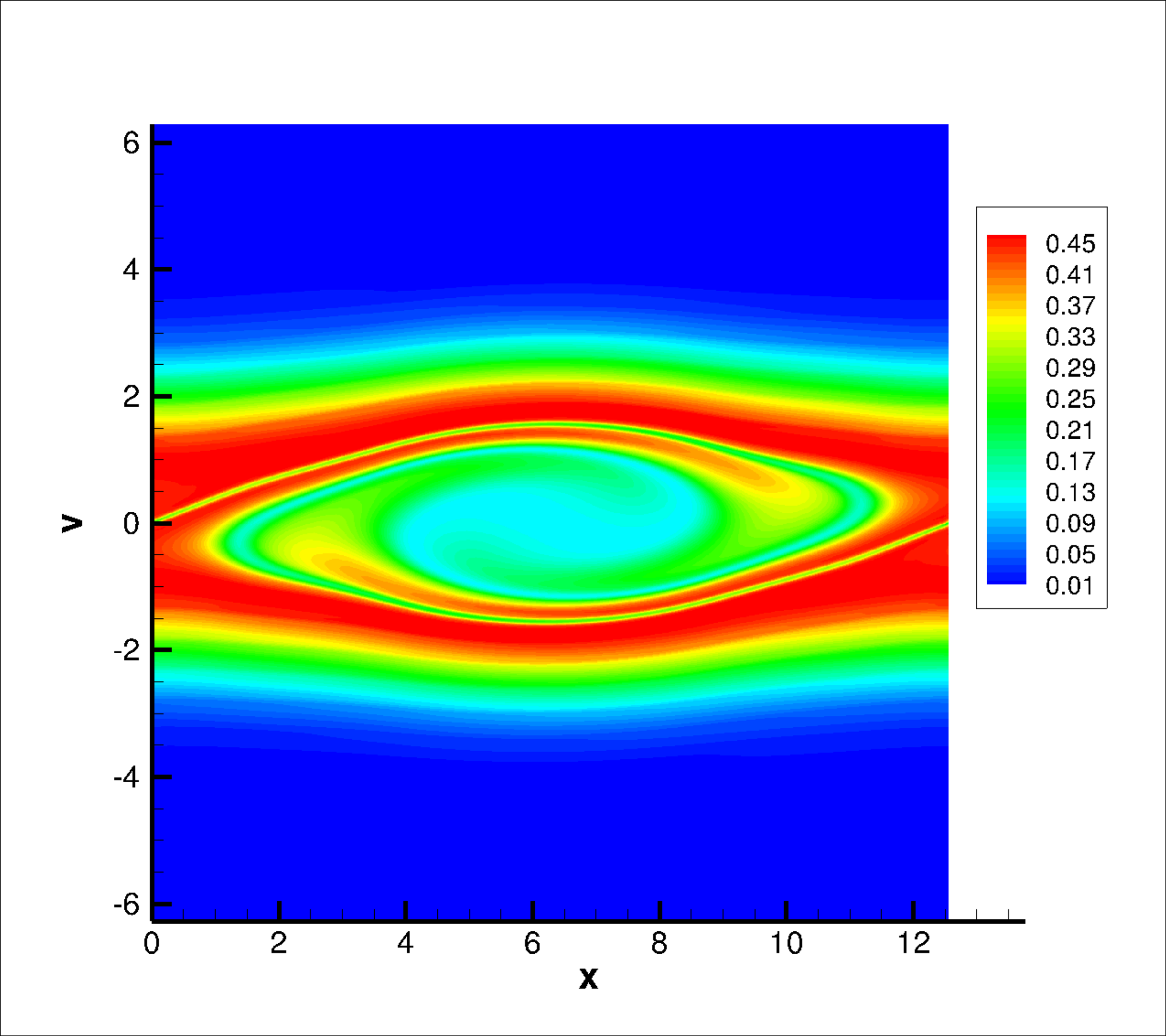}}
\subfigure[$T=40$]{
\includegraphics[width=0.4\textwidth]{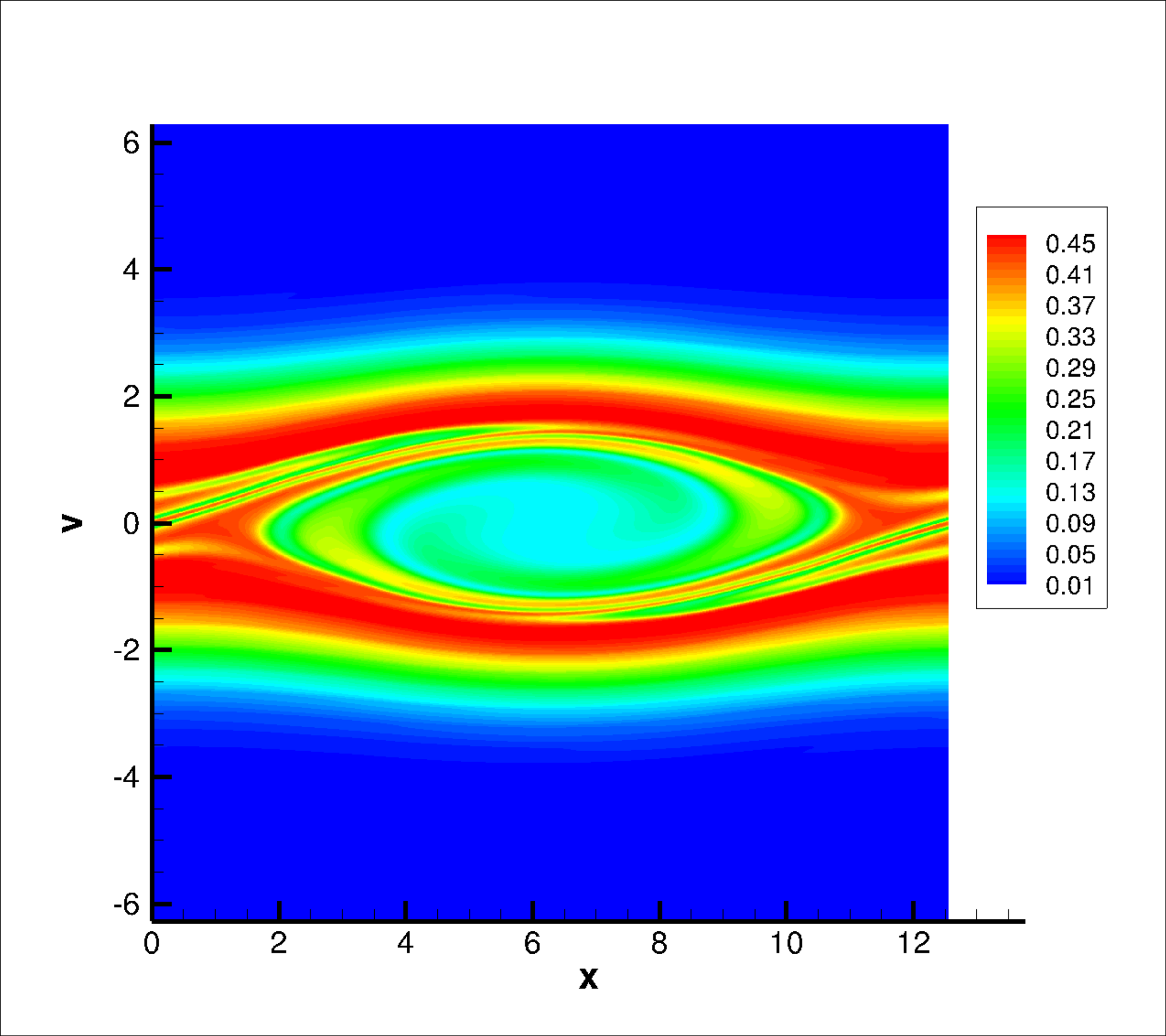}}
\caption{\em Two-stream instability I. WENO5, with $256\times1024$ grid points.}
\label{Fig6}
\end{figure}

\begin{figure}
\centering
\subfigure[$T=10$]{
\includegraphics[width=0.4\textwidth]{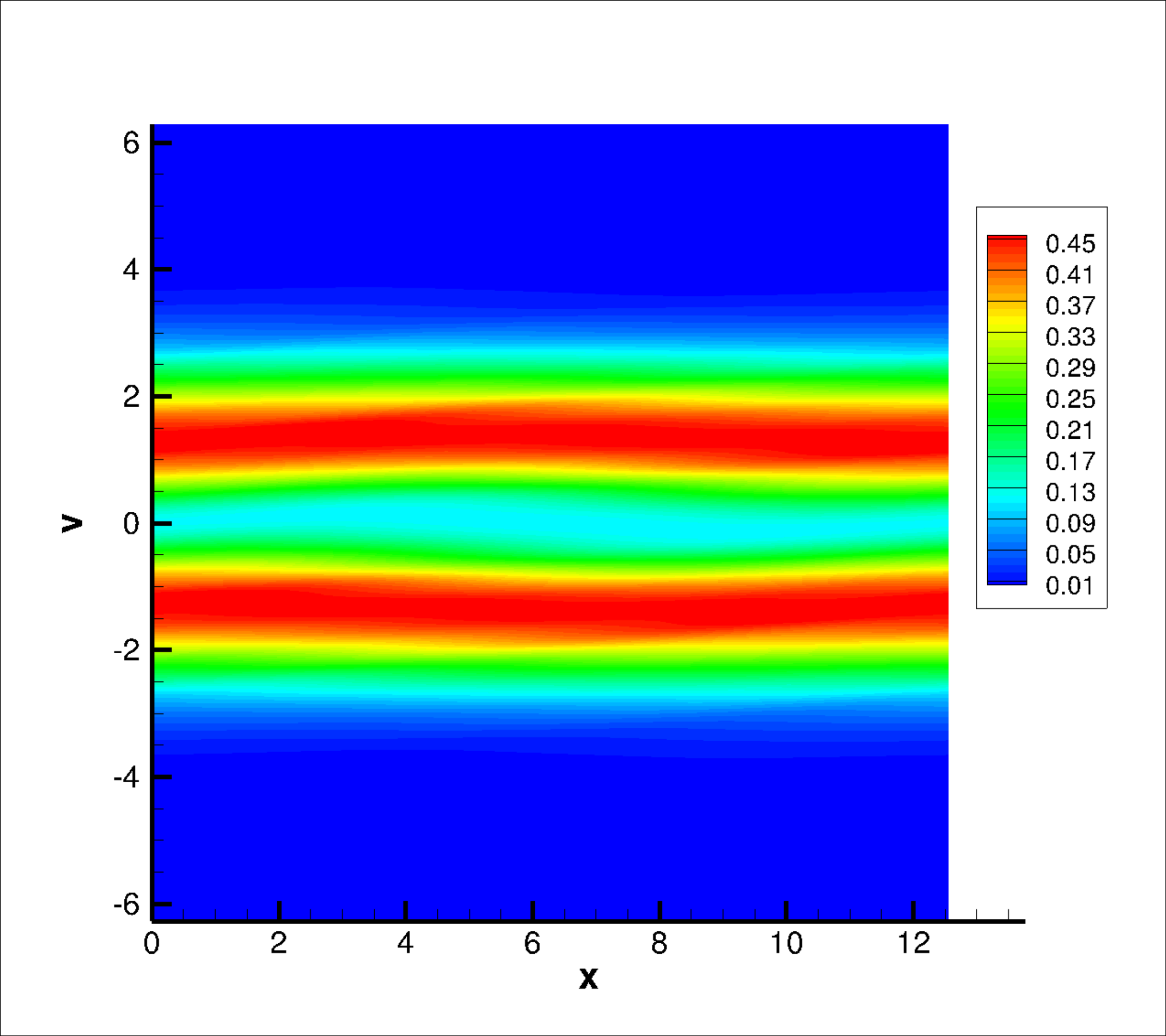}}
\subfigure[$T=20$]{
\includegraphics[width=0.4\textwidth]{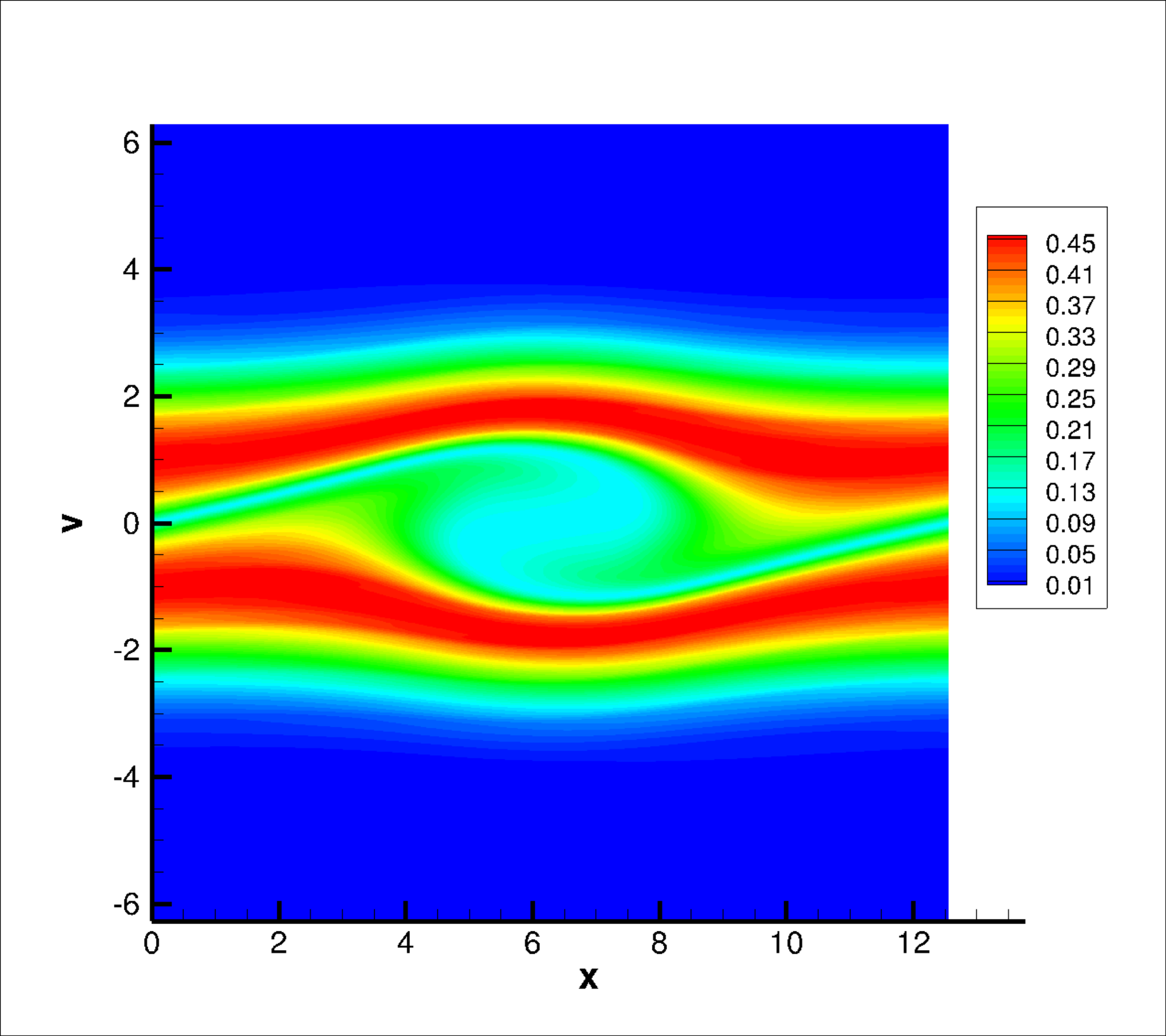}}
\subfigure[$T=30$]{
\includegraphics[width=0.4\textwidth]{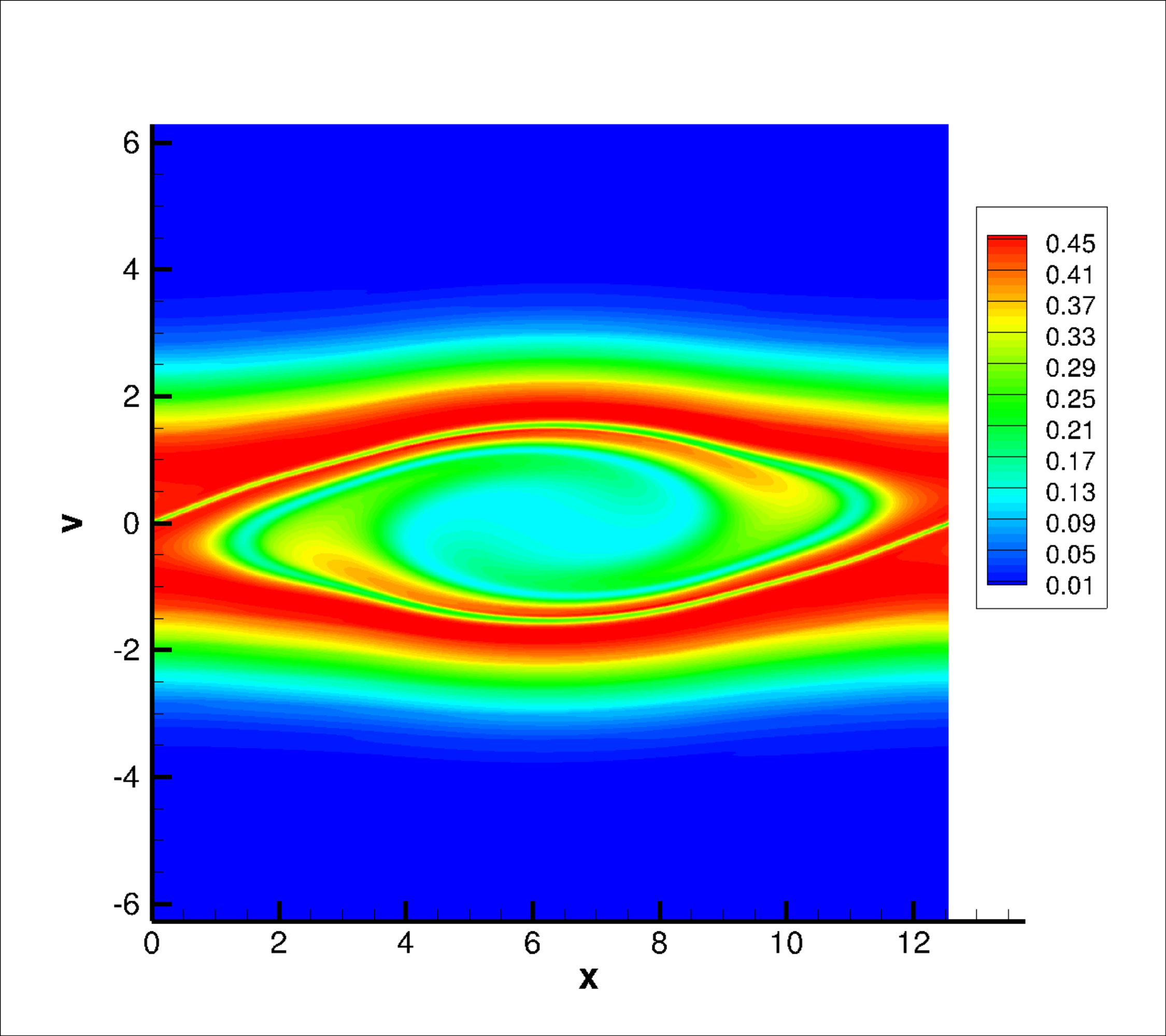}}
\subfigure[$T=40$]{
\includegraphics[width=0.4\textwidth]{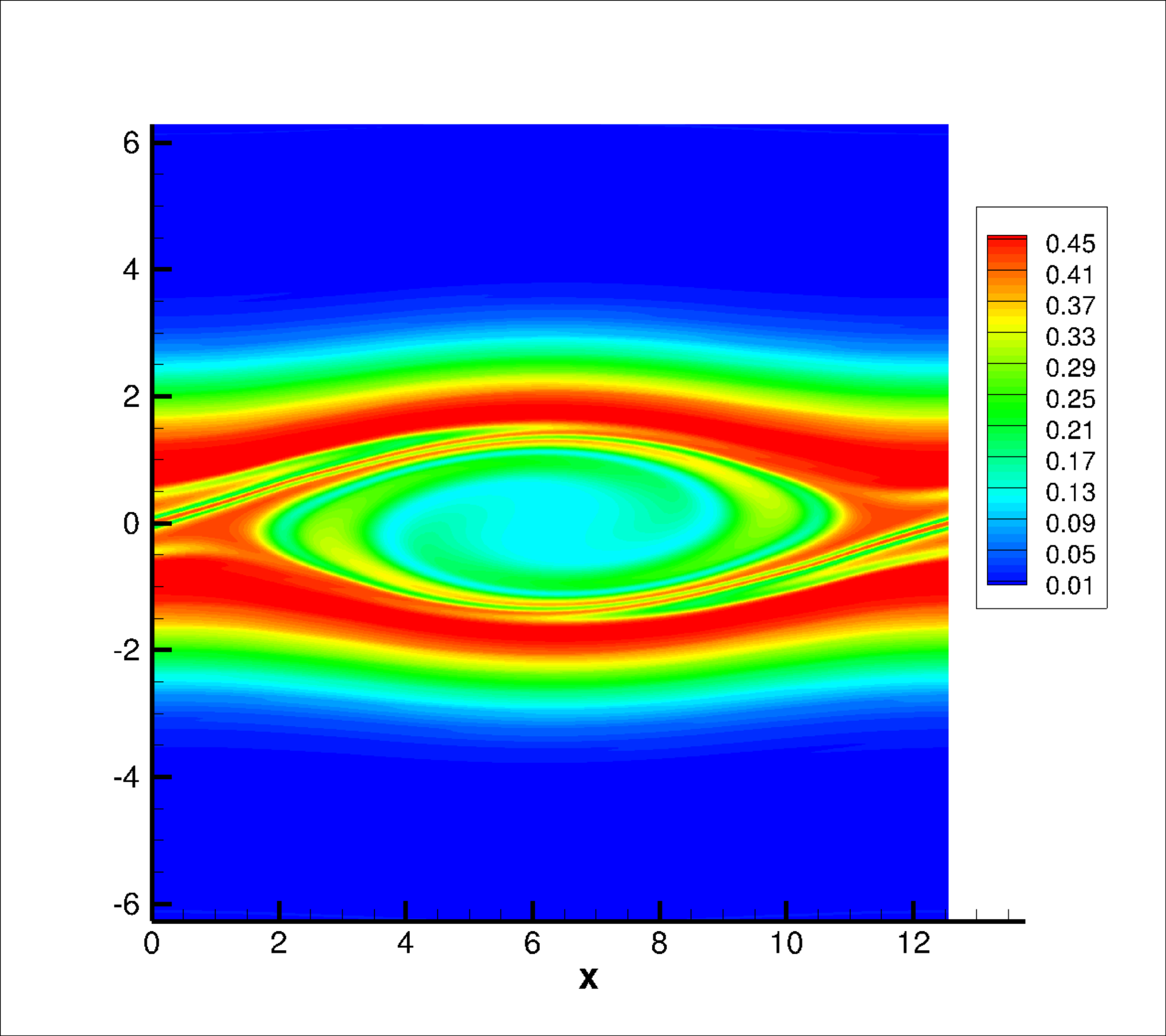}}
\caption{\em Two-stream instability I. Second order centel difference scheme with WENO5. $256\times1024$ grid points.}
\label{Fig17}
\end{figure}

\begin{figure}
\centering
\subfigure[$T=10$]{
\includegraphics[width=0.4\textwidth]{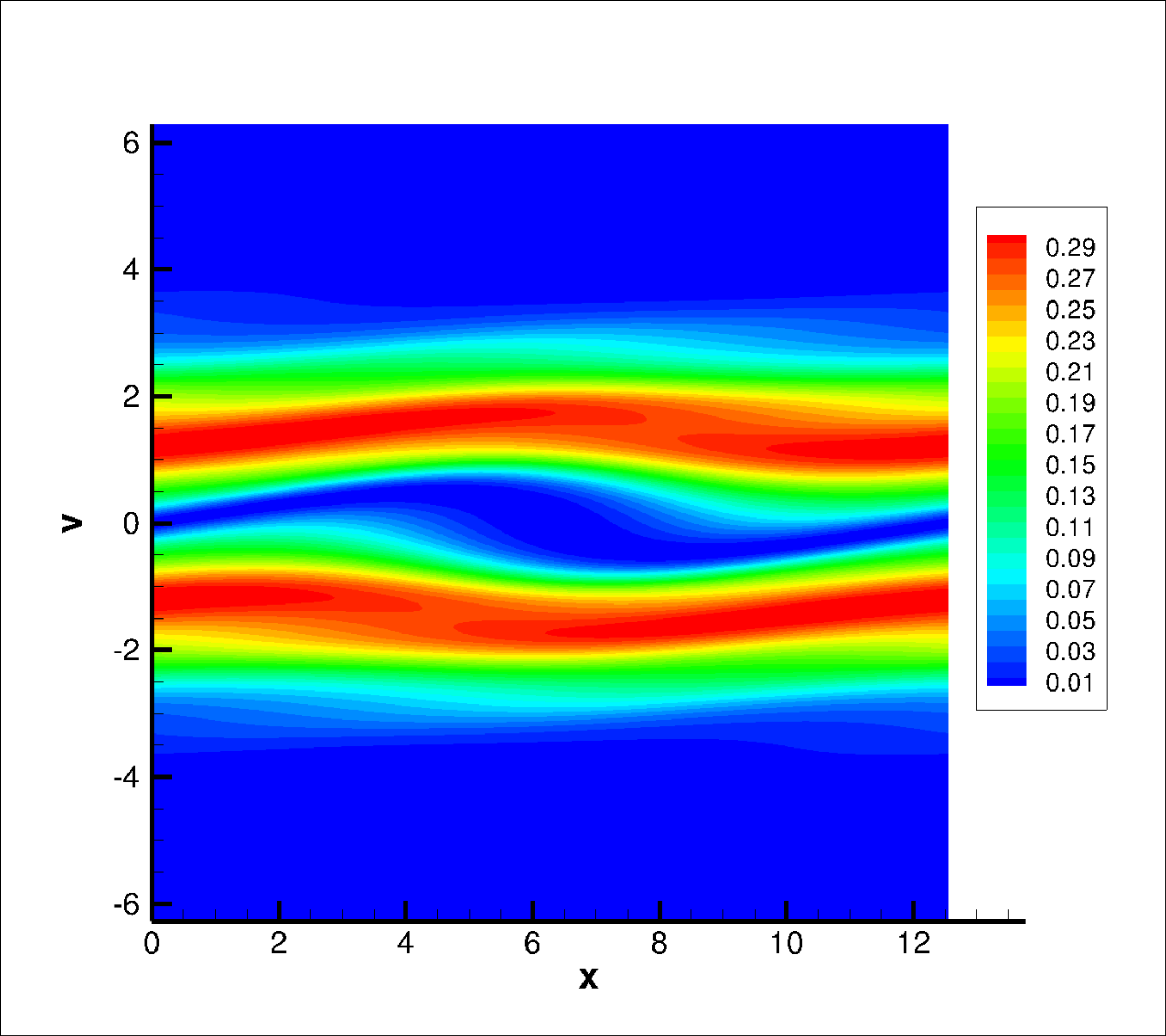}}
\subfigure[$T=20$]{
\includegraphics[width=0.4\textwidth]{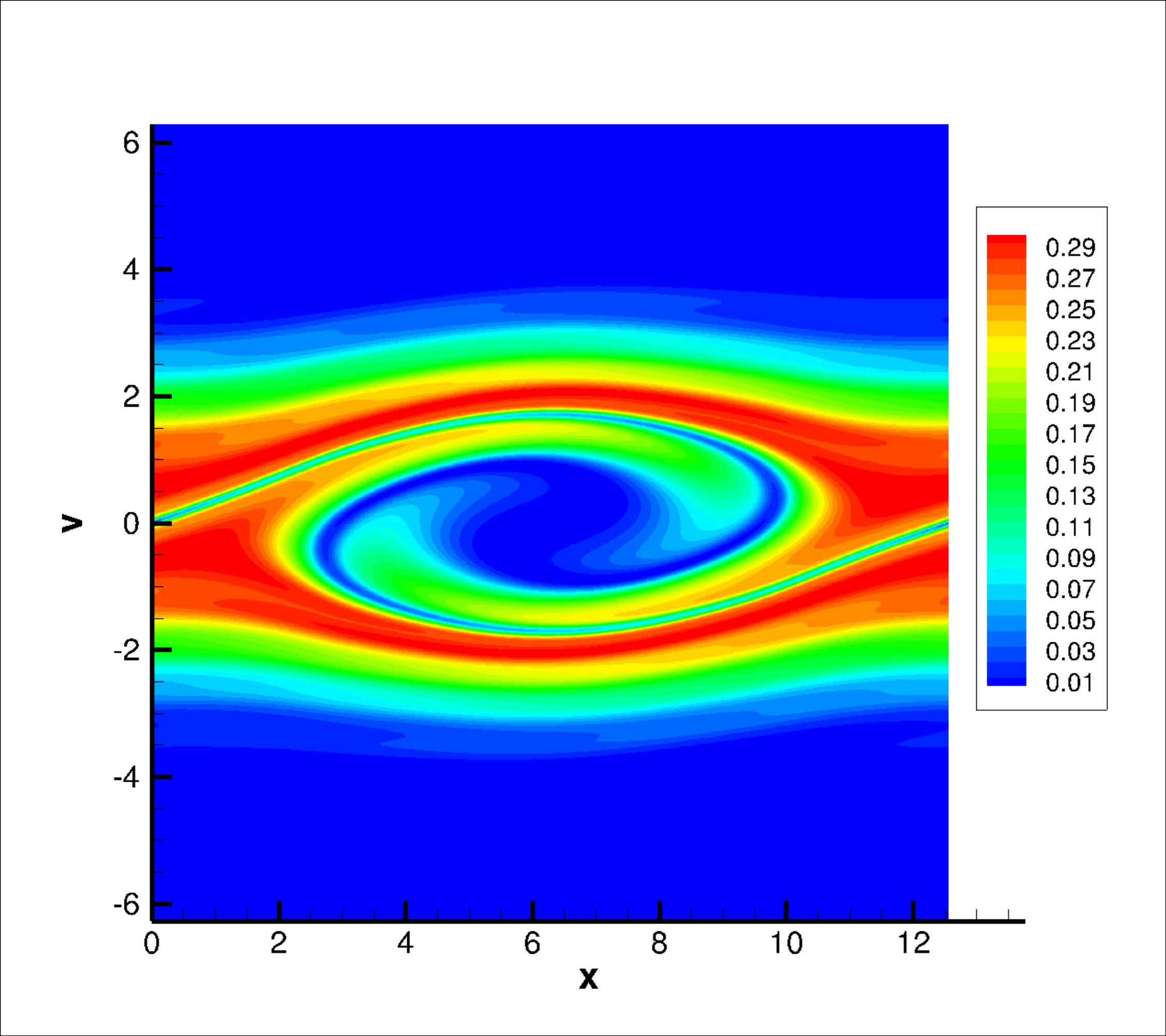}}
\subfigure[$T=30$]{
\includegraphics[width=0.4\textwidth]{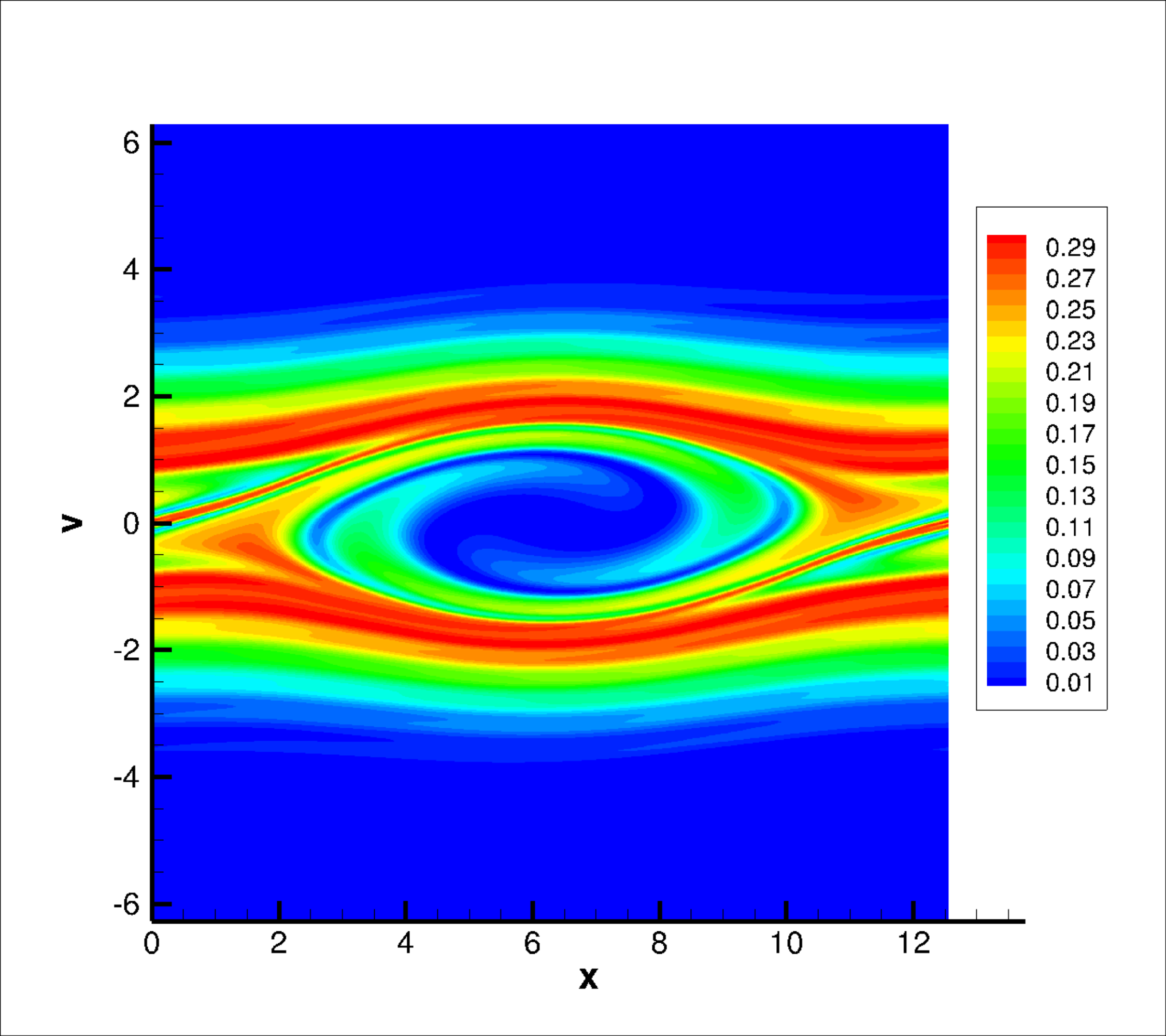}}
\subfigure[$T=40$]{
\includegraphics[width=0.4\textwidth]{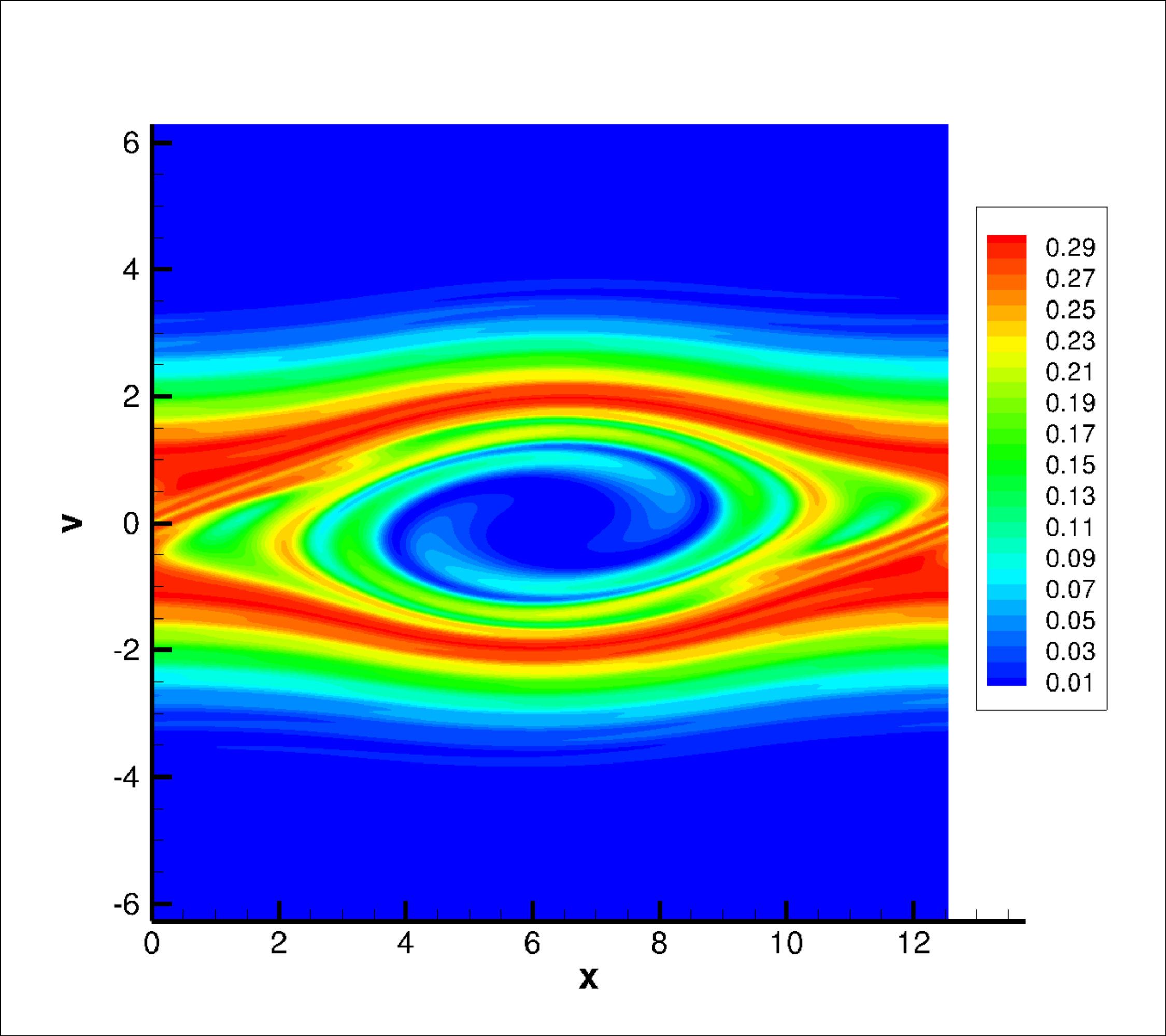}}
\caption{\em Two-stream instability II. WENO3, with $256\times1024$ grid points.}
\label{Fig7}
\end{figure}

\begin{figure}
\centering
\subfigure[$T=10$]{
\includegraphics[width=0.4\textwidth]{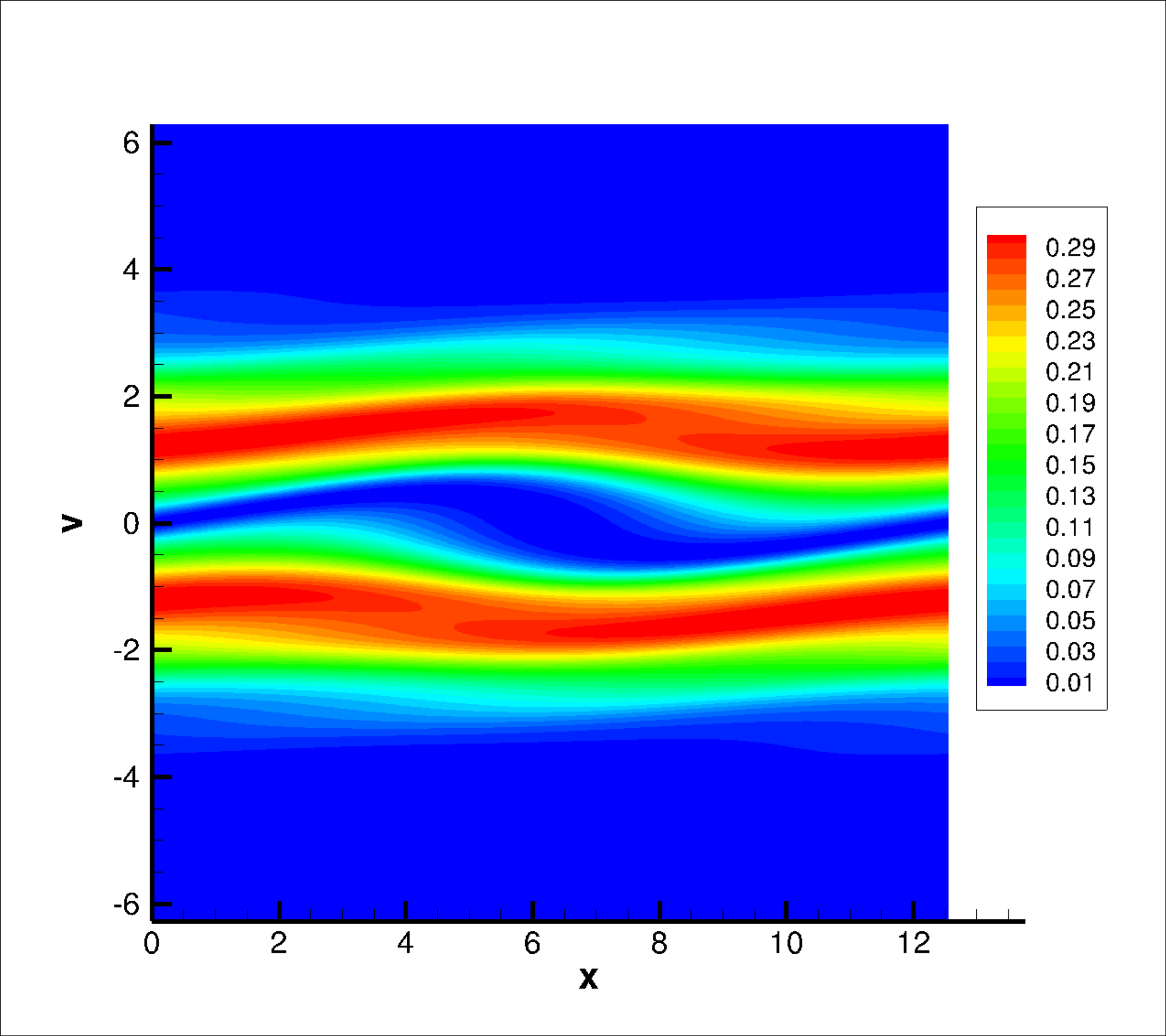}}
\subfigure[$T=20$]{
\includegraphics[width=0.4\textwidth]{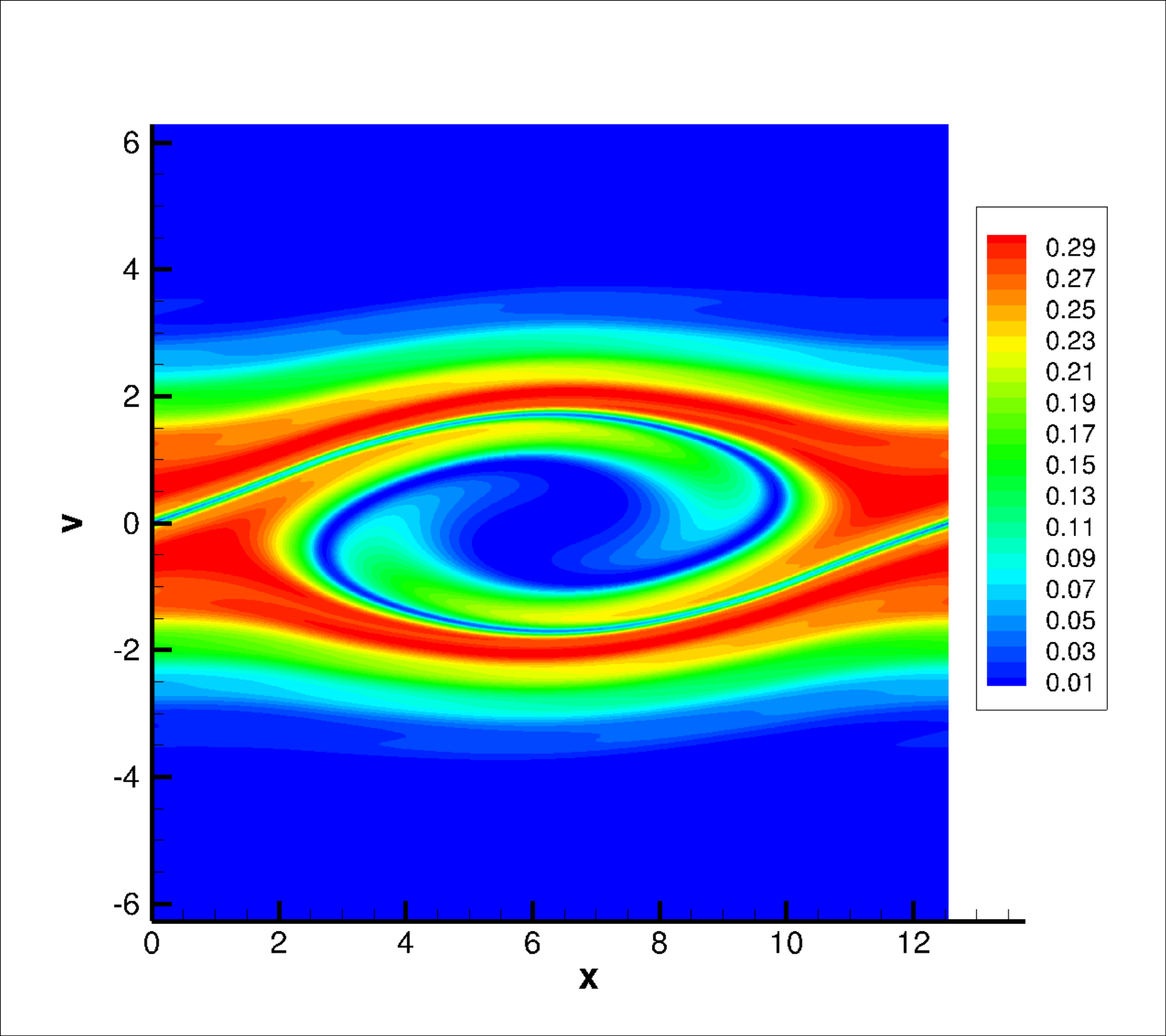}}
\subfigure[$T=30$]{
\includegraphics[width=0.4\textwidth]{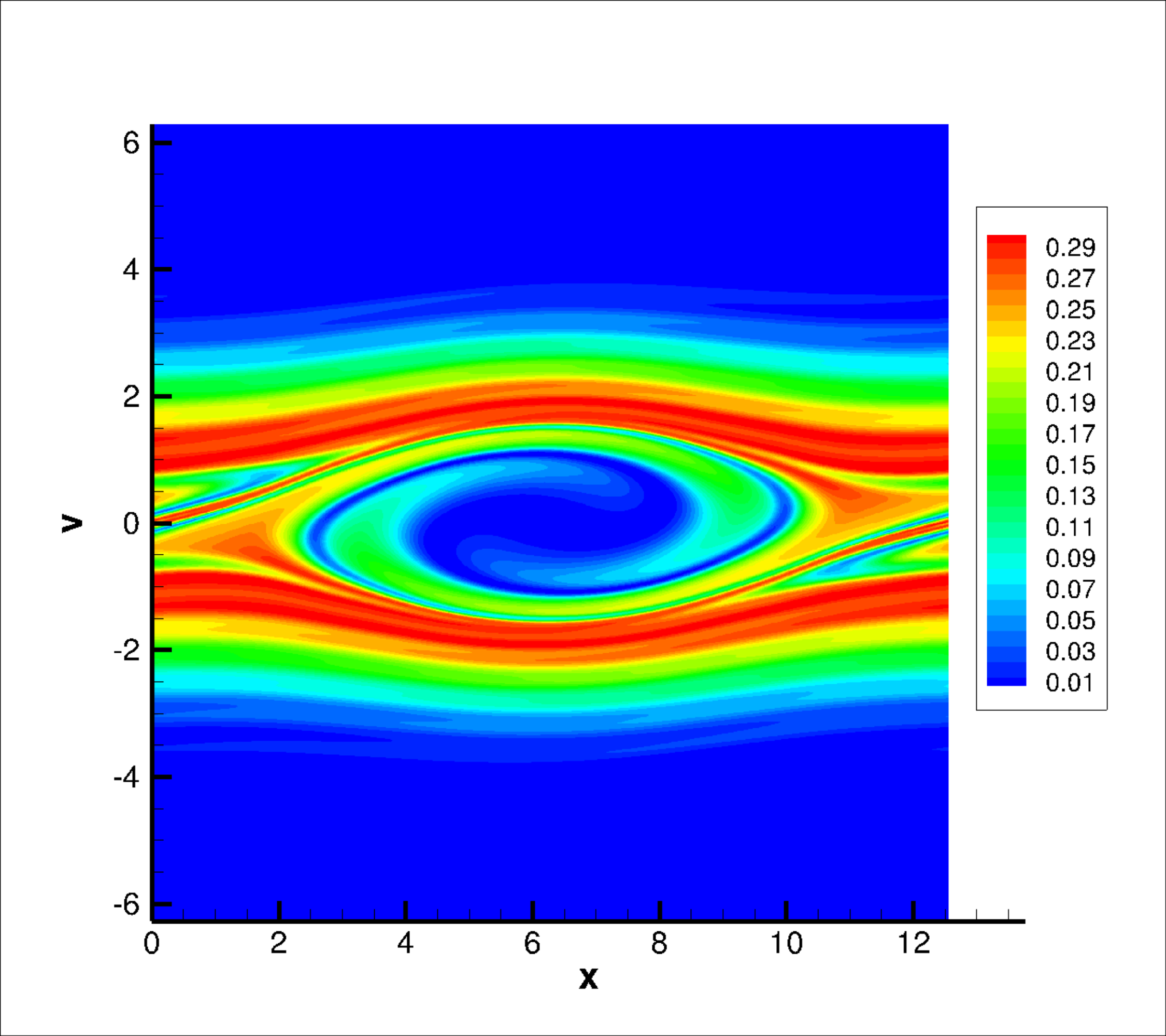}}
\subfigure[$T=40$]{
\includegraphics[width=0.4\textwidth]{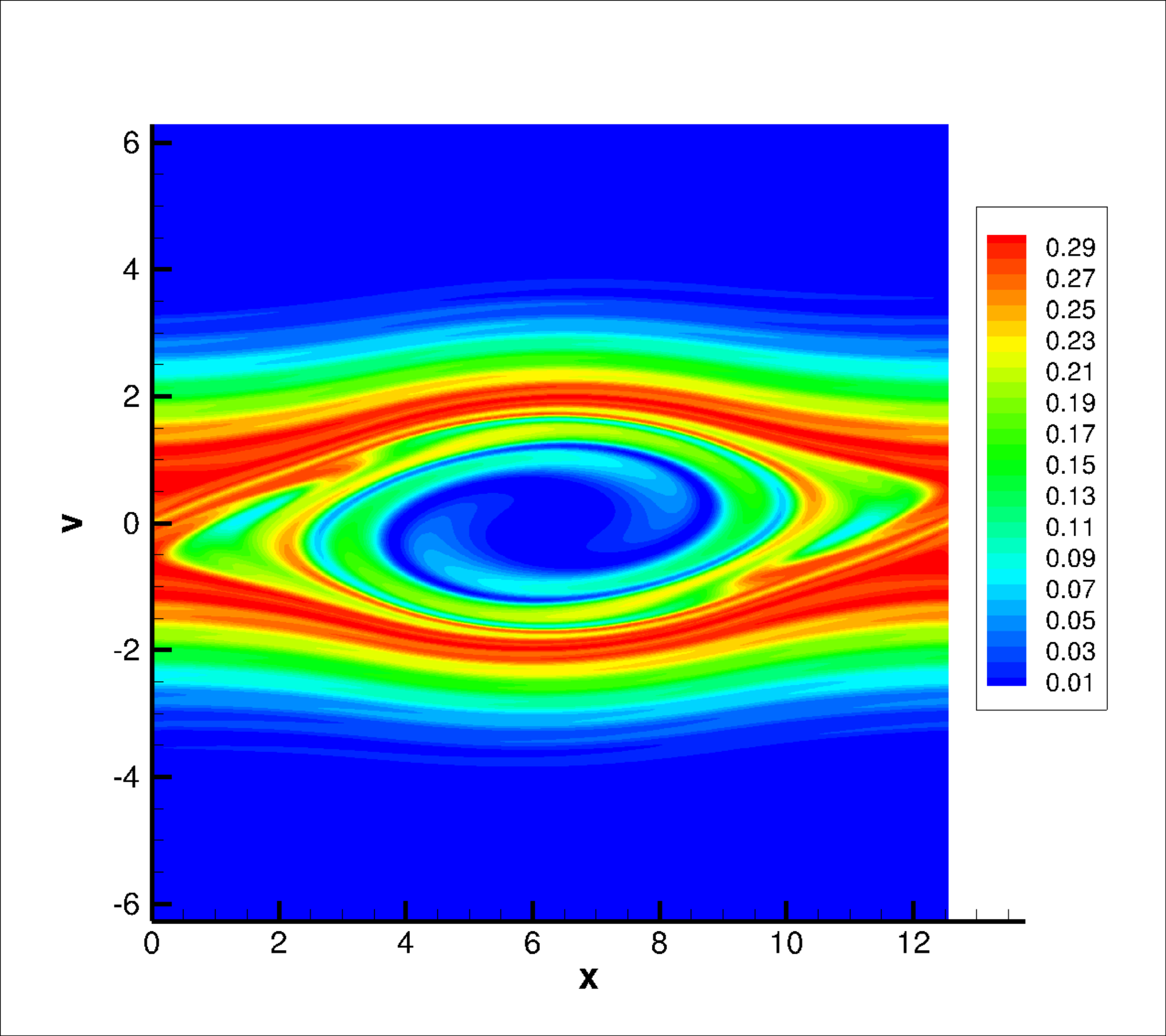}}
\caption{\em Two-stream instability II. WENO5, with $256\times1024$ grid points.}
\label{Fig8}
\end{figure}

\begin{figure}
\centering
\subfigure[$T=10$]{
\includegraphics[width=0.4\textwidth]{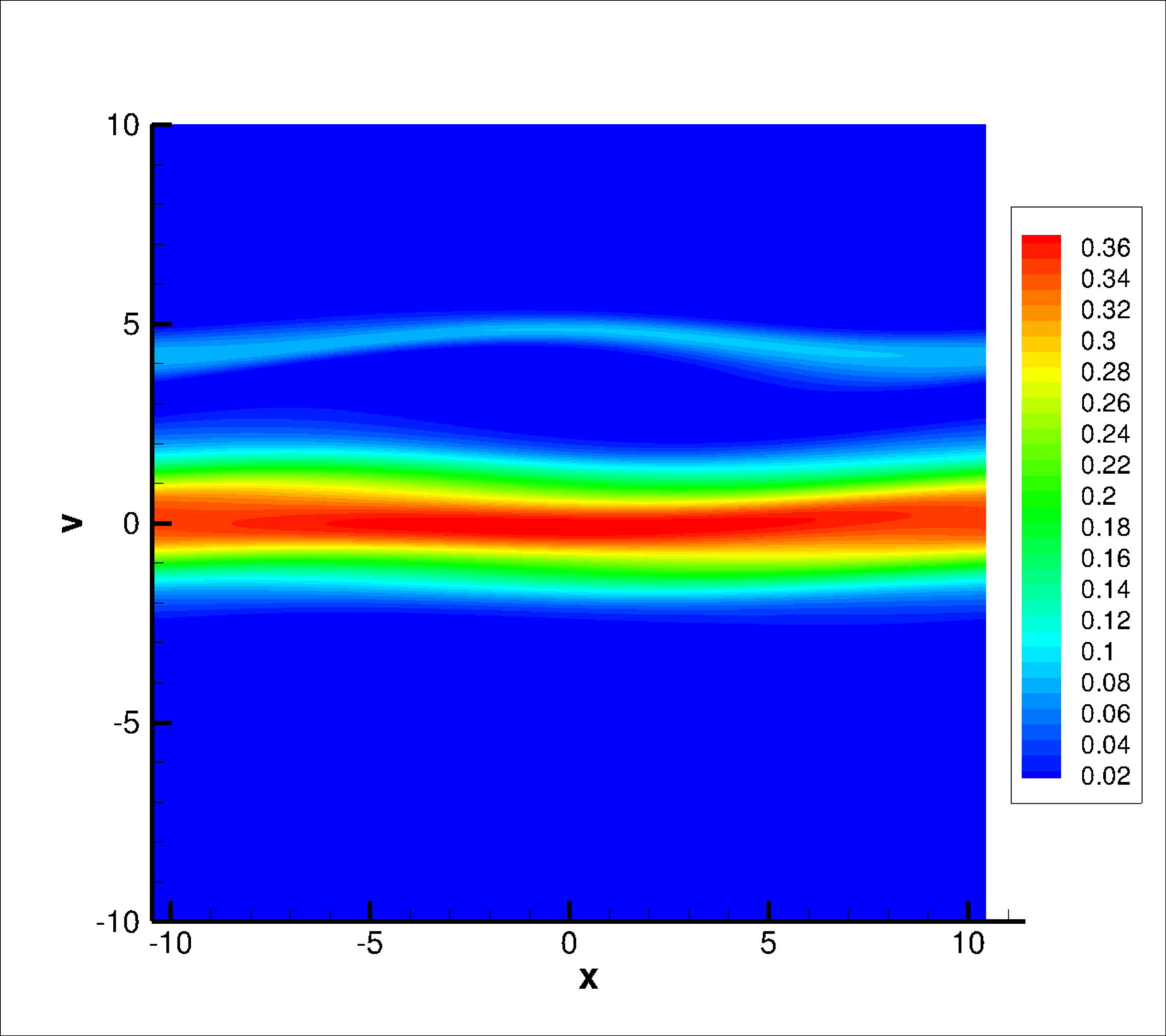}}
\subfigure[$T=20$]{
\includegraphics[width=0.4\textwidth]{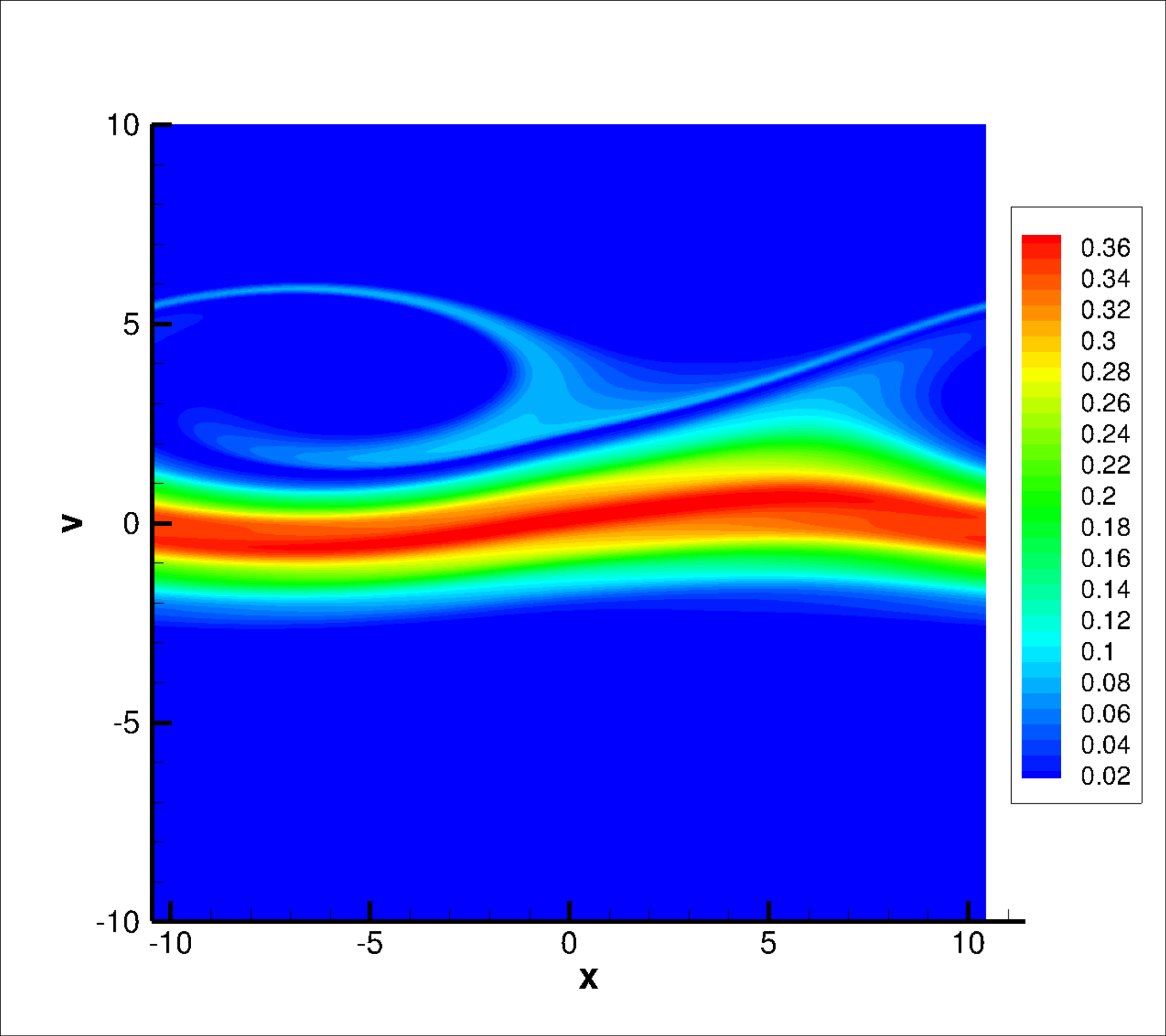}}
\subfigure[$T=30$]{
\includegraphics[width=0.4\textwidth]{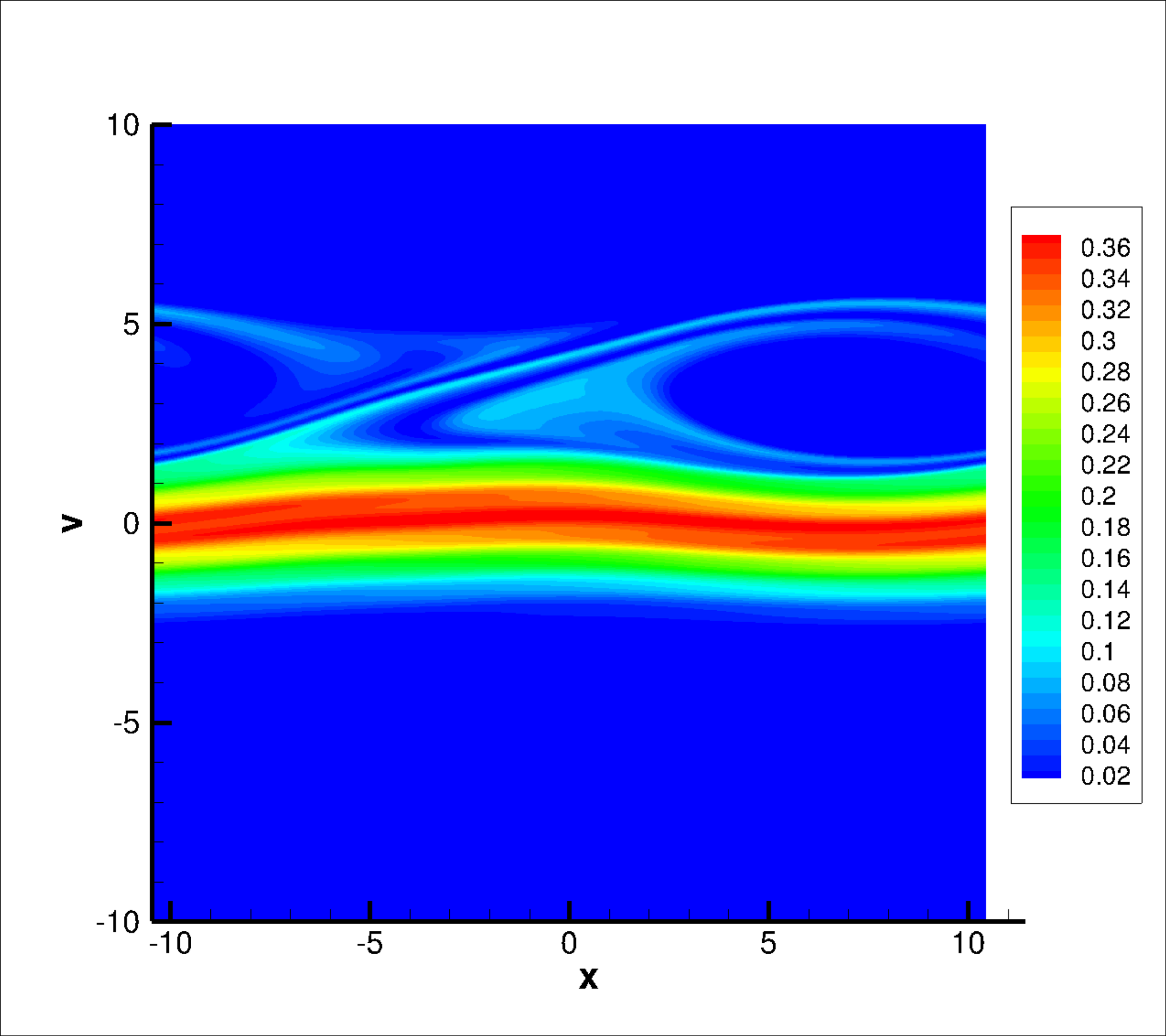}}
\subfigure[$T=40$]{
\includegraphics[width=0.4\textwidth]{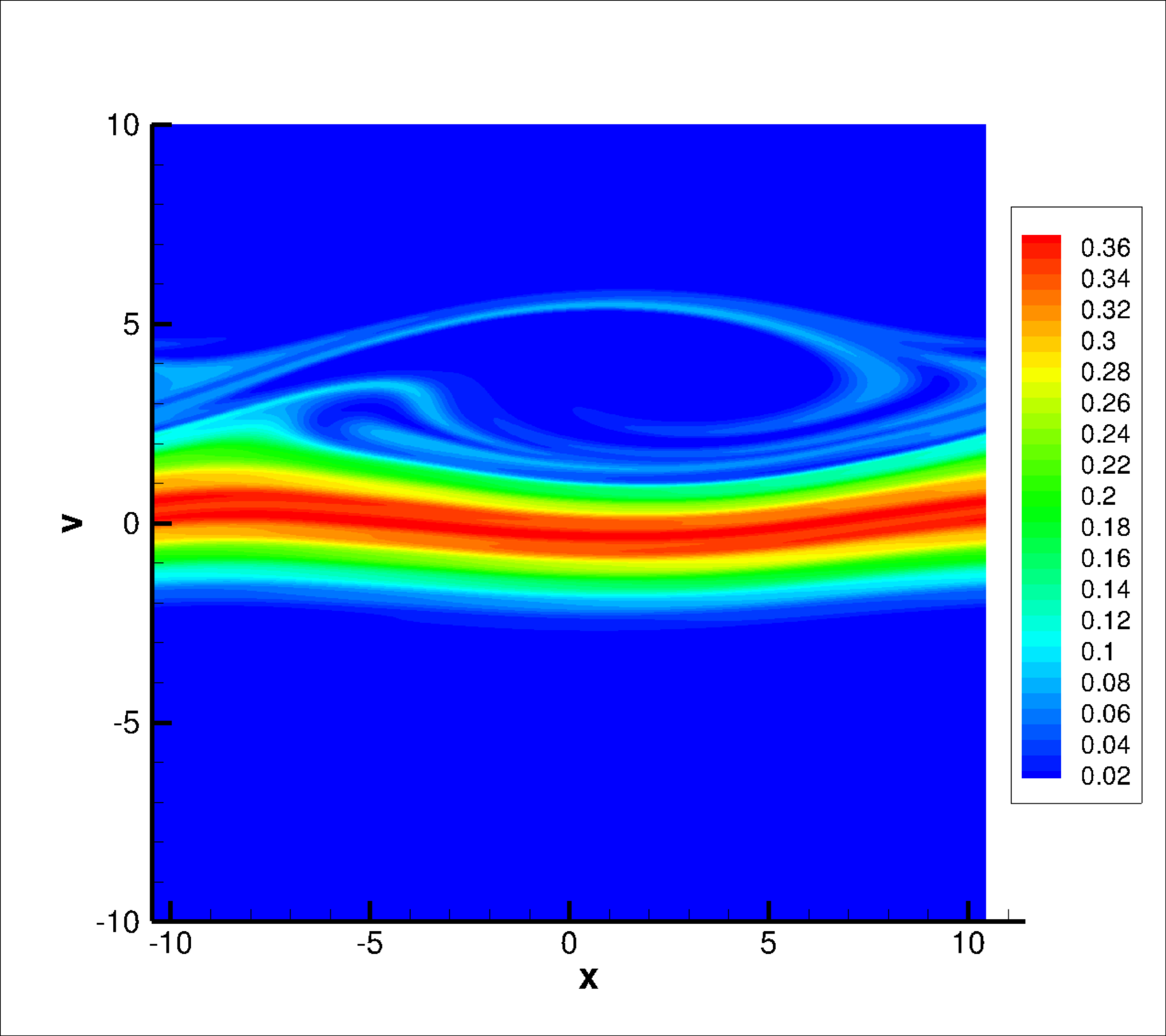}}
\subfigure[$T=50$]{
\includegraphics[width=0.4\textwidth]{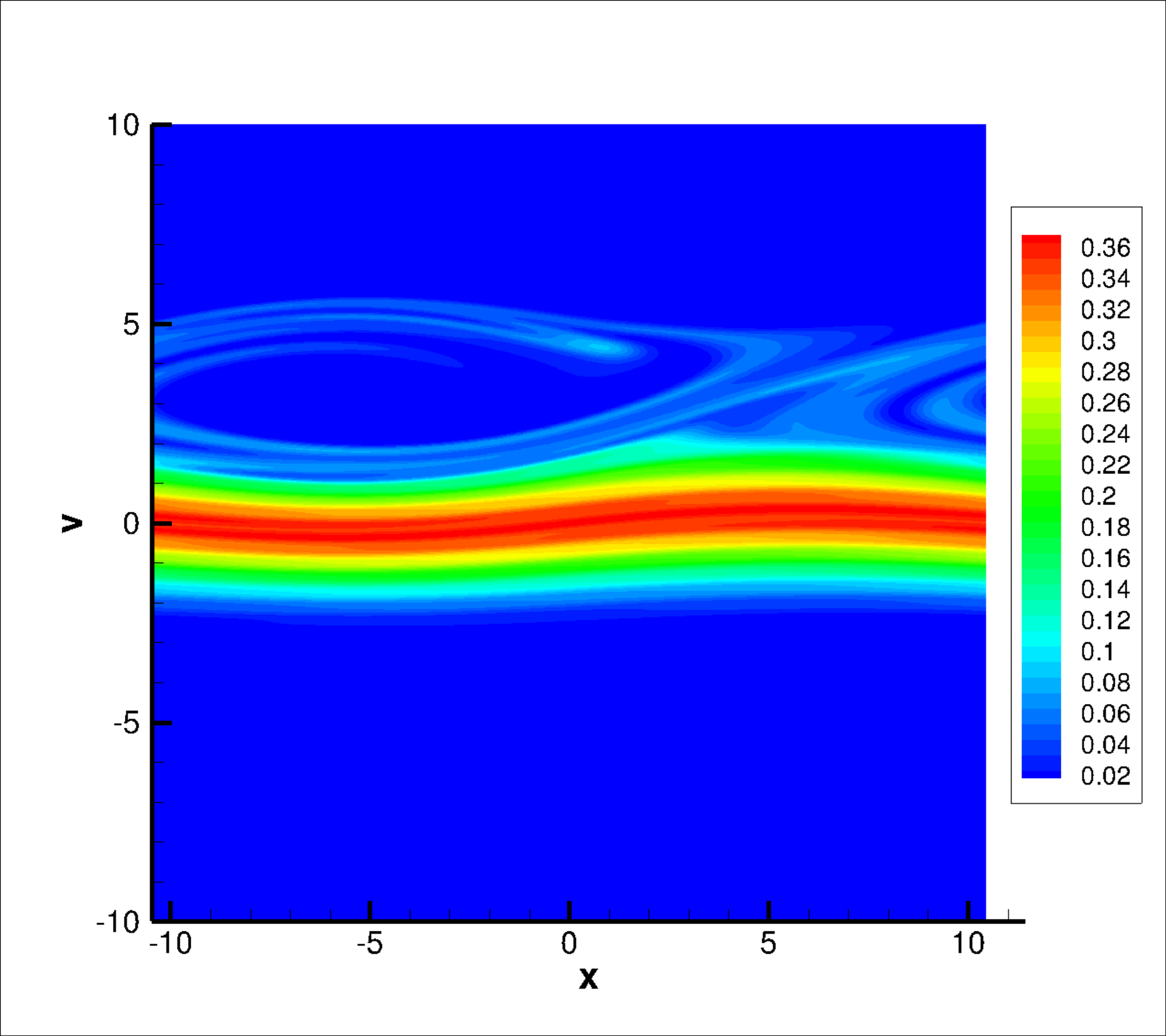}}
\subfigure[$T=60$]{
\includegraphics[width=0.4\textwidth]{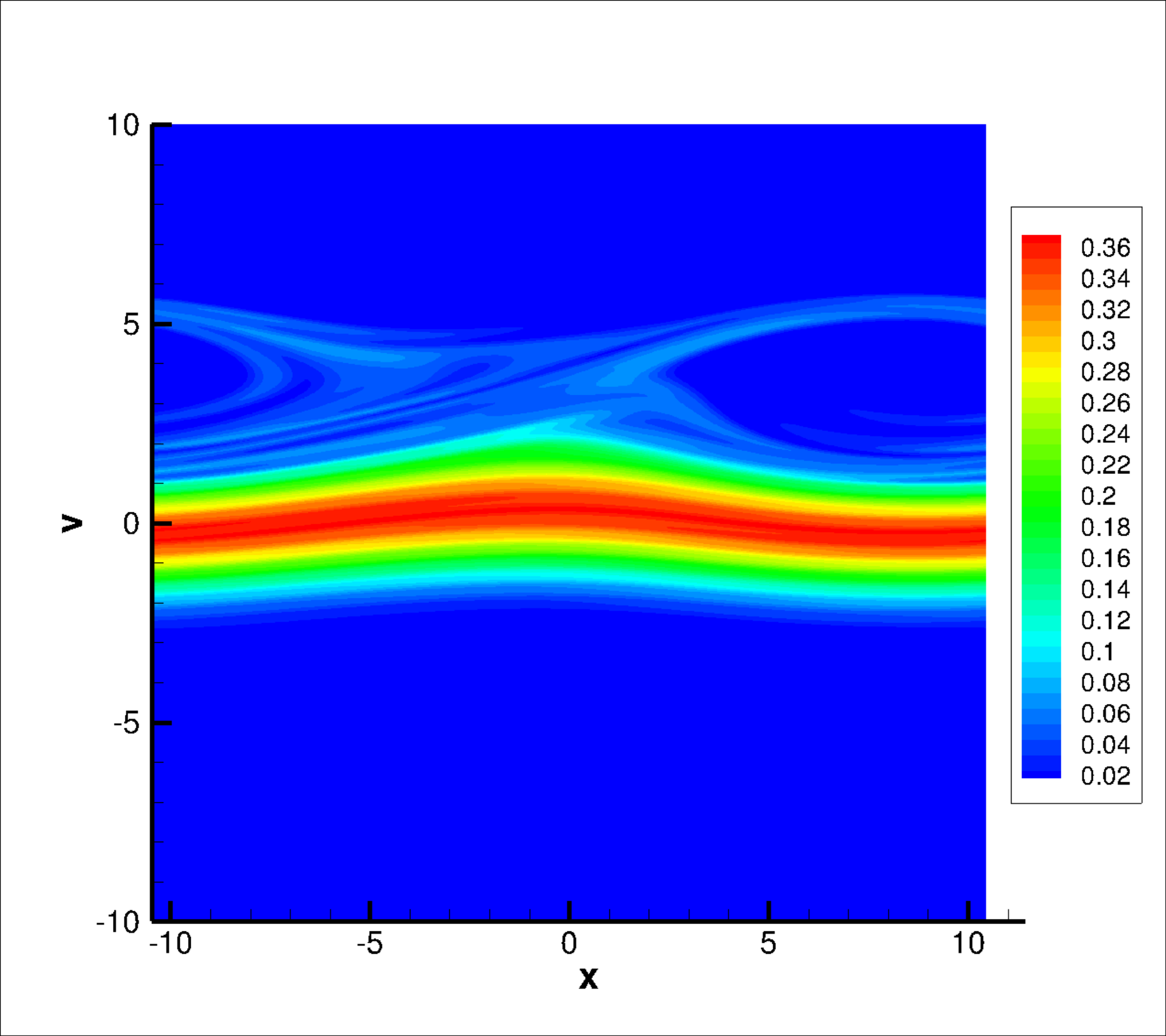}}
\caption{\em Bump-on-tail instability. WENO3, with $256\times1024$ grid points.}
\label{Fig9}
\end{figure}

\begin{figure}
\centering
\subfigure[$T=10$]{
\includegraphics[width=0.4\textwidth]{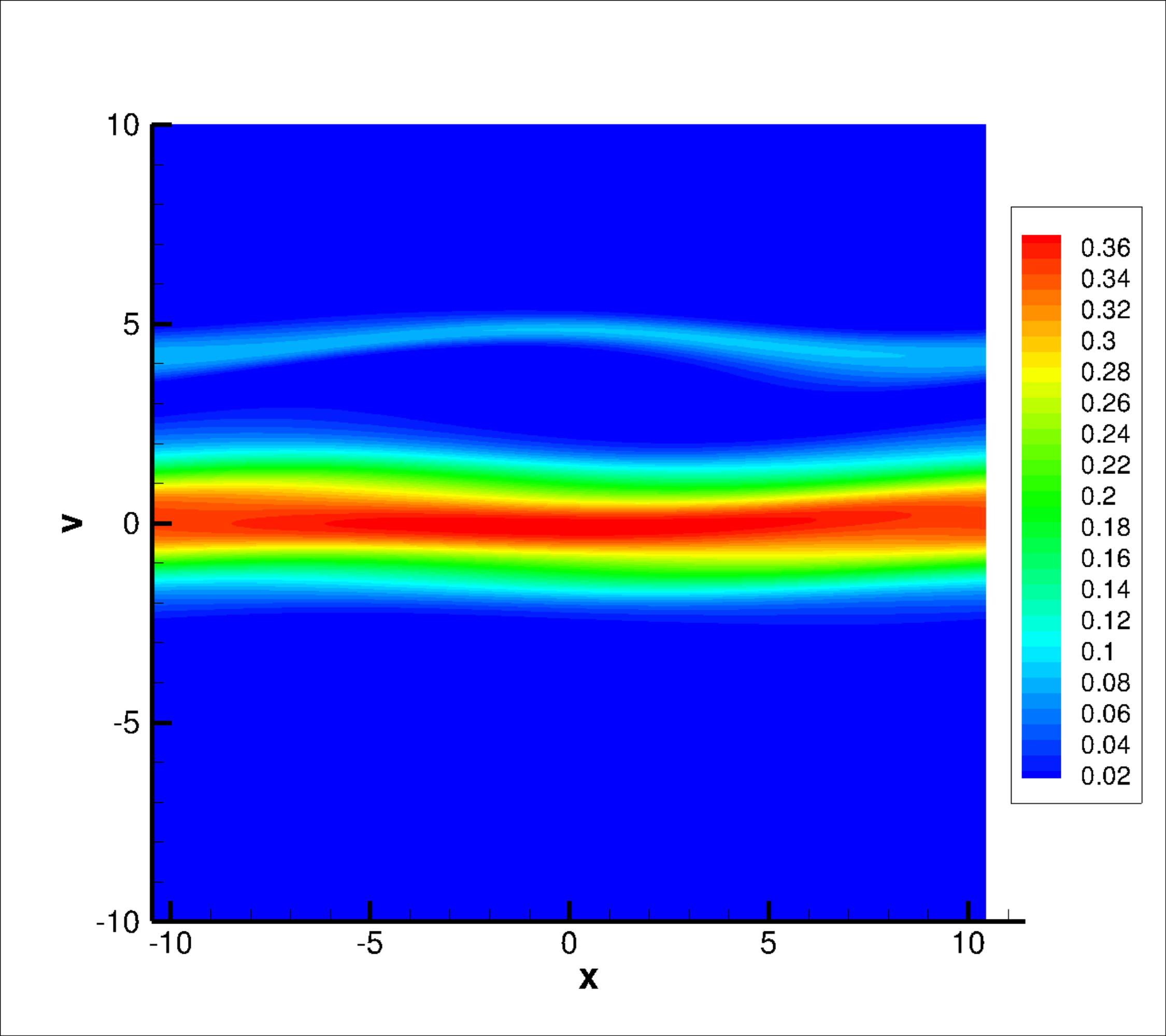}}
\subfigure[$T=20$]{
\includegraphics[width=0.4\textwidth]{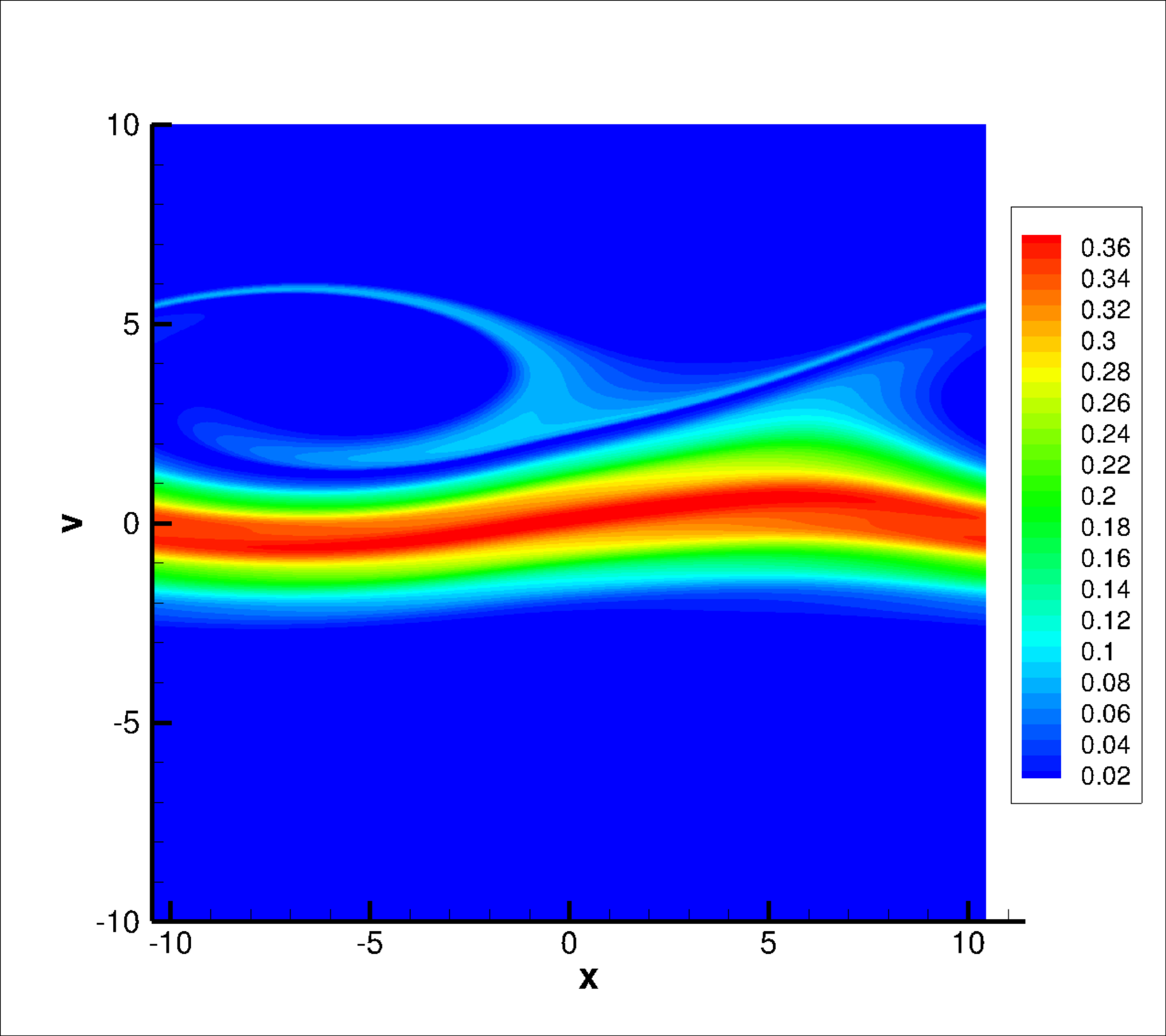}}
\subfigure[$T=30$]{
\includegraphics[width=0.4\textwidth]{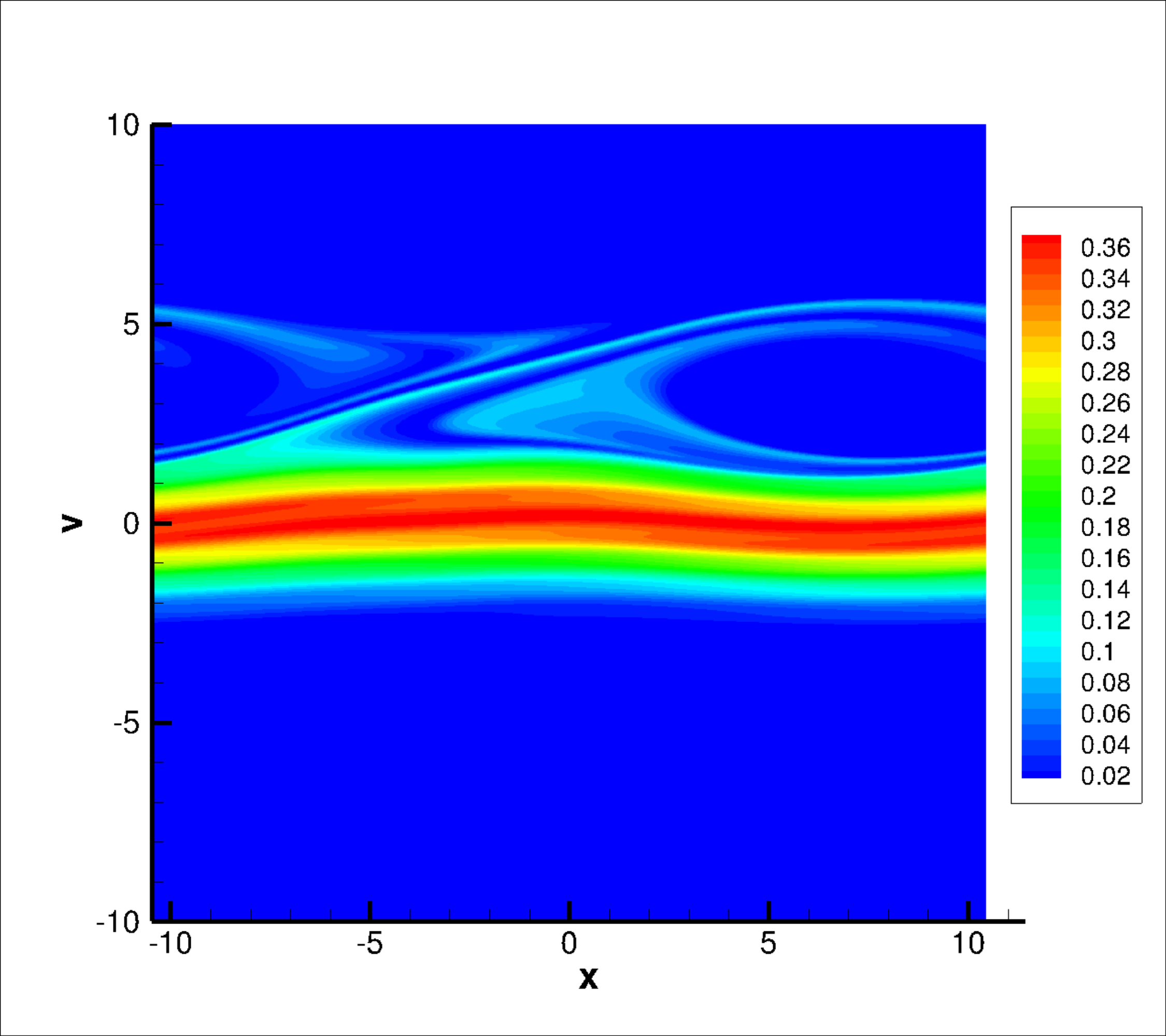}}
\subfigure[$T=40$]{
\includegraphics[width=0.4\textwidth]{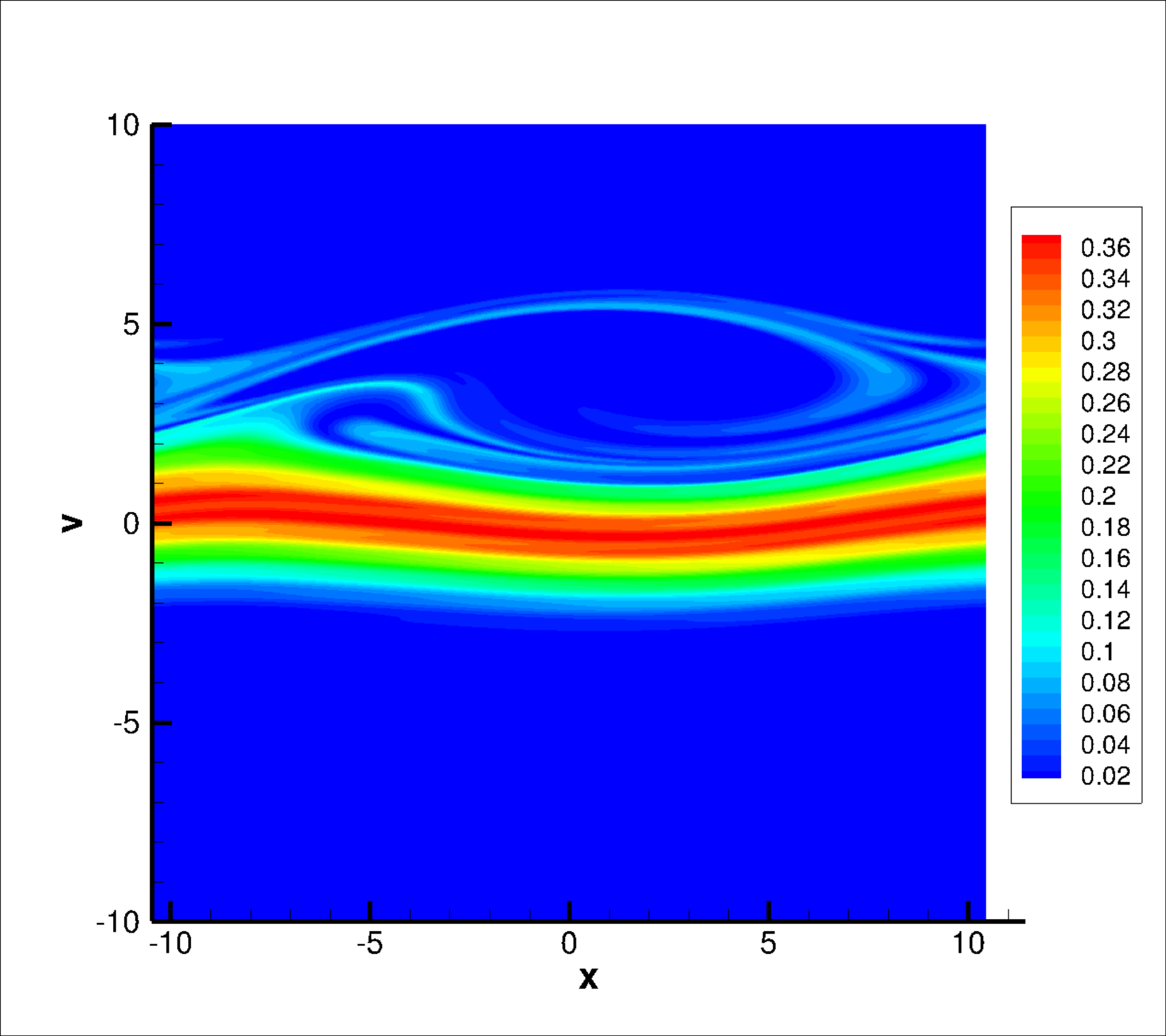}}
\subfigure[$T=50$]{
\includegraphics[width=0.4\textwidth]{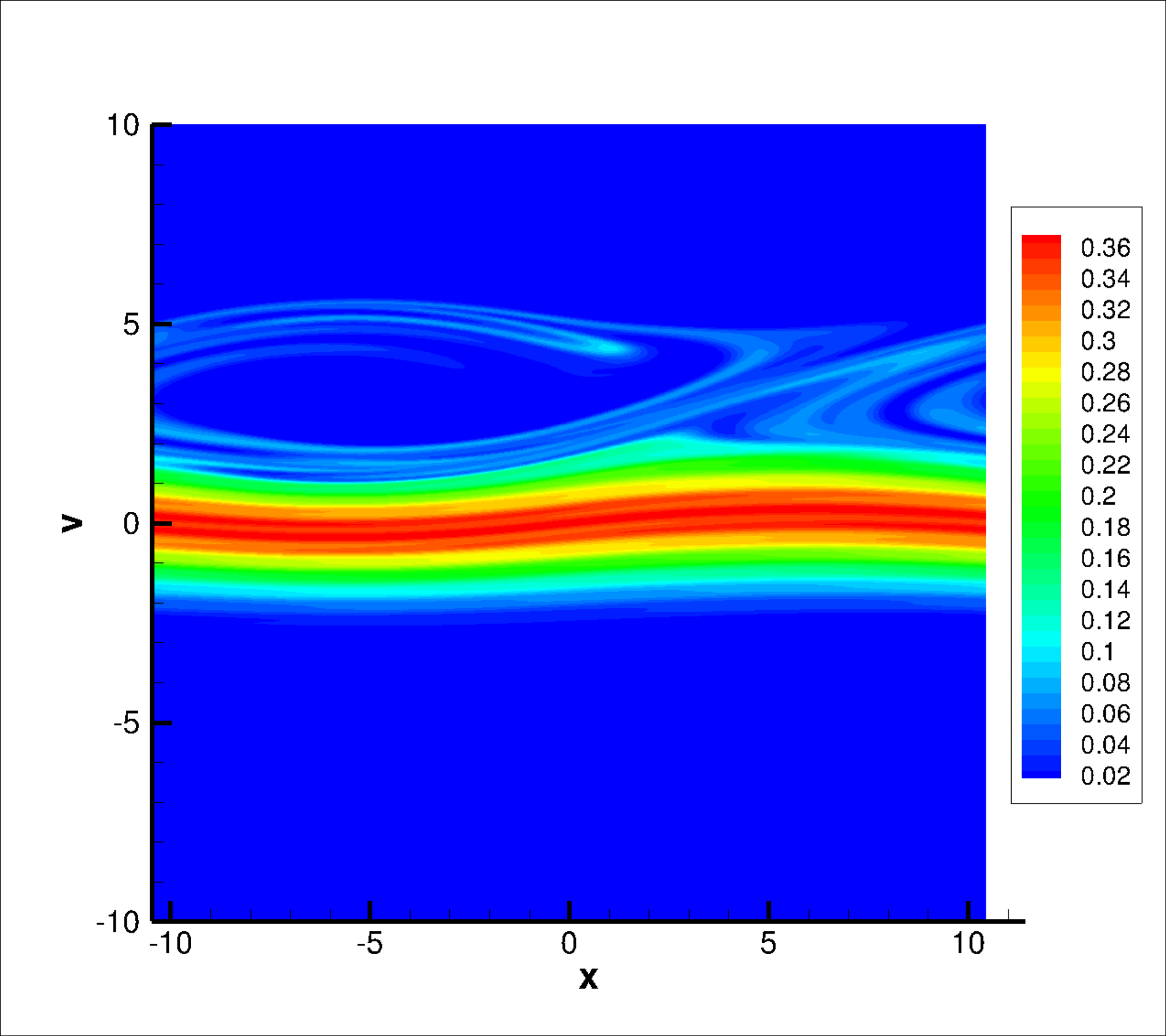}}
\subfigure[$T=60$]{
\includegraphics[width=0.4\textwidth]{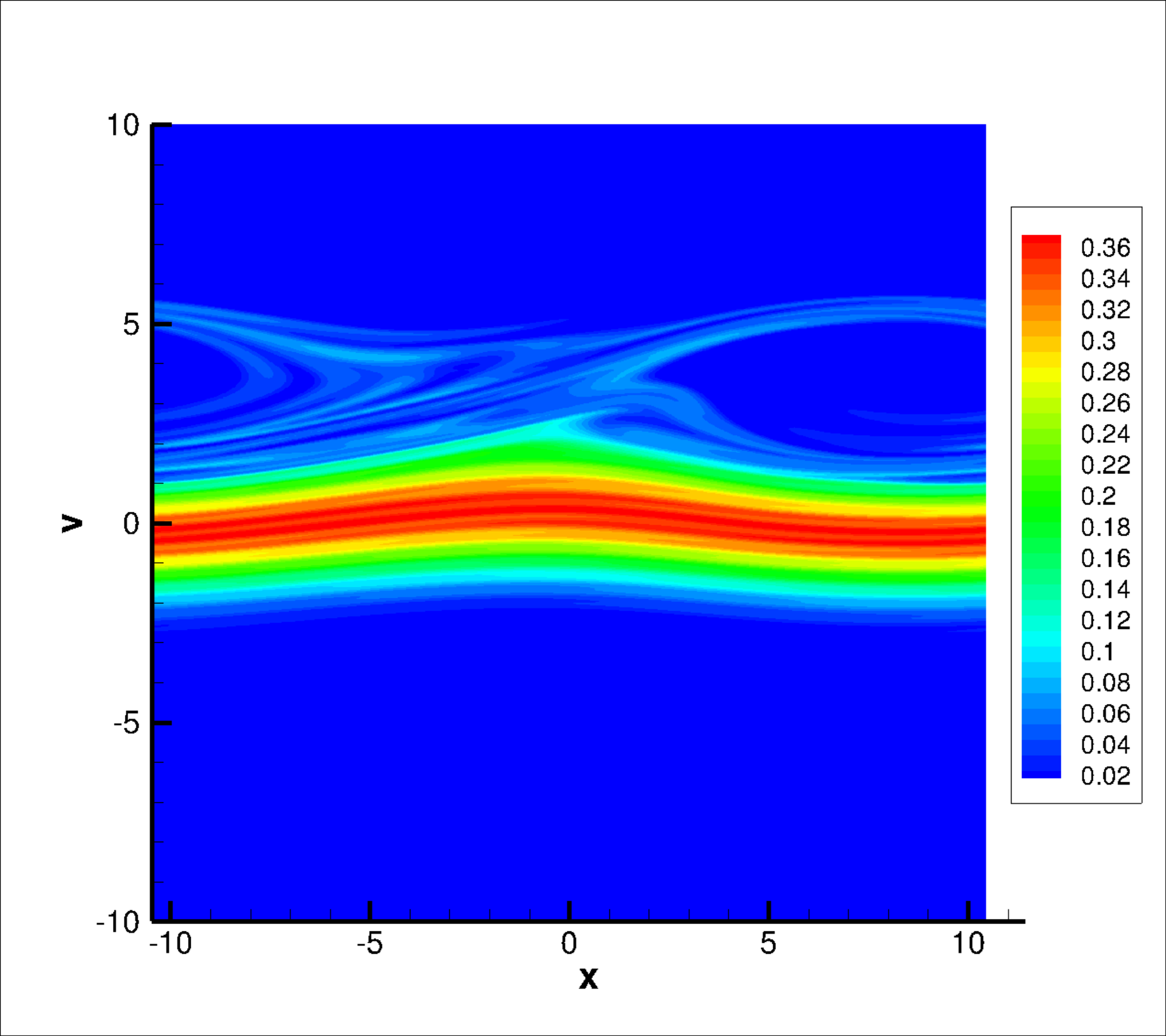}}
\caption{\em Bump-on-tail instability. WENO5, with $256\times1024$ grid points.}
\label{Fig10}
\end{figure}

For all problems, the WENO5 scheme performs better than the WENO3 scheme, showed in Figure \ref{Fig3}-\ref{Fig10}. It is observed that both schemes preserve $L_{1}$ norm (Figure \ref{Fig13}) up to the machine error, since the scheme are both conservative and positivity-preserving if the PP-limiter is applied.
While the WENO5 scheme does a better job in preserving the energy, $L_2$ norm and momentum than the WENO3 scheme. 


 We test Two-stream instability I with second order central difference scheme in time, which is non-conditionally stable and conditionally SSP.  CFL is taken as 10 such that the scheme is stable but not SSP. In space, we use WENO5. In Figure \ref{Fig17}, we can observe that the scheme simulates the structures well. We remark that the SSP property is not always desired for the VP simulations, especially when the solution is smooth.

\begin{table}[htb]
\caption{\label{tab5}\em Accuracy on strong Landau damping. $T=10$.}
\centering
\bigskip
\begin{small}
\begin{tabular}{|c|c|c|c|c|c|c|c|}
  \hline
 &  &  $N_x\times N_v$ &  $L_1$ errors &   order  &  $L_\infty$ error  & order & min value \\\hline
 \multirow{8}{1.5cm}{WENO3} & \multirow{4}{2cm}{Without\\ PP-limiters}  &
	       $32\times64$ &  1.49E+00  &   --     &  8.98E-02  &    --     & -1.00E-04  \\  		
 & &   $64\times128$ &  3.04E-01  &  2.30  &  2.37E-02  &  1.92  &  -8.67E-05  \\
 & & $128\times256$ &  3.78E-02  &  3.01  &  3.56E-03  &  2.73  &  -3.48E-05  \\
 & & $256\times512$ &  4.43E-03  &  3.09  &  4.05E-04  &  3.14  &  -5.84E-06  \\\cline{2-8}
 & \multirow{4}{2cm}{With\\ PP-limiters}  &
	       $32\times64$ &  1.49E+00  &   --     &  8.98E-02  &    --     &  9.99E-17  \\  		
 & &   $64\times128$ &  3.03E-01  &  2.30  &  2.37E-02  &  1.92  &  1.00E-16  \\
 & & $128\times256$ &  3.78E-02  &  3.01  &  3.56E-03  &  2.73  &  1.00E-16  \\
 & & $256\times512$ &  4.44E-03  &  3.09  &  4.05E-04  &  3.14  &  1.00E-16  \\\hline
 \multirow{8}{1.5cm}{WENO5} & \multirow{4}{2cm}{Without\\ PP-limiters}  &
    		   $32\times64$ &  5.70E-01  &   --     &  4.28E-02  &    --     &  -1.79E-04  \\  		
 & &   $64\times128$ &  4.32E-02  &  3.72  &  5.02E-03  &  3.09  &  -1.47E-05  \\
 & & $128\times256$ &  1.46E-03  &  4.89  &  1.87E-04  &  4.75  &  -2.17E-06  \\
 & & $256\times512$ &  3.98E-05  &  5.19  &  4.82E-06  &  5.28  &  -1.38E-07  \\\cline{2-8}
 & \multirow{4}{2cm}{With\\ PP-limiters}  &
    		   $32\times64$ &  5.70E-01  &    --    &  4.28E-02  &    --     &  1.00E-16  \\  		
 & &   $64\times128$ &  4.31E-02  &  3.72  &  5.02E-03  &  3.09  &  1.00E-16  \\
 & & $128\times256$ &  1.46E-03  &  4.88  &  1.87E-04  &  4.74  &  1.00E-16  \\
 & & $256\times512$ &  4.01E-05  &  5.19  &  4.83E-06  &  5.28  &  1.00E-16  \\\hline
\end{tabular}
\end{small}
\end{table}

\begin{table}[htb]
\caption{\label{tab6}\em Accuracy on two-stream instability I. $T=10$.}
\centering
\bigskip
\begin{small}
\begin{tabular}{|c|c|c|c|c|c|c|c|}
  \hline
 &  &  $N_x\times N_v$ &  $L_1$ errors &   order  &  $L_\infty$ error  & order & min value \\\hline
 \multirow{8}{1.5cm}{WENO3} & \multirow{4}{2cm}{Without \\ PP-limiters}  &
  		   $32\times64$ &  7.17E-02  &   --     &  3.56E-03  &    --     &  -3.47E-06  \\  	
 & &   $64\times128$ &  9.39E-03  &  2.93  &  4.03E-04  &  3.14  &  -3.87E-06  \\
 & & $128\times256$ &  1.05E-03  &  3.17  &  4.51E-05  &  3.16  &  -4.88E-07  \\
 & & $256\times512$ &  1.27E-04  &  3.05  &  5.44E-06  &  3.05  &  1.45E-08  \\\cline{2-8}
 & \multirow{4}{2cm}{With \\ PP-limiters}  &
  		   $32\times64$ &  7.18E-02  &   --     &  3.56E-03  &    --     &  1.00E-16  \\  	
 & &   $64\times128$ &  9.38E-03  &  2.94  &  4.03E-04  &  3.14  &  1.00E-16  \\
 & & $128\times256$ &  1.05E-03  &  3.16  &  4.51E-05  &  3.16  &  1.00E-16  \\
 & & $256\times512$ &  1.27E-04  &  3.05  &  5.44E-06  &  3.05  &  1.45E-08  \\
  \hline
 \multirow{8}{1.5cm}{WENO5} & \multirow{4}{2cm}{Without \\ PP-limiters}  &
  		   $32\times64$ &  3.27E-02  &   --     &  1.44E-03  &    --     &  2.58E-07  \\  		
 & &   $64\times128$ &  1.13E-03  &  4.86  &  5.75E-05  &  4.65  &  -1.21E-07  \\
 & & $128\times256$ &  2.51E-05  &  5.49  &  1.84E-06  &  4.97  &  8.42E-08  \\
 & & $256\times512$ &  5.61E-07  &  5.49  &  4.49E-08  &  5.35  &  6.75E-08  \\\cline{2-8}  
 & \multirow{4}{2cm}{With \\ PP-limiters}  &
  		   $32\times64$ &  3.27E-02  &   --     &  1.44E-03  &    --     &  2.58E-07  \\  		
 & &   $64\times128$ &  1.13E-03  &  4.86  &  5.75E-05  &  4.65  &  1.00E-16  \\
 & & $128\times256$ &  2.51E-05  &  5.49  &  1.84E-06  &  4.97  &  8.42E-08  \\
 & & $256\times512$ &  5.61E-07  &  5.49  &  4.49E-08  &  5.35  &  6.75E-08  \\  
\hline
\end{tabular}
\end{small}
\end{table}

\begin{table}[htb]
\caption{\label{tab7}\em Accuracy on two-stream instability II. $T=10$.}
\centering
\bigskip
\begin{small}
\begin{tabular}{|c|c|c|c|c|c|c|c|}
  \hline
 &  &  $N_x\times N_v$ &  $L_1$ errors &   order  &  $L_\infty$ error  & order & min value \\\hline
 \multirow{8}{1.5cm}{WENO3} & \multirow{4}{2cm}{Without \\ PP-limiters}  &
     	   $32\times64$ &  8.90E-02  &   --     &  4.56E-03  &    --     &  -9.96E-06  \\  	
 & &   $64\times128$ &  9.40E-03  &  3.24  &  5.06E-04  &  3.17  &  -6.65E-06  \\
 & & $128\times256$ &  1.05E-03  &  3.16  &  5.94E-05  &  3.09  &  -8.71E-07  \\
 & & $256\times512$ &  1.22E-04  &  3.10  &  7.06E-06  &  3.07  &  -2.15E-08  \\\cline{2-8}
 & \multirow{4}{2cm}{With \\ PP-limiters}  &
     	   $32\times64$ &  8.92E-02  &   --     &  4.56E-03  &    --     &  1.00E-16  \\  	
 & &   $64\times128$ &  9.38E-03  &  3.25  &  5.06E-04  &  3.17  &  1.00E-16  \\
 & & $128\times256$ &  1.05E-03  &  3.16  &  5.94E-05  &  3.09  &  1.00E-16  \\
 & & $256\times512$ &  1.22E-04  &  3.10  &  7.06E-06  &  3.07  &  1.00E-16  \\\hline
 \multirow{8}{1.5cm}{WENO5} & \multirow{4}{2cm}{Without \\ PP-limiters}  &
     	   $32\times64$ &  7.58E-03  &   --     &  6.37E-04  &    --     &  -3.94E-07  \\  	 	
 & &   $64\times128$ &  2.11E-04  &  5.17  &  2.36E-05  &  4.76  &  -9.53E-07  \\
 & & $128\times256$ &  5.29E-06  &  5.32  &  7.44E-07  &  4.99  &  1.70E-14  \\
 & & $256\times512$ &  1.43E-07  &  5.21  &  2.27E-08  &  5.04  &  2.43E-17  \\\cline{2-8}
 & \multirow{4}{2cm}{With \\ PP-limiters}  &
     	   $32\times64$ &  7.66E-03  &   --     &  6.37E-04  &    --     &  1.00E-16  \\  	 	
 & &   $64\times128$ &  2.17E-04  &  5.14  &  2.36E-05  &  4.76  &  1.00E-16  \\
 & & $128\times256$ &  5.31E-06  &  5.36  &  7.44E-07  &  4.99  &  5.24E-15  \\
 & & $256\times512$ &  1.43E-07  &  5.21  &  2.27E-08  &  5.04  &  1.24E-16  \\\hline
\end{tabular}
\end{small}
\end{table}

\begin{table}[htb]
\caption{\label{tab8}\em Accuracy on bump-on-tail instability. $T=10$. }
\centering
\bigskip
\begin{small}
\begin{tabular}{|c|c|c|c|c|c|c|c|}
  \hline
 &  &  $N_x\times N_v$ &  $L_1$ errors &   order  &  $L_\infty$ error  & order & min value \\\hline
 \multirow{8}{1.5cm}{WENO3} & \multirow{4}{2cm}{Without \\ PP-limiters}  &
 	       $32\times64$ &  5.78E-01  &    --     &  1.97E-02  &   --     &  -1.28E-05  \\  	
 & &   $64\times128$ &  1.10E-01  &  2.39  &  6.36E-03  &  1.63  &  -1.68E-05  \\
 & & $128\times256$ &  1.26E-02  &  3.13  &  8.16E-04  &  2.96  &  -5.76E-07  \\
 & & $256\times512$ &  1.44E-03  &  3.13  &  9.02E-05  &  3.18  &  -5.35E-07  \\\cline{2-8}
 & \multirow{4}{2cm}{With \\ PP-limiters}  &
 	       $32\times64$ &  5.78E-01  &    --     &  1.97E-01  &   --     &  1.00E-16  \\  	
 & &   $64\times128$ &  1.10E-01  &  2.39  &  6.36E-03  &  1.63  &  1.00E-16  \\
 & & $128\times256$ &  1.26E-02  &  3.13  &  8.16E-04  &  2.96  &  1.00E-16  \\
 & & $256\times512$ &  1.44E-03  &  3.13  &  9.02E-05  &  3.18  &  1.00E-16  \\\hline
 \multirow{8}{1.5cm}{WENO5} & \multirow{4}{2cm}{Without \\ PP-limiters}  &
    		   $32\times64$ &  2.25E-01  &    --    &  1.02E-02  &     --    &  -1.23E-06  \\  	 	
 & &   $64\times128$ &  1.66E-02  &  3.76  &  1.44E-03  &  2.83  &  -1.19E-06  \\
 & & $128\times256$ &  7.49E-04  &  4.47  &  7.75E-05  &  4.21  &  -1.31E-06  \\
 & & $256\times512$ &  2.08E-05  &  5.17  &  2.79E-06  &  4.80  &  -1.48E-07  \\\cline{2-8}
 & \multirow{4}{2cm}{With \\ PP-limiters}  &
    		   $32\times64$ &  2.25E-01  &   --     &  1.02E-02  &    --     &  1.00E-16  \\  	 	
 & &   $64\times128$ &  1.66E-02  &  3.76  &  1.44E-03  &  2.83  &  1.00E-16  \\
 & & $128\times256$ &  7.46E-04  &  4.47  &  7.75E-05  &  4.21  &  1.00E-16  \\
 & & $256\times512$ &  2.09E-05  &  5.16  &  2.79E-06  &  4.80  &  1.00E-16  \\\hline
\end{tabular}
\end{small}
\end{table}

\clearpage


\begin{figure}
\centering
\subfigure[Strong-Landau damping]{
\includegraphics[width=0.4\textwidth]{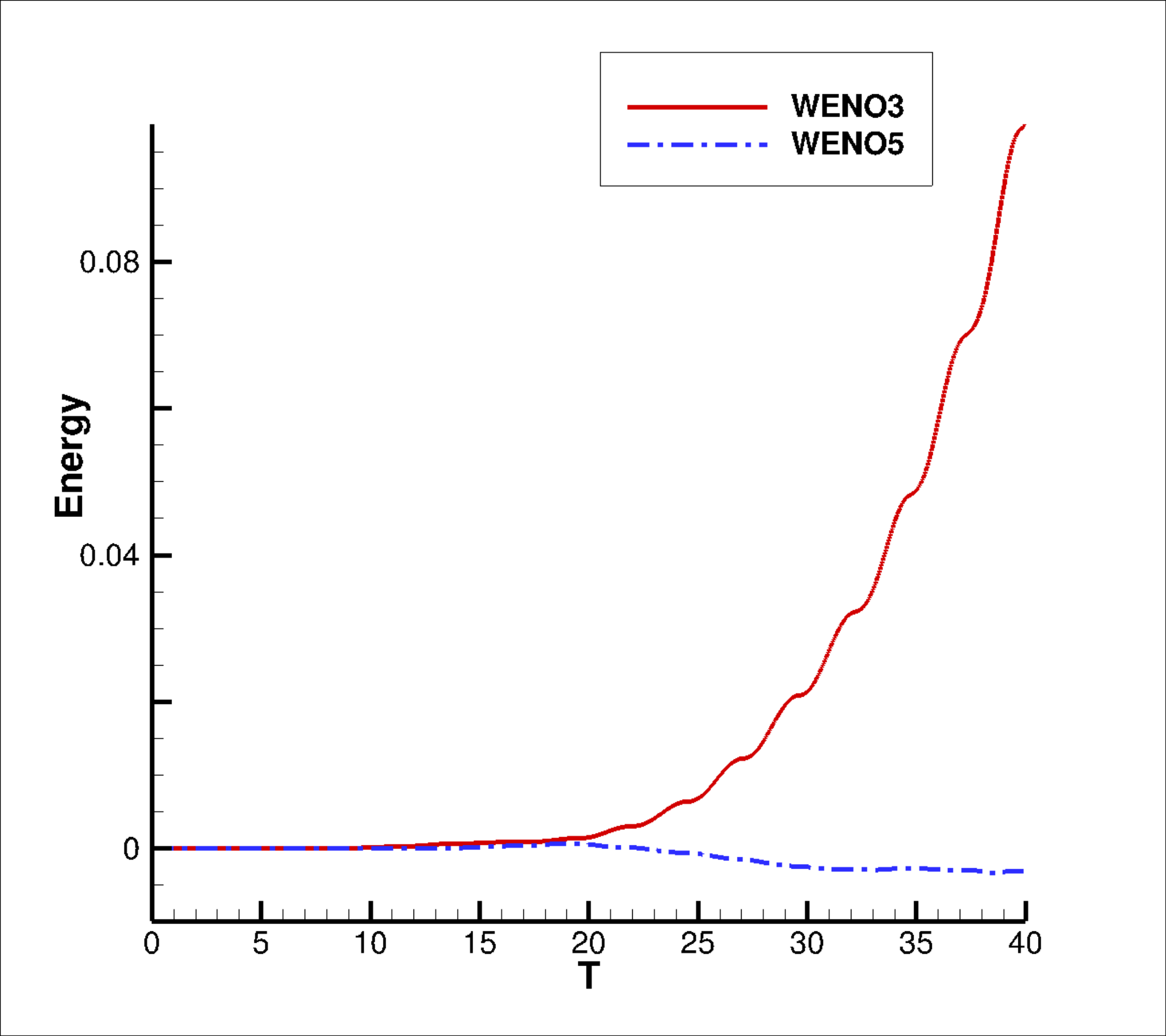}}
\subfigure[Two-stream instability I]{
\includegraphics[width=0.4\textwidth]{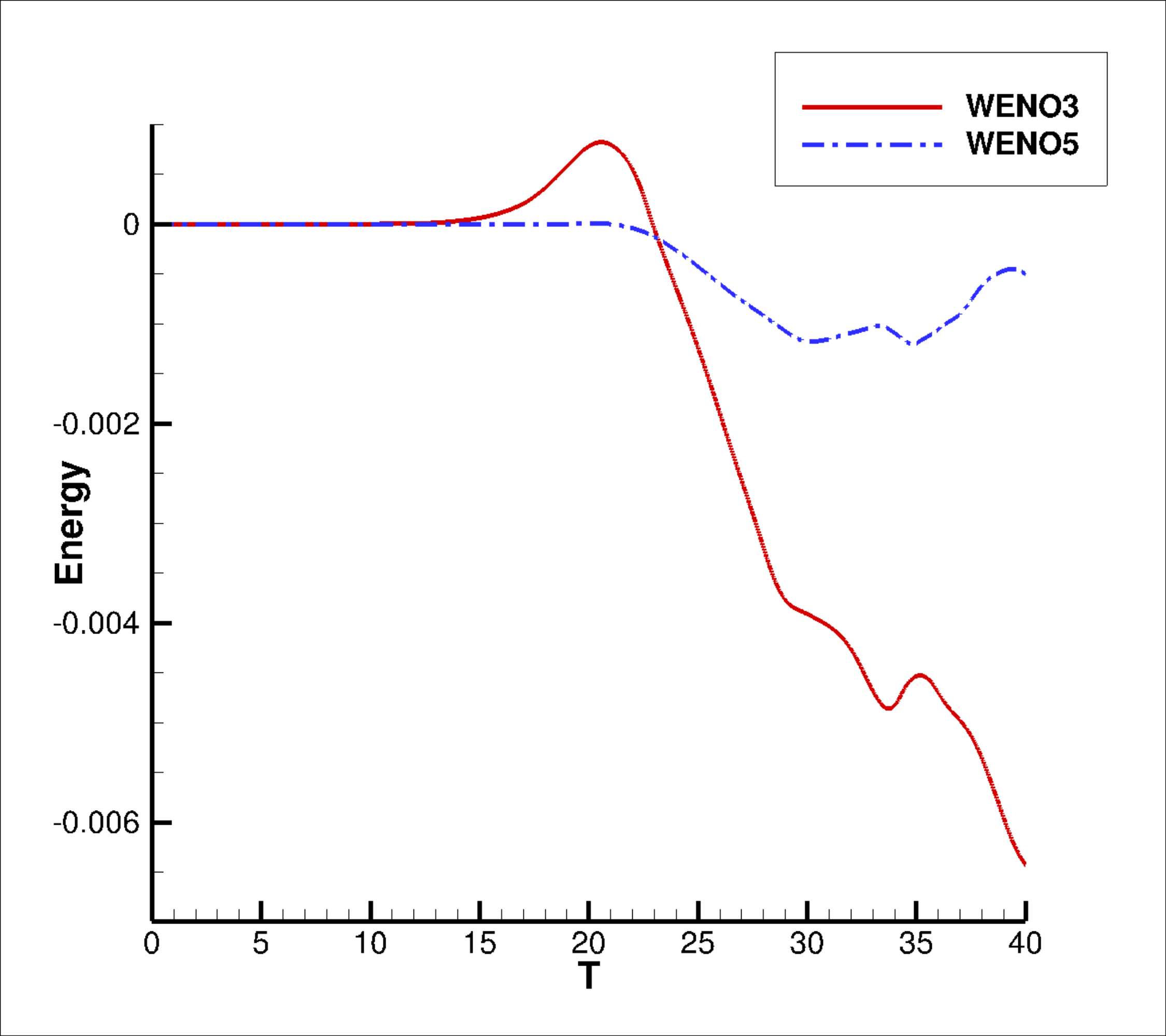}}
\subfigure[Two-stream instability II]{
\includegraphics[width=0.4\textwidth]{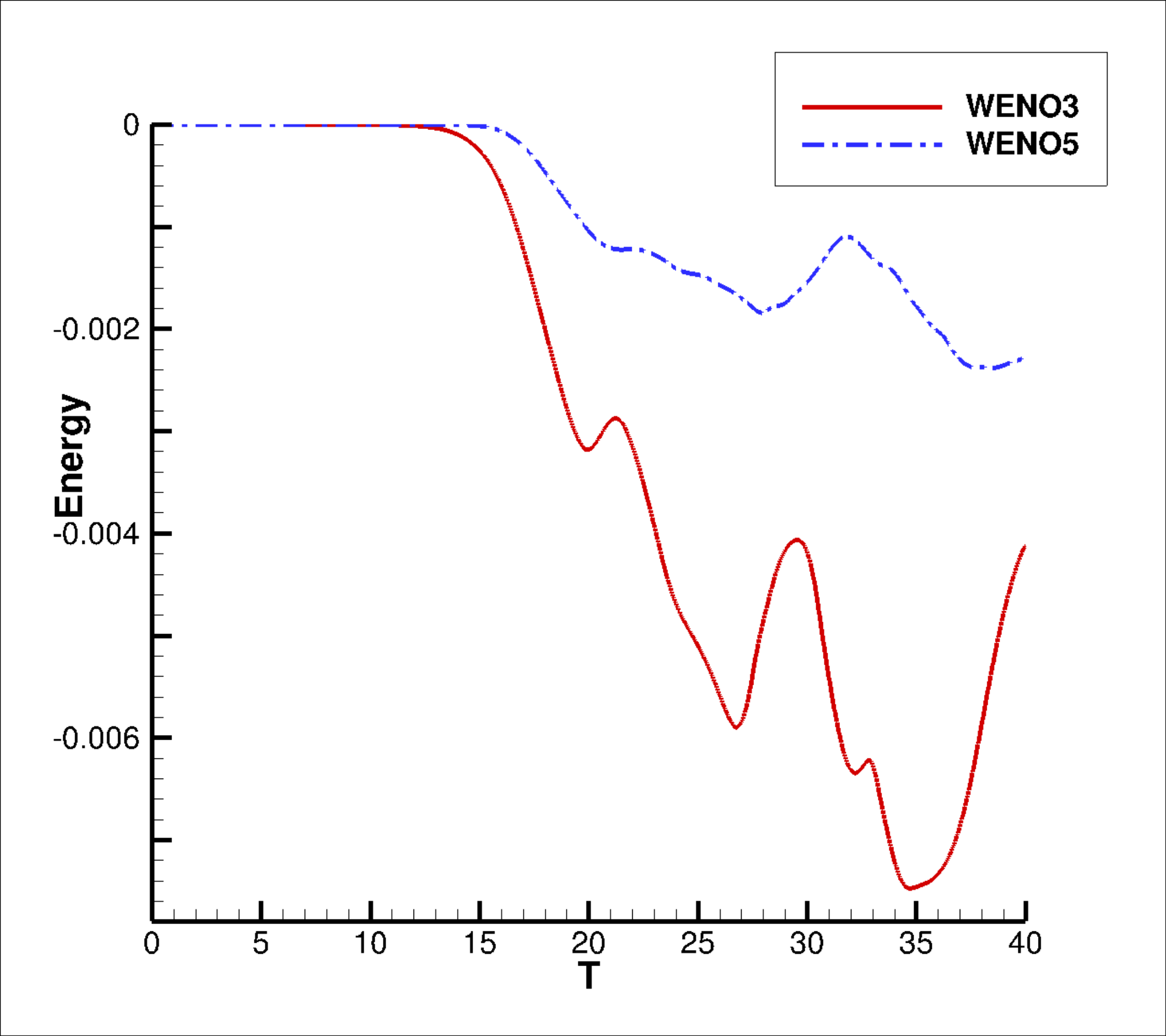}}
\subfigure[Bump-on-tail instability]{
\includegraphics[width=0.4\textwidth]{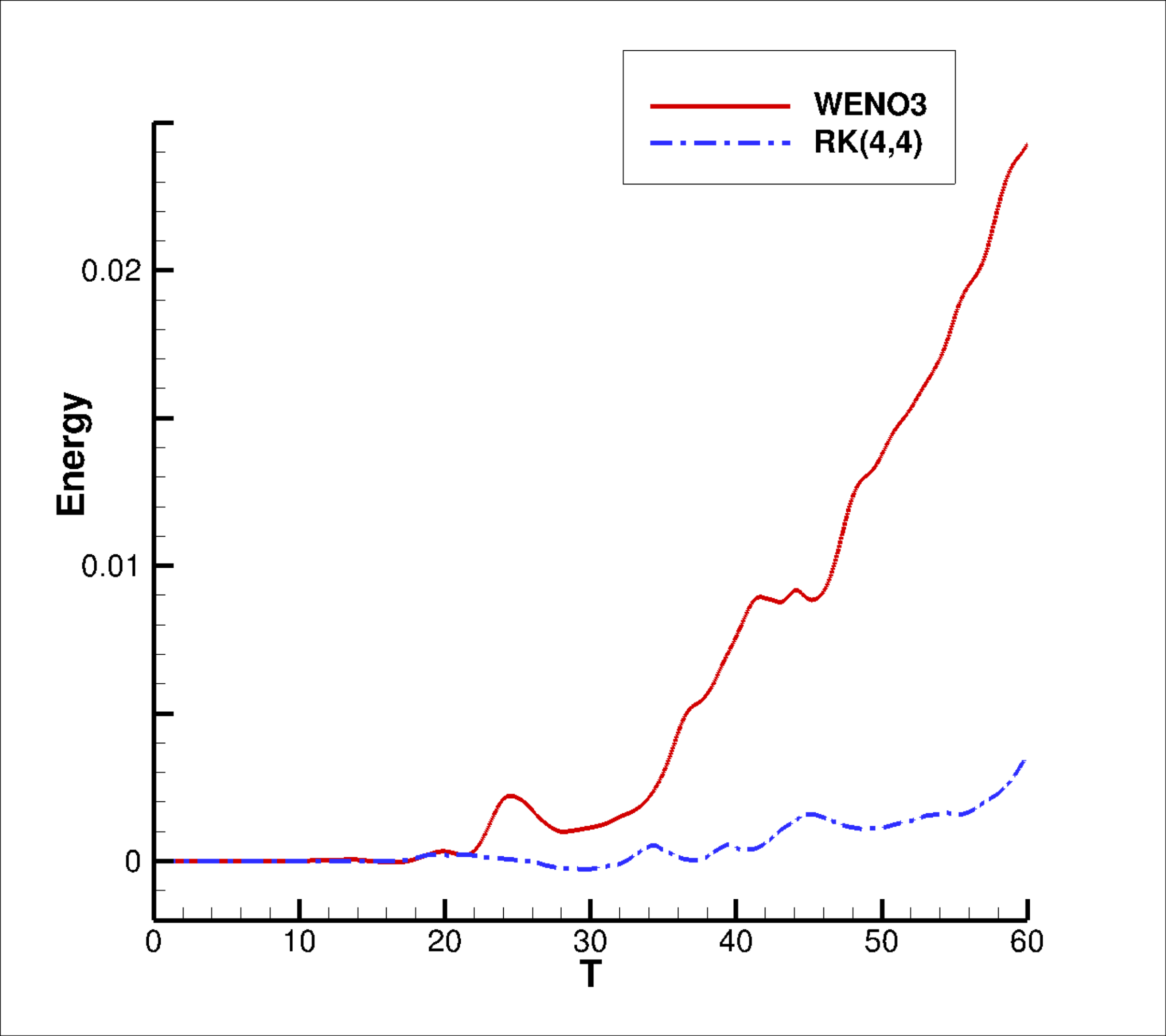}}
\caption{\em The time evolution of the relative deviation in total energy. $256\times1024$ grid points. }
\label{Fig12}
\end{figure}

\begin{figure}
\centering
\subfigure[Strong Landau damping]{
\includegraphics[width=0.4\textwidth]{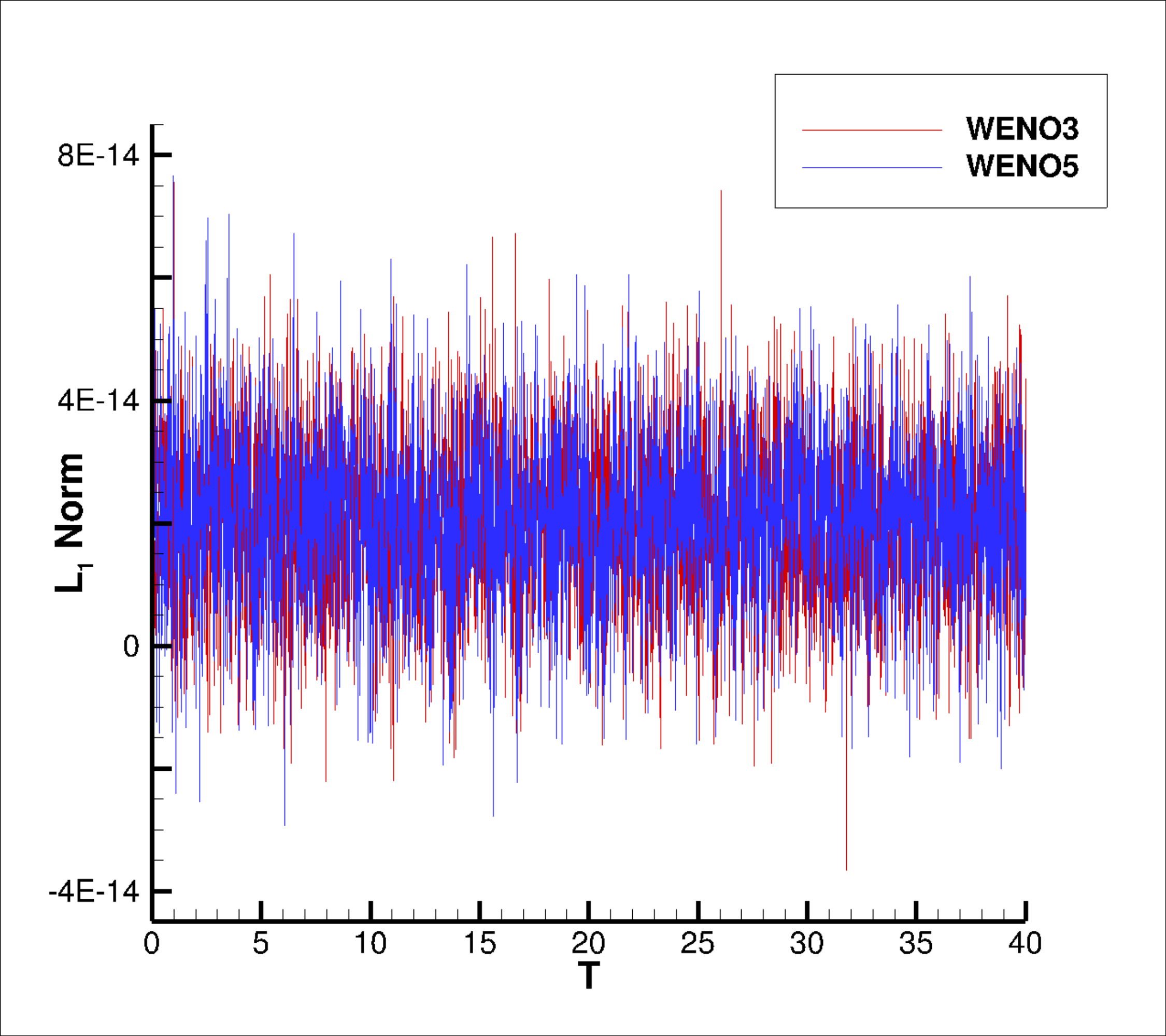}}
\subfigure[Two-stream instability I]{
\includegraphics[width=0.4\textwidth]{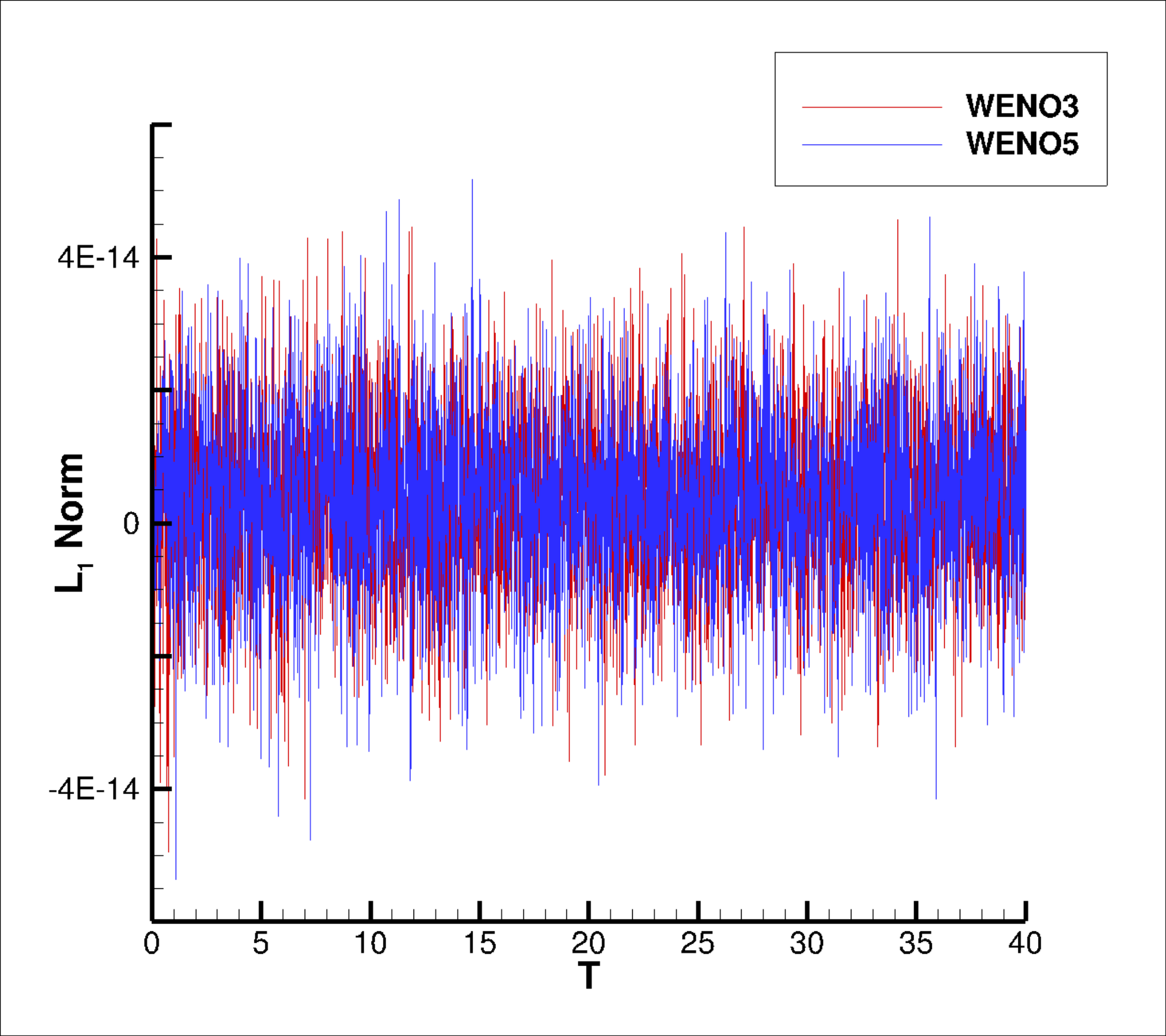}}
\subfigure[Two-stream instability II]{
\includegraphics[width=0.4\textwidth]{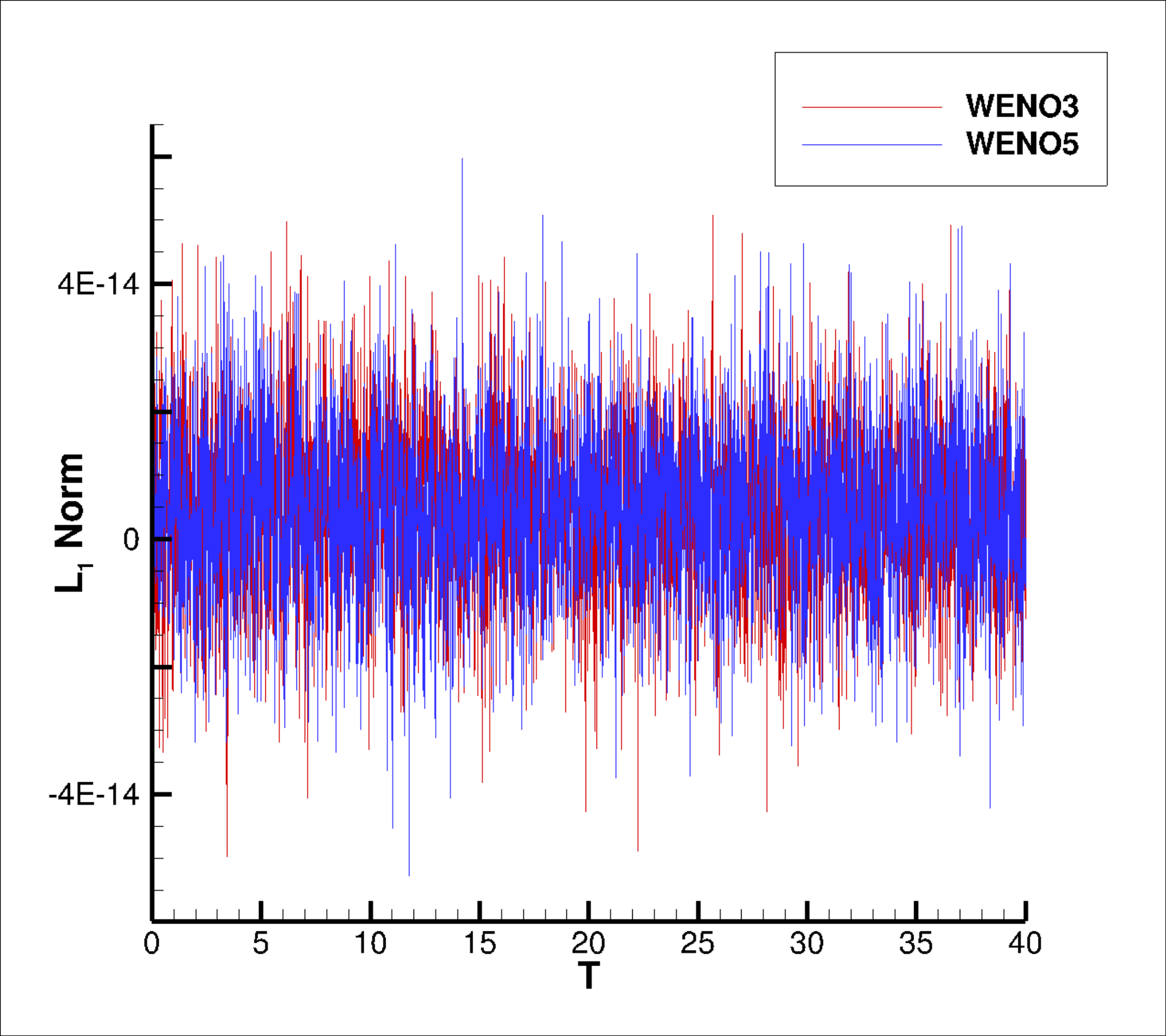}}
\subfigure[Bump-on-tail instability]{
\includegraphics[width=0.4\textwidth]{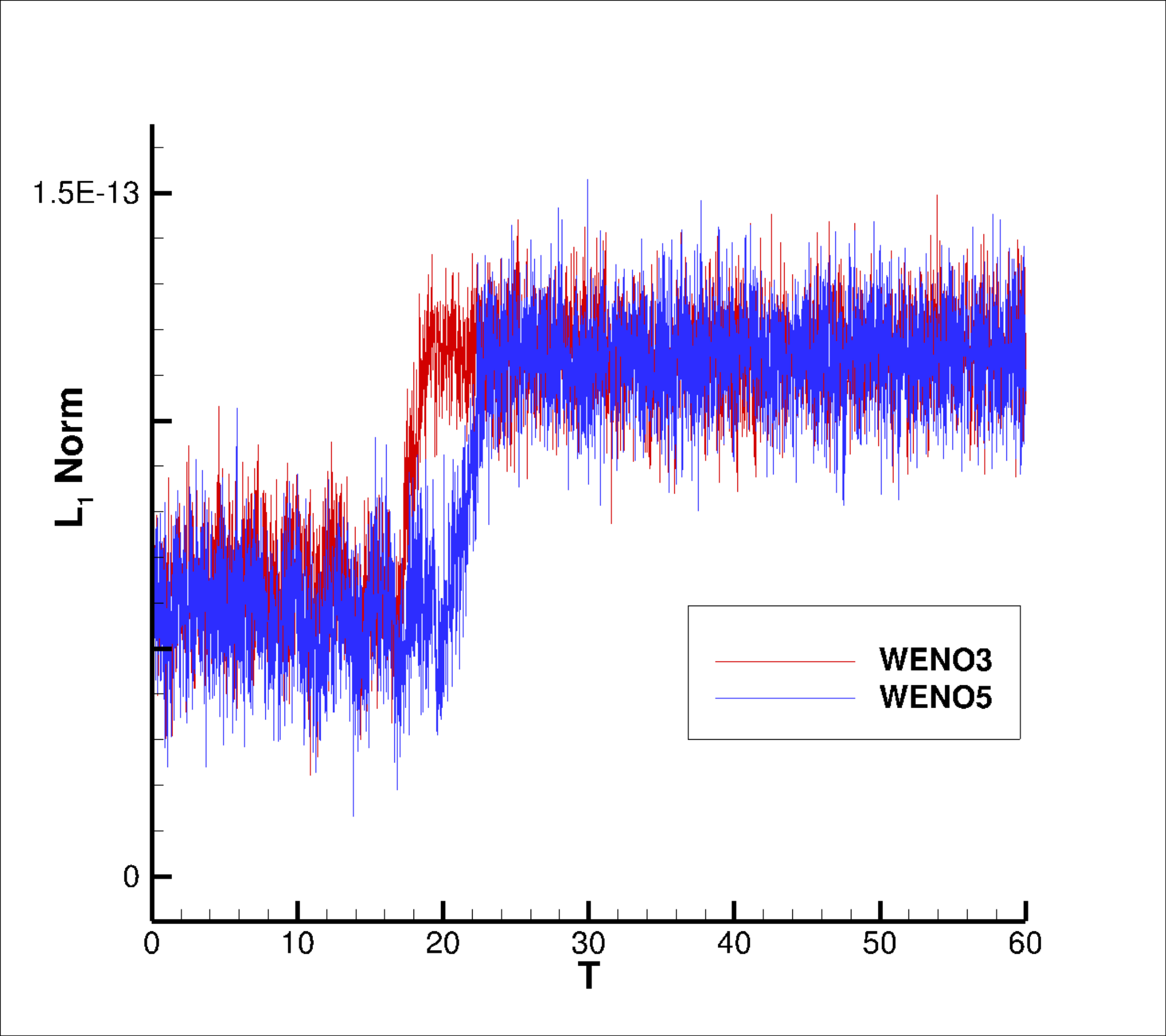}}
\caption{\em The time evolution of the relative deviation in $L_{1}$ norm. $256\times1024$ grid points. }
\label{Fig13}
\end{figure}

\begin{figure}
\centering
\subfigure[Strong Landau damping]{
\includegraphics[width=0.4\textwidth]{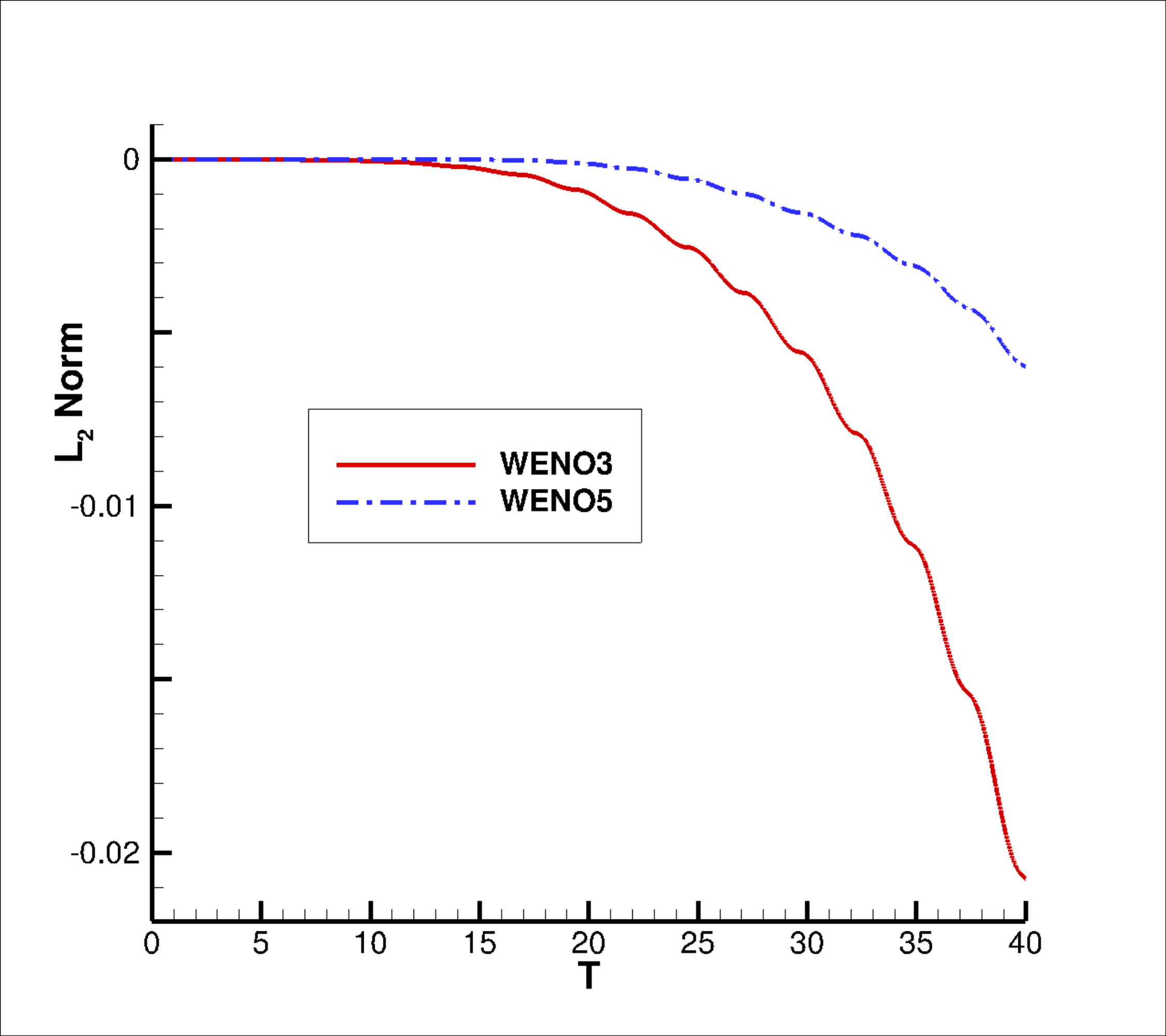}}
\subfigure[Two-stream instability I]{
\includegraphics[width=0.4\textwidth]{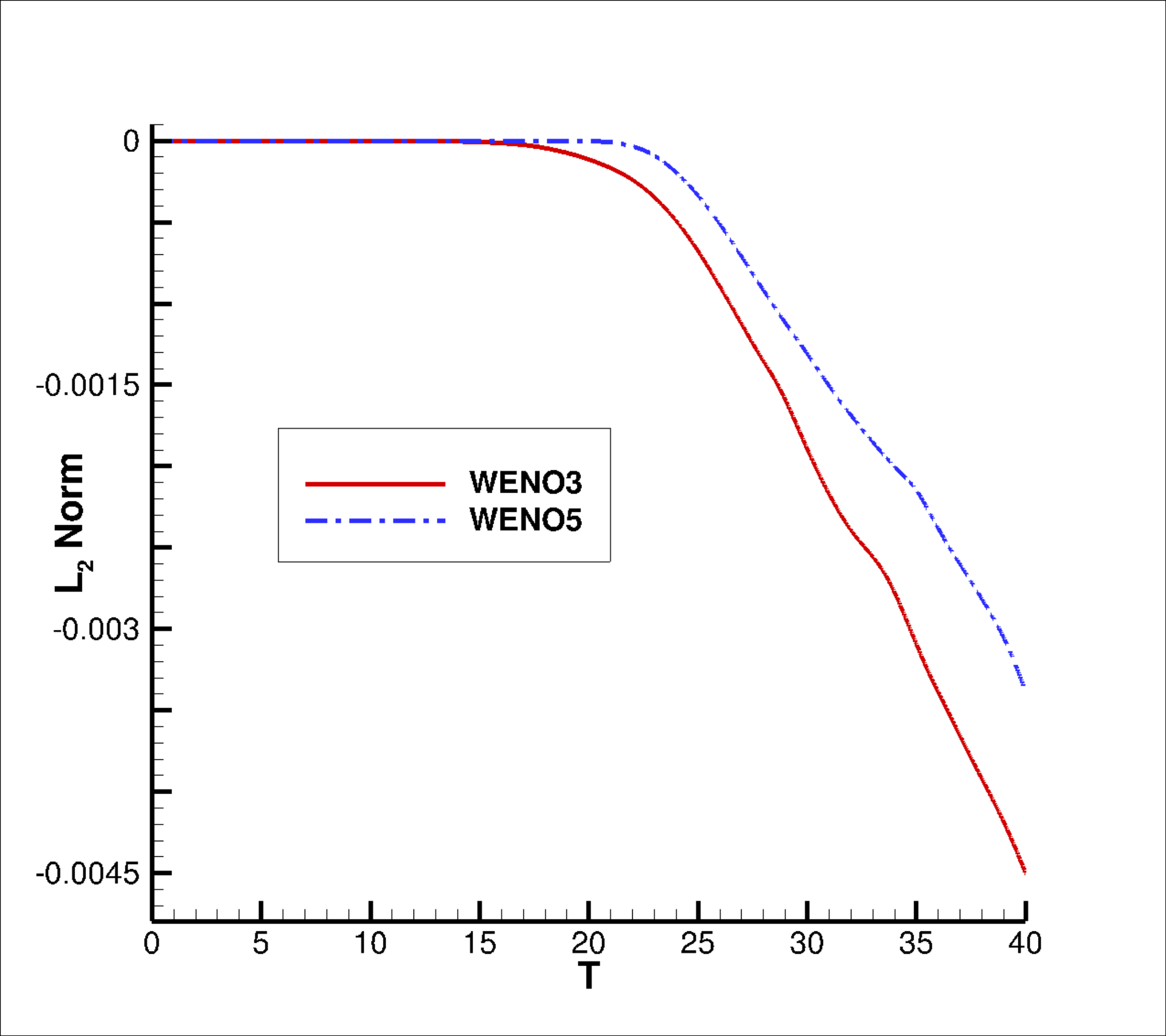}}
\subfigure[Two-stream instability II]{
\includegraphics[width=0.4\textwidth]{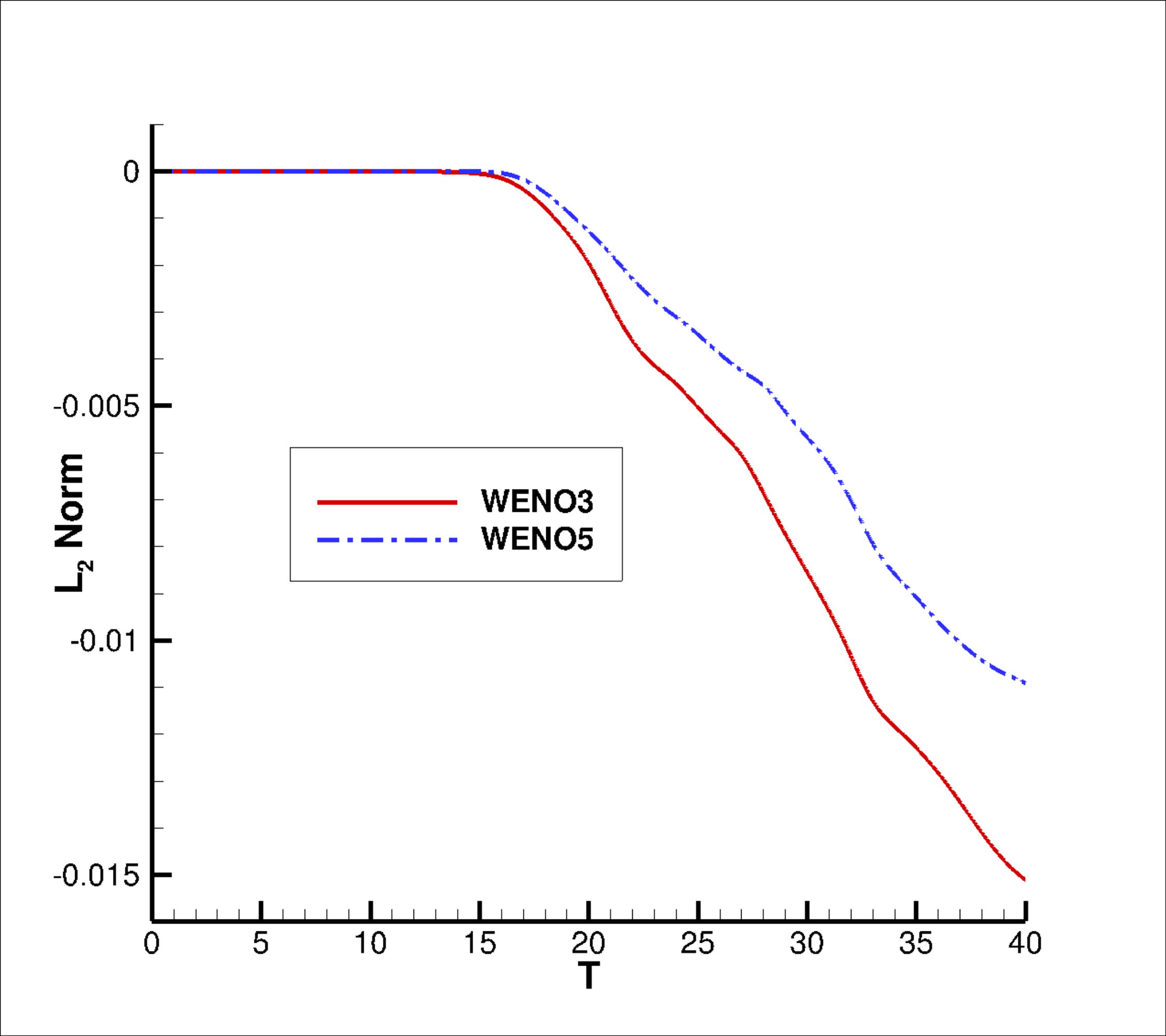}}
\subfigure[Bump-on-tail instability]{
\includegraphics[width=0.4\textwidth]{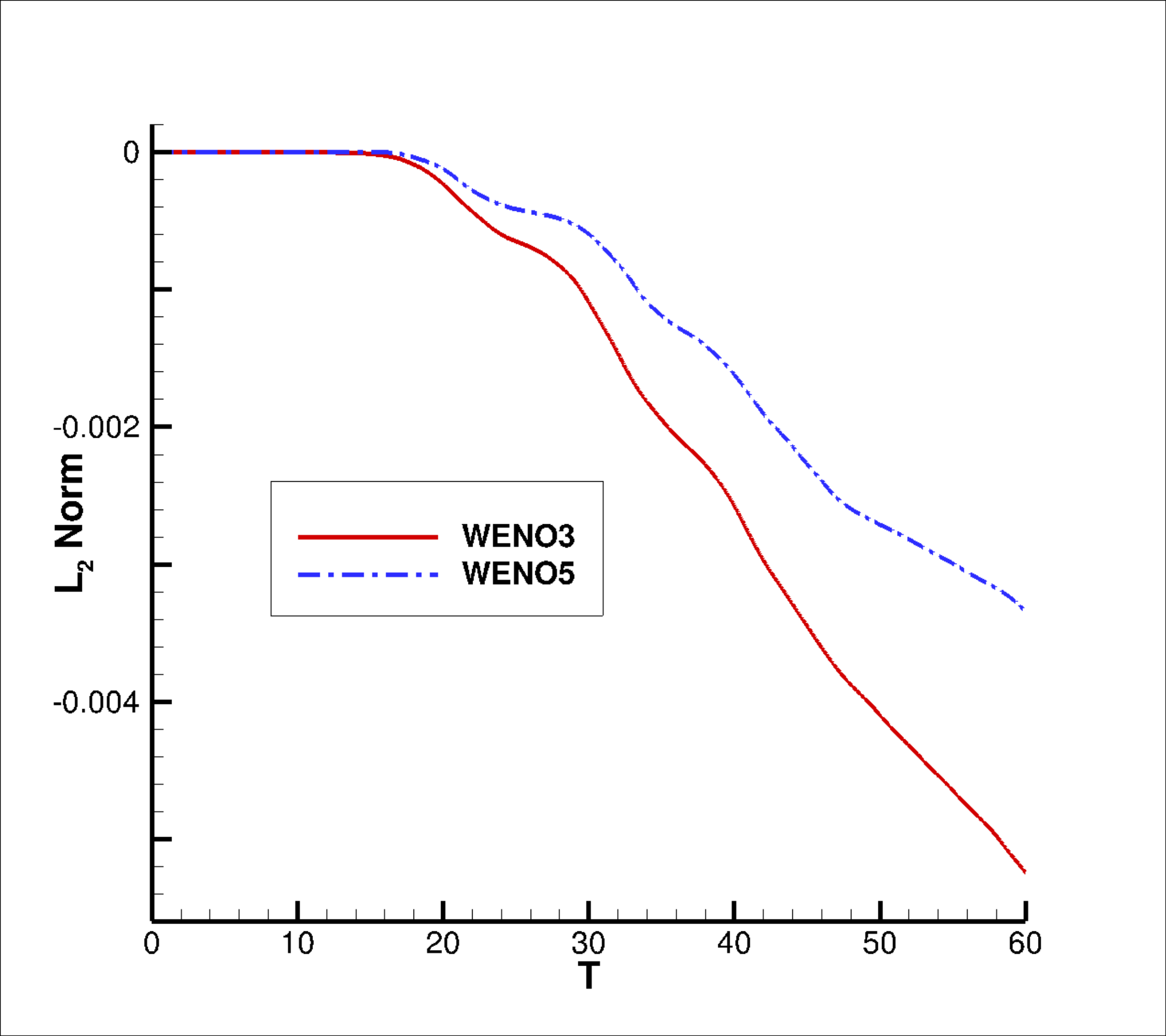}}
\caption{\em The time evolution of the relative deviation in $L_{2}$ norm. $256\times1024$ grid points. }
\label{Fig14}
\end{figure}

\begin{figure}
\centering
\subfigure[Strong Landau damping]{
\includegraphics[width=0.4\textwidth]{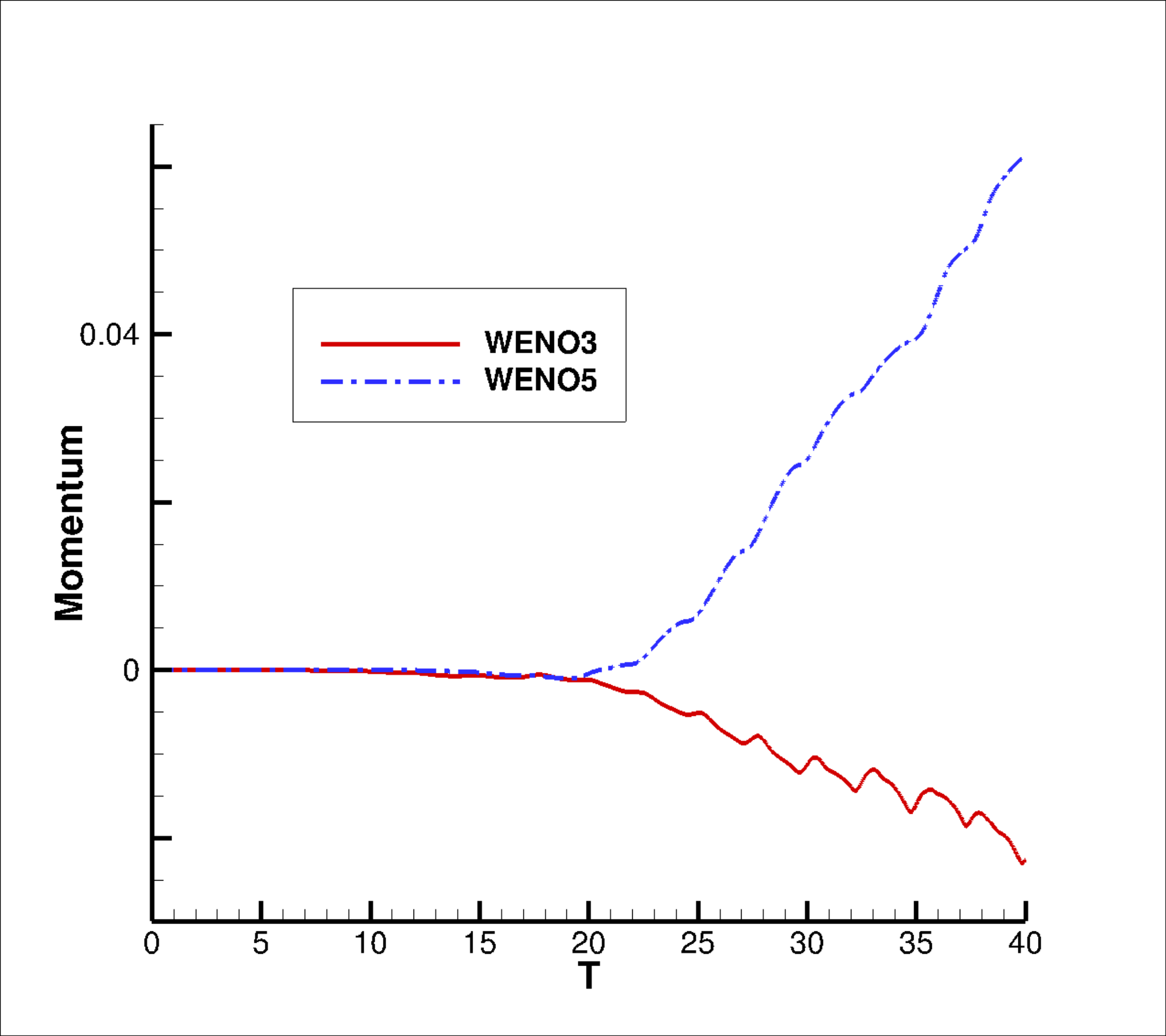}}
\subfigure[Two-stream instability I]{
\includegraphics[width=0.4\textwidth]{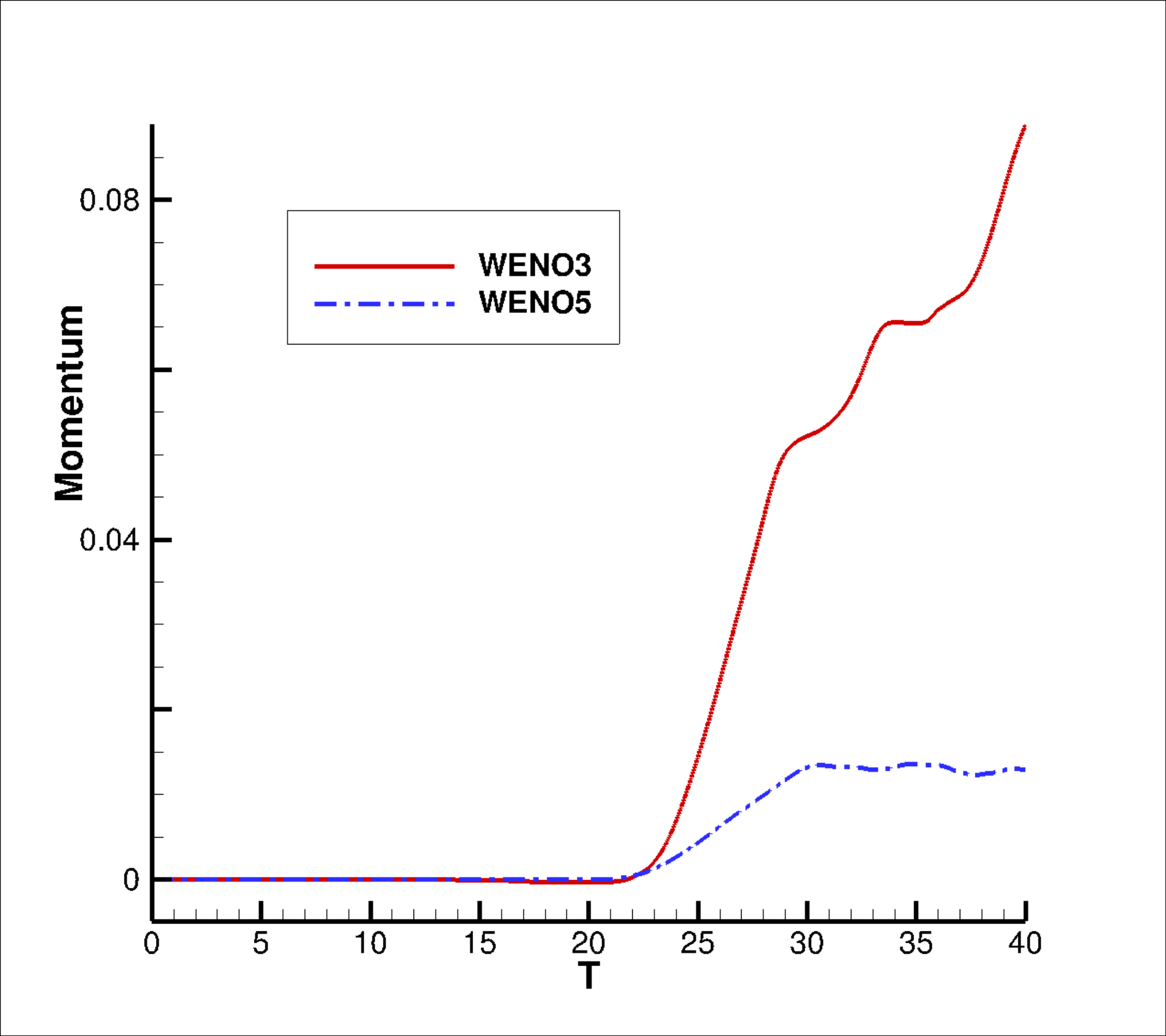}}
\subfigure[Two-stream instability II]{
\includegraphics[width=0.4\textwidth]{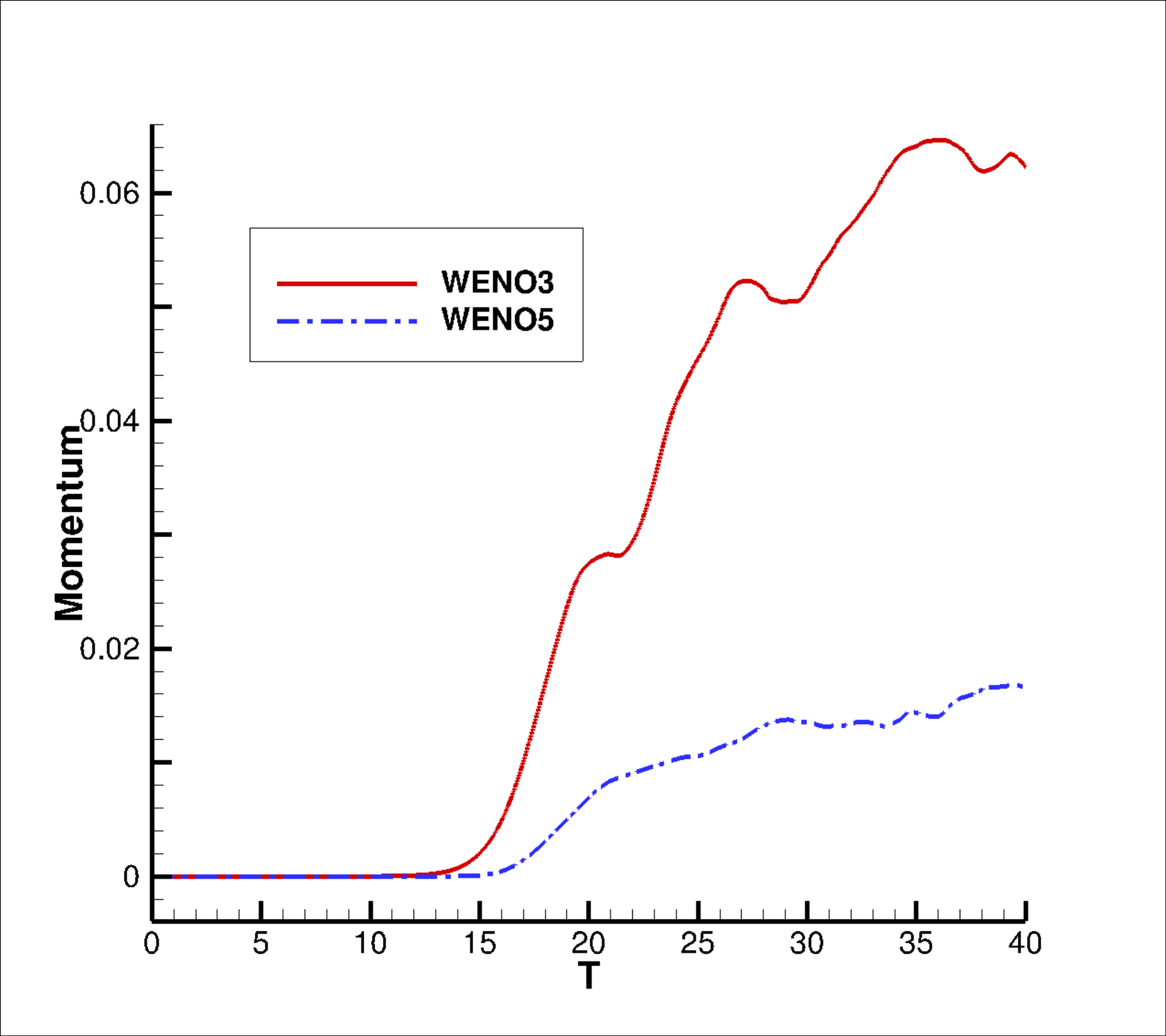}}
\subfigure[Bump-on-tail instability]{
\includegraphics[width=0.4\textwidth]{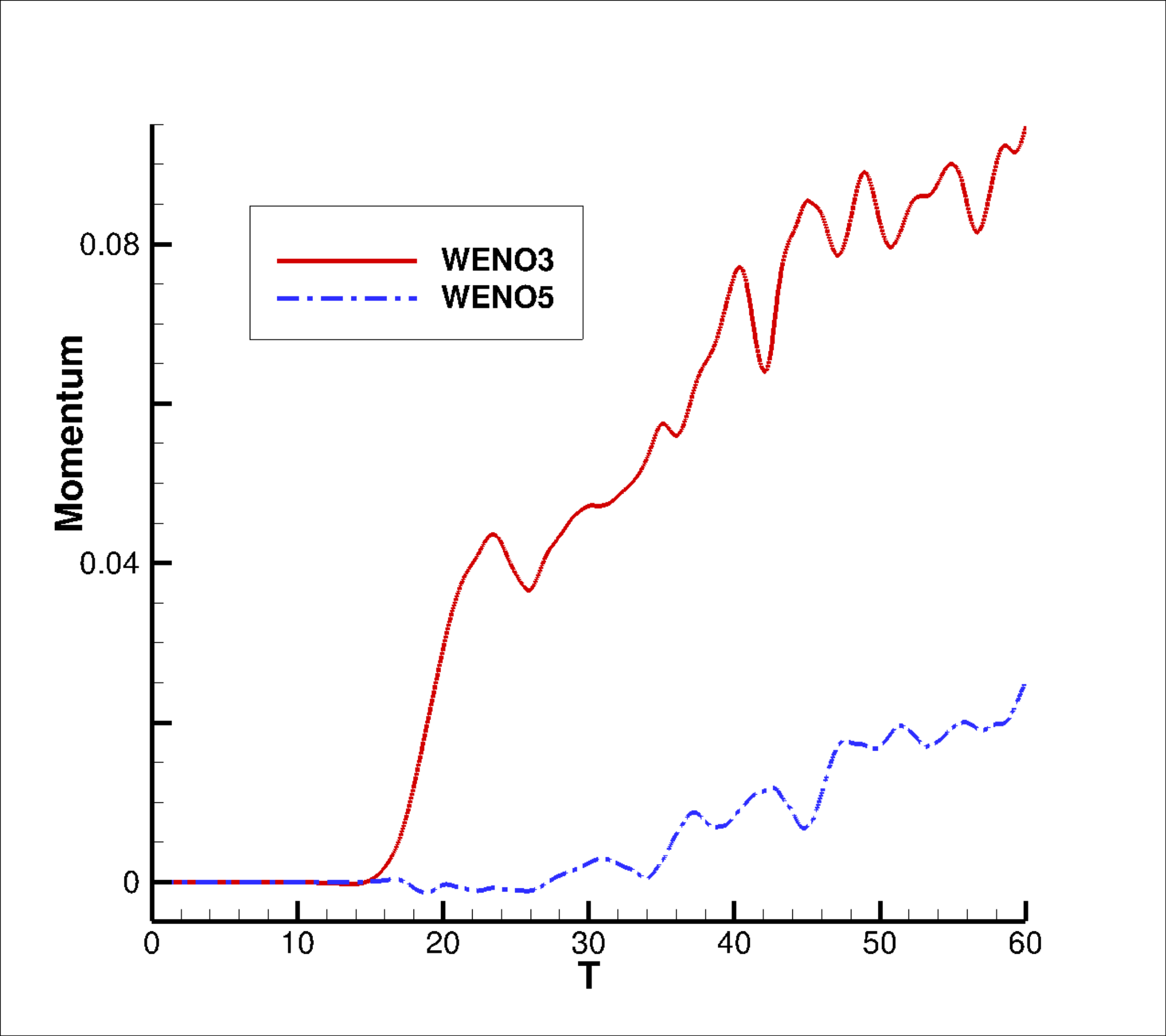}}
\caption{\em The time evolution of momentum. $256\times1024$ grid points. }
\label{Fig15}
\end{figure}

\subsection{Plasma theath}

\begin{figure}
\centering
\subfigure[WENO3]{
\includegraphics[width=0.8\textwidth]{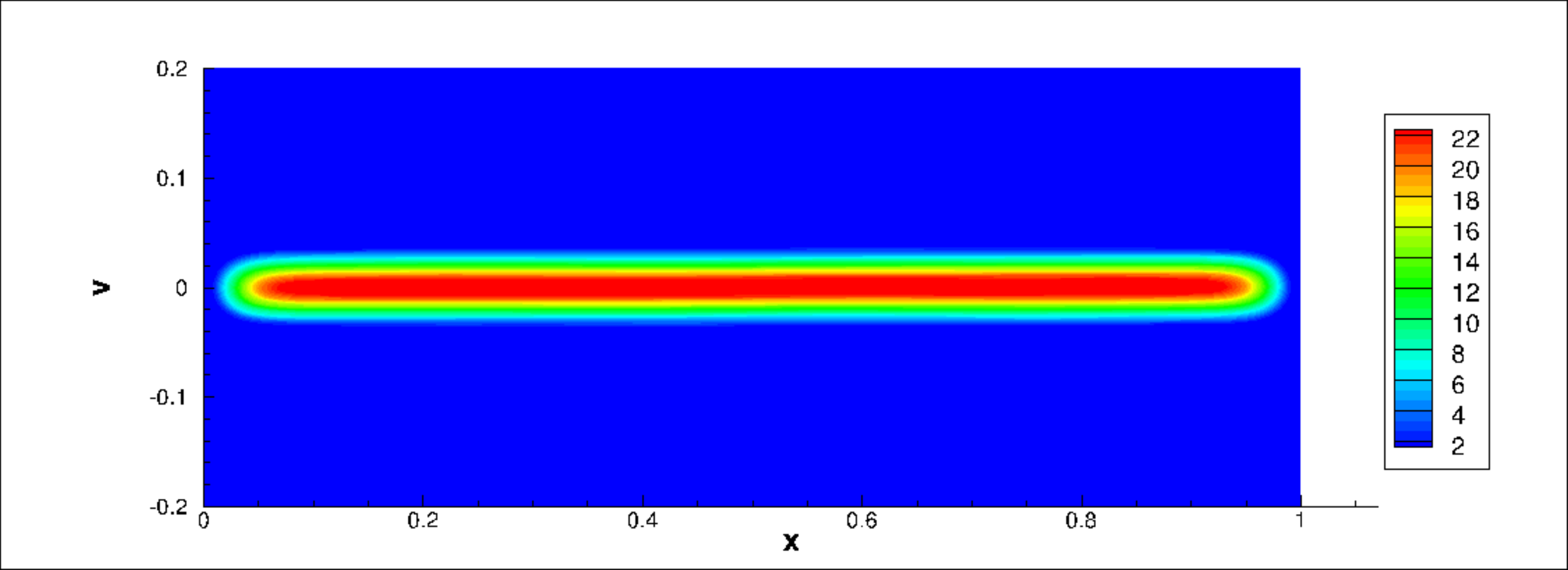}}
\subfigure[WENO5]{
\includegraphics[width=0.8\textwidth]{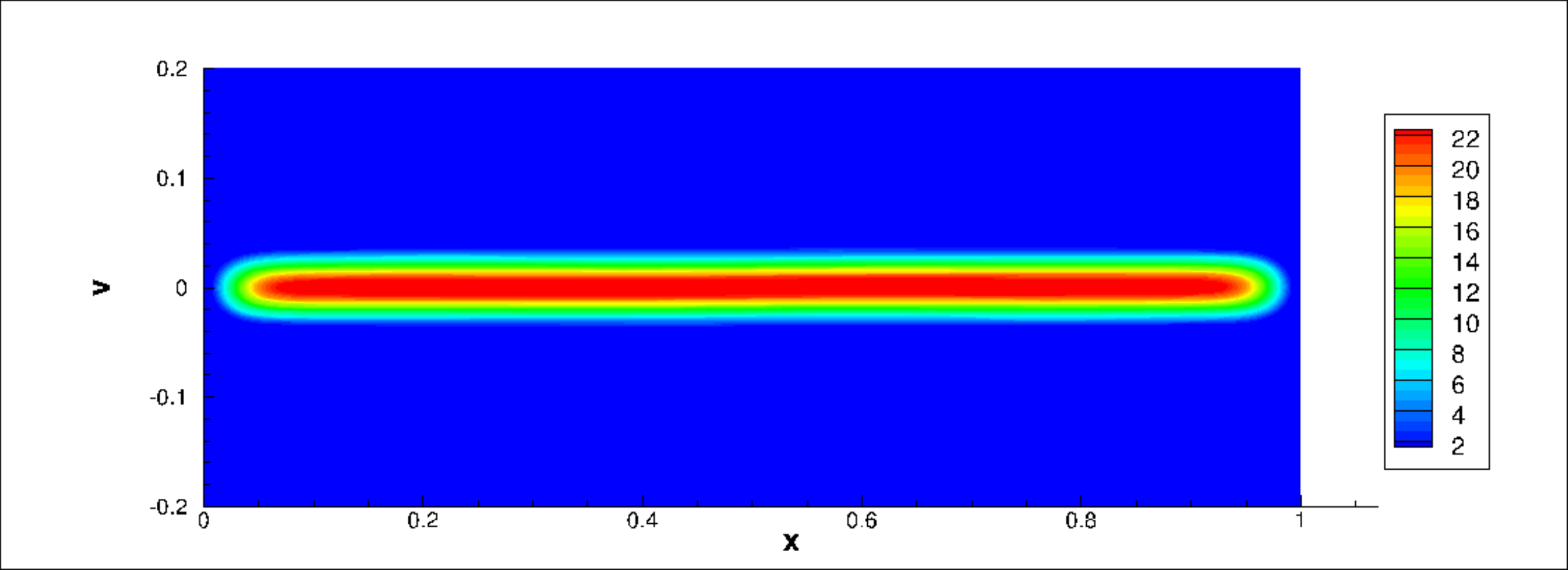}}
\caption{\em Plasma sheath. $256\times1024$ grid points. $T=140$}
\label{Fig16}
\end{figure}

Next, we consider the so-called plasma sheath problem, which discribes VP system 
with non-periodic boundary condition in x-direction. The initial condition is
\begin{align}
f(x,v,0)=\frac{1}{\sqrt{\alpha\pi}}\exp(-\frac{1}{\alpha}v^2) , \ \ \ x\in[0,L], \ \ \ v\in[-V_{c},V_{c}]
\end{align}
where $\alpha=0.0005526350206$, $L=1$ and $V_{c}=0.2$.
We assume that no electrons can enter the domain, meaning that we take zero-inflow boundary condition: 
\begin{align}
f(x= 0,t, v\geq 0) = f( x= L,t, v\leq 0) = 0.
\end{align}
The electric potential is fixed to a constant value at $x=0$ and $x=L$. Without loss of generality, 
we take $\phi(x=0,t)=\phi(x=L,t)=0$. The solutions at the final time $T=140$ are shown in Figure \ref{Fig16}, with $256\times1024$ grid points. Our results are consistent with that in \cite{seal2012discontinous} .

\section{Conclusion}

In this paper, we proposed a high order implicit Method of Lines Transpose (MOL$^{T}$) method for linear advection equations, which was further applied to the Vlasov-Poisson (VP) system for plasma simulations. The scheme coupled a WENO methodology and the diagonally implicit strong-stability-preserving  Runge-Kutta method in the MOL$^{T}$ framework, in order to attain high order accuracy in smooth region and avoid oscillation near the discontinuities with large time step evolution. Further, the scheme is able to efficiently deal with general boundary conditions, as a advantage over the popular semi-Lagrangian schemes. A simple high order positivity-preserving (PP)-limiter was developed and analyzed to ensure the desirable PP property for numerical solutions. The algorithm was extended to multi-dimensional problems and the VP system via dimensional splitting. A collection of numerical tests demonstrated good performance of the proposed scheme in resolving complex solutions' structures. Future work consists of extending the scheme to high-dimensional VP simulations.

\appendix
\section{Formulation of WENO integration and extrapolation}

When approximating $J^{L}_{i}$ with third order accuracy, the big stencil is $S(i)=\{x_{i-2},x_{i-1},x_{i},x_{i+1}\}$
and the small stencils are 
$S_{1}(i)=\{x_{i-2},x_{i-1},x_{i}\}$ and $S_{0}(i)=\{x_{i-1},x_{i},x_{i+1}\}$.
We denote $\nu=\alpha \Delta x$, then we have
\begin{align*}
c^{(0)}_{i-1}=& \frac{2+\nu-(2+3\nu+2\nu^2 )e^{-\nu}}{2\nu^2}, \\
c^{(0)}_{i}=& -\frac{2-\nu^2-(2+2\nu)e^{-\nu}}{\nu^2}, \\
c^{(0)}_{i}=& \frac{2-\nu-(2+\nu)e^{-\nu}}{2\nu^2}, \\
c^{(1)}_{i-2}=& \frac{2-\nu-(2+\nu)e^{-\nu}}{2\nu^2}, \\
c^{(1)}_{i-1}=& -\frac{2-2\nu-(2-\nu^2)e^{-\nu}}{\nu^2}, \\
c^{(1)}_{i}=& \frac{2-3\nu+2\nu^2-(2-\nu)e^{-\nu}}{2\nu^2}.
\end{align*}
And the linear weights are
\begin{align*}
& d_{0}=\frac{-( 6 - 6\nu + 2\nu^2) -(6-\nu^2)e^{-\nu} }{3\nu ( 2-\nu-(2+\nu)e^{-\nu} )},\\
& d_{1}=\frac{ 6-\nu^2-(6+6\nu+2\nu^2)e^{-\nu}}{3\nu (2-\nu -(2 + \nu)e^{-\nu}) }.
\end{align*}
In this case, the smoothness indicators are
\begin{align*}
& \beta_{0}=\frac{13}{12}(v_{i-1}-2v_{i}+v_{i+1})^2 + (v_{i-1}-v_{i})^2, \\ 
& \beta_{1}=\frac{13}{12}(v_{i-2}-2v_{i-1}+v_{i})^2 +(v_{i-1}-v_{i})^2.
\end{align*}

When approximating $J^{L}_{i}$ with fifth order accuracy, we choose the large stencil as 
$S(i)=\{x_{i-3},x_{i-2},x_{i-1},x_{i},x_{i+1}x_{i+2}\}$, and the three small stencils are
$S_{0}(i)=\{x_{i-1},x_{i},x_{i+1},x_{i+2}\}$, $S_{1}(i)=\{x_{i-2},x_{i-1},x_{i},x_{i+1}\}$, and $S_{2}(i)=\{x_{i-3},x_{i-2},x_{i-1},x_{i}\}$. Then we can obtain
\begin{align*}
c^{(0)}_{i-1} =& \frac{6+6\nu+2\nu^2-(6+12\nu+11\nu^2+6\nu^3)e^{-\nu}}{6\nu^3}, \\
c^{(0)}_{i} =& -\frac{6+4\nu-\nu^2-2\nu^3-(6+10\nu+6\nu^2)e^{-\nu}}{2\nu^3}, \\
c^{(0)}_{i+1}=& \frac{6+2\nu-2\nu^2-(6+8\nu+3\nu^2)e^{-\nu}}{2\nu^3},\\
c^{(0)}_{i+2}=& -\frac{6-\nu^2-(6+6\nu+2\nu^2)e^{-\nu}}{6\nu^3} u_{i+2},\\
c^{(1)}_{i-2} =& \frac{6-\nu^2-(6+6\nu+2\nu^2)e^{-\nu}}{6\nu^3},\\
c^{(1)}_{i-1} =&-\frac{6-2\nu-2\nu^2-(6+4\nu-\nu^2-2\nu^3)e^{-nu}}{2\nu^3}, \\
c^{(1)}_{i} =& \frac{6-4\nu-\nu^2+2\nu^3-(6+2\nu-2\nu^2)e^{-\nu}}{2\nu^3}, \\
c^{(1)}_{i+1}=& -\frac{6-6\nu+2\nu^2-(6-\nu^2)e^{-\nu}}{6\nu^3}, \\
c^{(2)}_{i-3} =& \frac{6-6\nu+2\nu^2-(6-\nu^2 )e^{-\nu}}{6\nu^3}, \\
c^{(2)}_{i-2} =& -\frac{6-8\nu+3\nu^2-(6-2\nu-2\nu^2)e^{-\nu}}{2\nu^3}, \\
c^{(2)}_{i-1} =& \frac{6-10\nu+6\nu^2-(6-4\nu-\nu^2+2\nu^3)e^{-\nu}}{2\nu^3},\\
c^{(2)}_{i} =& -\frac{6-12\nu+11\nu^2-6\nu^3-(6-6\nu+2\nu^2)e^{-\nu}}{6\nu^3}.
\end{align*}
And the linear weights are
\begin{align*}
& d_{0}=\frac{ 60-60\nu+15\nu^2+5\nu^3-3\nu^4-(60-15\nu^2+2\nu^4)e^{-\nu}}{10\nu^2(6-\nu^2-(6+6\nu+2\nu^2)e^{-\nu})}\\
& d_{2} =\frac{6-\nu^2-(6 + 6\nu + 2\nu^2 )e^{-\nu} }{3\nu ((2-\nu)-(2+\nu)e^{-\nu})}\\
& d_{1}=1-d_{0}-d_{2}
\end{align*}
The smoothness indicators are
\begin{align*}
& \beta_{0}= \frac{781}{120}(-v_{i-1}+3v_{i}-3v_{i+1}+v_{i+2})^2 
+ \frac{13}{48}(-3v_{i-1}+7v_{i}-5v_{i+1}-v_{i+2})^2 + (v_{i-1}-v_{i})^2 \\
& \beta_{1}= \frac{781}{120}(-v_{i-2}+3v_{i-1}-3v_{i}+v_{i+1})^2 
+ \frac{13}{48}(v_{i-2}-v_{i-1}-v_{i}+v_{i+1})^2 + (v_{i-1}-v_{i})^2\\
& \beta_{2}= \frac{781}{120}(-v_{i-3}+3v_{i-2}-3v_{i-1}+v_{i})^2 
+ \frac{13}{48}(v_{i-3}-5v_{i-2}+7v_{i-1}-3v_{i})^2 + (v_{i-1}-v_{i})^2 
\end{align*}

When doing WENO extrapolation on left boundary, the smoothness indicators are listed below:
\begin{align*}
 \beta_{1} =& \Delta x^2\\
 \beta_{2} =& ( v_{0} - v_{1} )^2\\
 \beta_{3} =& \frac{13}{12}( v_{0} - 2v_{1} + v_{2} )^2 + ( 2v_{0} - 3v_{1} +v_{2} )^2\\
 \beta_{4} =& \frac{781}{720} (v_{0}-3v_{1}+3v_{2} -v_{3} )^2 
+ \frac{13}{48}( 5v_{0}-13v_{1}+11v_{2}-3v_{3} )^2 
+ ( 3v_{0} - 6v_{1} + 4v_{2} - v_{3} )^2\\
 \beta_{5} =& \frac{1421461}{1310400} ( v_{0}-4v_{1}+ 6v_{2}- 4v_{3} + v_{4} )^2  
+ \frac{781}{720} ( -3v_{0}+11v_{1}-15v_{2}+9v_{3} - 2v_{4} )^2  \\
& + \frac{13}{7300800}( 3379v_{0}-10786v_{1} +12864v_{2}-6886v_{3} + 1429v_{4} )^2  \\
& + ( -4v_{0} + 10v_{1} - 10v_{2} + 5v_{3} - v_{4} )^2
\end{align*}

For $J^{R}_{i}$ and extrapolation on the right boundary, we can have the coefficients due to the symmetry.

\section{Formulation of RK(2,3) and RK(4,4)}

For $RK(2,3)$, the matrixes $\textbf{A}$ and $\textbf{b}$ are given as
\begin{align*}
& \textbf{A}=\begin{pmatrix}
\frac{1}{2}(1-\frac{1}{\sqrt{3}}) & 0 \\
\frac{1}{\sqrt{3}} & \frac{1}{2}(1-\frac{1}{\sqrt{3}}) \\
\end{pmatrix}, 
& \textbf{b}=\begin{pmatrix}
\frac{1}{2} \\
\frac{1}{2} \\
\end{pmatrix}
\end{align*}
When $c>0$, $RK(2,3)$ can be rewritten as
\begin{align*}
 u^{(1)}(x) = & I^L\left[u^n,\alpha \right](x) + A^{rk,1}e^{-\alpha(x-a)}, \\
 u^{(2)}(x) = & I^L\left[-\sqrt{3}u^n+(1+\sqrt{3})u^{(1)},\alpha \right](x) + A^{rk,2}e^{-\alpha(x-a)}\\
 u^{n+1}(x )= & (1+\sqrt{3})u^{n}(x) - \frac{3}{2}(1+\sqrt{3})u^{(1)}(x) +\frac{1}{2}(3+\sqrt{3})u^{(2)}(x). 
\end{align*}
where $\alpha = \frac{3+\sqrt{3}}{c\Delta t}$. And the corresponding time on each stage is $t^{1}_{rk}=t_{n}+\frac{1}{2}(1-\frac{1}{\sqrt{3}})\Delta t$ and $t^{2}_{rk}=t_{n}+\frac{1}{2}(1+\frac{1}{\sqrt{3}})\Delta t$. 

\bigskip

For $RK(4,4)$, the matrixes $\textbf{A}$ and $\textbf{b}$ are given as 
\begin{align*}
\text{A}& =\begin{pmatrix}
0.087475824368378 & 0 & 0 & 0 \\
0.306653000581791 & 0.106634669130071 & 0 & 0 \\
0.306653000581791 & 0.325811845343484 & 0.106634688637712 & 0 \\
0.306049667930486 & 0.220166571892301 & 0.220166585074543 & 0.087475807723977 \\
\end{pmatrix}\\
\text{b}^{T} & =\begin{pmatrix}
0.306092539007907 &
0.204522170534763 &
0.204522182780312 &
0.284863107677018 \\
\end{pmatrix}
\end{align*}
When $c>0$, we can wirte $RK(4,4)$ as
\begin{align*}
u^{(1)}(x)=& I^{L}\left[ u^{n},\alpha_{1}\right](x) +A^{rk,1}e^{-\alpha_{1}(x-a)},\\
u^{(2)}(x)=& I^{L}[-2.50557428633555u^{n} + 3.5055742863355546u^{(1)},\alpha_{2}] +A^{rk,2}e^{-\alpha_{2}(x-a)}\\
u^{(3)}(x)=& I^{L}[5.149963903410283 u^{n} - 7.205366503722802 u^{(1)} + 3.0554026003125183 u^{(2)}, \alpha_{3}] \\
& +A^{rk,3} e^{-\alpha_{3}(x-a)} \\
u^{(4)}(x)=& I^{L}[-7.958496756460987 u^{n} + 11.1375658967502 u^{(1)} - 4.243749854297317 u^{(2)}\\
&  + 2.0646807140081034 u^{(3)},\alpha_{4}] +A^{rk,4} e^{-\alpha_{4}(x-a)} \\
u^{n+1}=& 18.345647630787816 u^{n} - 25.673987908728908 u^{(1)} + 9.877479922993285 u^{(2)} \\
 & - 4.805618380017192 u^{(3)} + 3.256478734964997 u^{(4)}
\end{align*}
where 
\begin{align*}
& \alpha_{1}=1/(0.087475824368378 c\Delta t),\\
& \alpha_{2}=1/(0.106634669130071 c\Delta t),\\
& \alpha_{3}=1/(0.106634688637712 c\Delta t),\\
& \alpha_{4}=1/(0.087475807723977 c\Delta t).
\end{align*}
The corresponding time on each stage is 
\begin{align*}
& t^{1}_{rk}=t_{n}+0.087475824368378\Delta t,\\
& t^{2}_{rk}=t_{n}+0.413287669711862\Delta t, \\
& t^{3}_{rk}=t_{n}+0.739099534562987\Delta t, \\
& t^{4}_{rk}=t_{n}+0.833858632621307\Delta t.
\end{align*}

In the case $c<0$, we have the similar form. 

\bibliographystyle{abbrv}
\bibliography{ref}

\begin{thebibliography}{10}

\bibitem{Isaia_rkbc}
I.~Alonso-Mallo.
\newblock Runge-kutta methods without order reduction for linear initial
  boundary value problems.
\newblock {\em Numerische Mathematik}, 91(4):577--603.

\bibitem{begue1999two}
M.~Begue, A.~Ghizzo, and P.~Bertrand.
\newblock {Two-dimensional Vlasov simulation of Raman scattering and plasma
  beatwave acceleration on parallel computers}.
\newblock {\em J. Comput. Phys.}, 151(2):458--478, 1999.

\bibitem{besse2005semi}
N.~Besse, J.~Segr{\'e}, and E.~Sonnendr{\"u}cker.
\newblock Semi-lagrangian schemes for the two-dimensional vlasov-poisson system
  on unstructured meshes.
\newblock {\em Transport Theory and Statistical Physics}, 34(3-5):311--332,
  2005.

\bibitem{birdsall2005plasma}
C.~K. Birdsall and A.~B. Langdon.
\newblock {\em Plasma Physics Via Computer Simulaition}.
\newblock CRC Press, 2005.

\bibitem{causley2014higher}
M.~Causley and A.~Christlieb.
\newblock {Higher order A-stable schemes for the wave equation using a
  successive convolution approach}.
\newblock {\em SIAM Journal on Numerical Analysis}, 52(1):220--235, 2014.

\bibitem{causley}
M.~Causley, A.~Christlieb, B.~Ong, and L.~Van~Groningen.
\newblock Method of lines transpose: An implicit solution to the wave equation.
\newblock {\em Mathematics of Computation}, 83(290):2763--2786, 2014.

\bibitem{Causley2015method}
M.~F. Causley, H.~Cho, A.~J. Christlieb, and D.~C. Seal.
\newblock Method of lines transpose: High order l-stable
  $\backslash$mathcalo(n) schemes for parabolic equations using successive
  convolution.
\newblock {\em SIAM Journal on Numerical Analysis}, 54(3):1635--1652, 2016.

\bibitem{causley2013method}
M.~F. Causley, A.~Christlieb, Y.~G\"{u}\c{c}l\"{u}, and E.~Wolf.
\newblock Method of lines transpose: A fast implicit wave propagator.
\newblock {\em arXiv preprint arXiv:1306.6902}, 2013.

\bibitem{cheng1975integration}
C.~Cheng and G.~Knorr.
\newblock {The integration of the Vlasov equation in configuration space}.
\newblock {\em J. Comput. Phys.}, 22(3):330--351, 1976.

\bibitem{Cheng2015asymptotic}
Y.~Cheng, A.~J. Christlieb, W.~Guo, and B.~Ong.
\newblock {An asymptotic preserving Maxwell solver resulting in the Darwin
  limit of electrodynamics}.
\newblock {\em submitted}, 2015.

\bibitem{cheng2014ts}
Y.~Cheng, A.~J. Christlieb, and X.~Zhong.
\newblock Numerical study of the two-species {Vlasov-Amp\'{e}re} system:
  energy-conserving schemes and the current-driven ion-acoustic instability.
\newblock 2014.
\newblock submitted.

\bibitem{cheng2012study}
Y.~Cheng, I.~M. Gamba, and P.~J. Morrison.
\newblock Study of conservation and recurrence of runge--kutta discontinuous
  {G}alerkin schemes for vlasov--poisson systems.
\newblock {\em J. Sci. Comput.}, pages 1--31, 2012.

\bibitem{chou2006high}
C.-S. Chou and C.-W. Shu.
\newblock High order residual distribution conservative finite difference weno
  schemes for steady state problems on non-smooth meshes.
\newblock {\em J. Comput. Phys.}, 214(2):698--724, 2006.

\bibitem{chou2007high}
C.-S. Chou and C.-W. Shu.
\newblock High order residual distribution conservative finite difference weno
  schemes for convection--diffusion steady state problems on non-smooth meshes.
\newblock {\em J. Comput. Phys.}, 224(2):992--1020, 2007.

\bibitem{christlieb2014high}
A.~Christlieb, W.~Guo, M.~Morton, and J.-M. Qiu.
\newblock {A high order time splitting method based on integral deferred
  correction for semi-Lagrangian Vlasov simulations}.
\newblock {\em J. Comput. Phys.}, 267:7--27, 2014.

\bibitem{crouseilles2010conservative}
N.~Crouseilles, M.~Mehrenberger, and E.~Sonnendr{\"u}cker.
\newblock Conservative semi-lagrangian schemes for vlasov equations.
\newblock {\em Journal of Computational Physics}, 229(6):1927--1953, 2010.

\bibitem{crouseilles2009forward}
N.~Crouseilles, T.~Respaud, and E.~Sonnendr{\"u}cker.
\newblock A forward semi-lagrangian method for the numerical solution of the
  vlasov equation.
\newblock {\em Computer Physics Communications}, 180(10):1730--1745, 2009.

\bibitem{FilbetS}
F.~Filbet and E.~Sonnendr{\"u}cker.
\newblock {Comparison of Eulerian Vlasov solvers}.
\newblock {\em Comput. Phys. Commun.}, 150(3):247--266, 2003.

\bibitem{filbet2001conservative}
F.~Filbet, E.~Sonnendr{\"u}cker, and P.~Bertrand.
\newblock {Conservative numerical schemes for the Vlasov equation}.
\newblock {\em J. Comput. Phys.}, 172(1):166--187, 2001.

\bibitem{forest1990fourth}
E.~Forest and R.~Ruth.
\newblock Fourth-order symplectic integration.
\newblock {\em Physica D: Nonlinear Phenomena}, 43(1):105--117, 1990.

\bibitem{gottlieb2001strong}
S.~Gottlieb, C.-W. Shu, and E.~Tadmor.
\newblock Strong stability-preserving high-order time discretization methods.
\newblock {\em SIAM review}, 43(1):89--112, 2001.

\bibitem{guo2013hybrid}
W.~Guo and J.-M. Qiu.
\newblock {Hybrid semi-Lagrangian finite element-finite difference methods for
  the Vlasov equation}.
\newblock {\em J. Comput. Phys.}, 234:108--132, 2013.

\bibitem{hockney2010computer}
R.~W. Hockney and J.~W. Eastwood.
\newblock {\em Computer simulation using particles}.
\newblock CRC Press, 2010.

\bibitem{jiang1996efficient}
G.-s. Jiang, C.-W. Shu, et~al.
\newblock Efficient implementation of weighted eno schemes.
\newblock In {\em J. Comput. Phys.}, 1996.

\bibitem{ketcheson2009optimal}
D.~I. Ketcheson, C.~B. Macdonald, and S.~Gottlieb.
\newblock Optimal implicit strong stability preserving runge--kutta methods.
\newblock {\em Applied Numerical Mathematics}, 59(2):373--392, 2009.

\bibitem{liang2014parametrized}
C.~Liang and Z.~Xu.
\newblock Parametrized maximum principle preserving flux limiters for high
  order schemes solving multi-dimensional scalar hyperbolic conservation laws.
\newblock {\em J. Sci. Comput.}, 58(1):41--60, 2014.

\bibitem{liu1994weighted}
X.-D. Liu, S.~Osher, and T.~Chan.
\newblock Weighted essentially non-oscillatory schemes.
\newblock {\em J. Comput. Phys.}, 115(1):200--212, 1994.

\bibitem{liu2009positivity}
Y.-y. Liu, C.-w. Shu, and M.-p. Zhang.
\newblock On the positivity of linear weights in weno approximations.
\newblock {\em Acta Mathematicae Applicatae Sinica, English Series},
  25(3):503--538, 2009.

\bibitem{parker1997convected}
G.~J. Parker and W.~N.~G. Hitchon.
\newblock Convected scheme simulations of the electron distribution function in
  a positive column plasma.
\newblock {\em Japanese journal of applied physics}, 36(7S):4799, 1997.

\bibitem{qiu2011positivity}
J.~Qiu and C.~Shu.
\newblock {Positivity preserving semi-Lagrangian discontinuous {G}alerkin
  formulation: Theoretical analysis and application to the Vlasov-Poisson
  system}.
\newblock {\em J. Comput. Phys.}, 230(23):8386--8409, 2011.

\bibitem{Qiu_Christlieb}
J.-M. Qiu and A.~Christlieb.
\newblock {A Conservative high order semi-Lagrangian WENO method for the Vlasov
  equation}.
\newblock {\em J. Comput. Phys.}, 229:1130--1149, 2010.

\bibitem{Qiu_Shu2}
J.-M. Qiu and C.-W. Shu.
\newblock {Conservative semi-Lagrangian finite difference WENO formulations
  with applications to the Vlasov equation}.
\newblock {\em Commun. Comput. Phys.}, 10(4):979--1000, 2011.

\bibitem{Qiu20118386}
J.-M. Qiu and C.-W. Shu.
\newblock {Positivity preserving semi-Lagrangian discontinuous {G}alerkin
  formulation: Theoretical analysis and application to the Vlasov-Poisson
  system}.
\newblock {\em J. Comput. Phys.}, 230(23):8386 -- 8409, 2011.

\bibitem{rossmanith2011positivity}
J.~Rossmanith and D.~Seal.
\newblock {A positivity-preserving high-order semi-Lagrangian discontinuous
  {G}alerkin scheme for the Vlasov-Poisson equations}.
\newblock {\em J. Comput. Phys.}, 230(16):6203--6232, 2011.

\bibitem{salazar2000theoretical}
A.~J. Salazar, M.~Raydan, and A.~Campo.
\newblock {Theoretical analysis of the Exponential Transversal Method of Lines
  for the diffusion equation}.
\newblock {\em Numerical Methods for Partial Differential Equations},
  16(1):30--41, 2000.

\bibitem{schemann1998adaptive}
M.~Schemann and F.~Bornemann.
\newblock {An adaptive Rothe method for the wave equation}.
\newblock {\em Computing and Visualization in Science}, 1(3):137--144, 1998.

\bibitem{seal2012discontinous}
D.~C. Seal.
\newblock {\em Discontinous Galerkin methods for Vlasov models of plasma}.
\newblock PhD thesis, UNIVERSITY OF WISCONSIN--MADISON, 2012.

\bibitem{shu2009high}
C.-W. Shu.
\newblock High order weighted essentially nonoscillatory schemes for convection
  dominated problems.
\newblock {\em SIAM review}, 51(1):82--126, 2009.

\bibitem{sonnendrucker1999semi}
E.~Sonnendr{\"u}cker, J.~Roche, P.~Bertrand, and A.~Ghizzo.
\newblock The semi-lagrangian method for the numerical resolution of the vlasov
  equation.
\newblock {\em Journal of computational physics}, 149(2):201--220, 1999.

\bibitem{tan2012efficient}
S.~Tan, C.~Wang, C.-W. Shu, and J.~Ning.
\newblock Efficient implementation of high order inverse lax--wendroff boundary
  treatment for conservation laws.
\newblock {\em Journal of Computational Physics}, 231(6):2510--2527, 2012.

\bibitem{verboncoeur2005particle}
J.~P. Verboncoeur.
\newblock Particle simulation of plasmas: review and advances.
\newblock {\em Plasma Phys. Contr. F.}, 47(5A):A231–--A260, 2005.

\bibitem{xiong2013parametrized}
T.~Xiong, J.-M. Qiu, and Z.~Xu.
\newblock A parametrized maximum principle preserving flux limiter for finite
  difference rk-weno schemes with applications in incompressible flows.
\newblock {\em Journal of Computational Physics}, 252:310--331, 2013.

\bibitem{xiong2014parametrized}
T.~Xiong, J.-M. Qiu, and Z.~Xu.
\newblock Parametrized positivity preserving flux limiters for the high order
  finite difference weno scheme solving compressible euler equations.
\newblock {\em arXiv preprint arXiv:1403.0594}, 2014.

\bibitem{xiong2014high}
T.~Xiong, J.-M. Qiu, Z.~Xu, and A.~Christlieb.
\newblock High order maximum principle preserving semi-lagrangian finite
  difference weno schemes for the vlasov equation.
\newblock {\em Journal of Computational Physics}, 273:618--639, 2014.

\bibitem{xu2014parametrized}
Z.~Xu.
\newblock Parametrized maximum principle preserving flux limiters for high
  order schemes solving hyperbolic conservation laws: one-dimensional scalar
  problem.
\newblock {\em Math. Comput.}, 83(289):2213--2238, 2014.

\bibitem{zhang2010maximum}
X.~Zhang and C.-W. Shu.
\newblock On maximum-principle-satisfying high order schemes for scalar
  conservation laws.
\newblock {\em Journal of Computational Physics}, 229(9):3091--3120, 2010.

\bibitem{zhang2010positivity}
X.~Zhang and C.-W. Shu.
\newblock On positivity-preserving high order discontinuous galerkin schemes
  for compressible euler equations on rectangular meshes.
\newblock {\em Journal of Computational Physics}, 229(23):8918--8934, 2010.

\bibitem{zhang2011maximum}
X.~Zhang and C.-W. Shu.
\newblock Maximum-principle-satisfying and positivity-preserving high-order
  schemes for conservation laws: survey and new developments.
\newblock In {\em Proceedings of the Royal Society of London A: Mathematical,
  Physical and Engineering Sciences}, page rspa20110153. The Royal Society,
  2011.

\bibitem{Zhang_JSC_2010}
X.~Zhang, Y.~Xia, and C.-W. Shu.
\newblock Maximum-principle-satisfying and positivity-preserving high order
  discontinuous {Galerkin} schemes for conservation laws on triangular meshes.
\newblock {\em J. Sci. Comput.}

\bibitem{zhou2001numerical}
T.~Zhou, Y.~Guo, and C.-W. Shu.
\newblock {Numerical study on Landau damping}.
\newblock {\em Phys. D}, 157(4):322--333, 2001.

\end{thebibliography}

\end{document}